 \documentclass{article}
 \usepackage{graphicx} 
 \newcommand{\matr}[4]{
\left[\begin{array}{cc}
#1&#2\\
#3&#4
\end{array}
\right]}

\def\J#1#2#3{ \left\{ #1,#2,#3 \right\} }

 \newcommand{\kbta}{K\overline{\otimes}A}

\def\cl{\hbox{$\mathcal L$}}
  \def\ad#1{ \hbox{ad}\, ({#1})}

 \newcommand{\tp}[3]{\{#1#2#3\}}

   \newcommand{\tpc}[3]{\{#1,#2,#3\}}

\newcommand{\csa}{$C^*$-algebra}

\newcommand{\tr}{\mbox{tr}\,}

\newcommand{\ip}[2]{\mbox{$(#1|#2)$}}

\newtheorem{lemma}{Lemma}[section]
\newtheorem{definition}[lemma]{Definition}
\newtheorem{theorem}[lemma]{Theorem}
\newcommand{\pf}{\noindent{\it Proof}.}
\newtheorem{proposition}[lemma]{Proposition}
\newtheorem{remark}[lemma]{Remark}
\newtheorem{corollary}[lemma]{Corollary}
\newtheorem{problem}[lemma]{Problem}

\newcommand{\jbst}{$JB^*$-triple}
\newcommand{\jcst}{$JC^*$-triple}
\newcommand{\jwst}{$JW^*$-triple}

\title{DERIVATIONS AND PROJECTIONS ON JORDAN TRIPLES\\
An introduction to nonassociative algebra, continuous cohomology, and quantum functional analysis }
\author{Bernard Russo}

\newcommand{\RR}{{\bf R}}

\newcommand{\CC}{{\bf C}}

\newcommand{\mn}{\hbox{$M_n(\IC)$}}
\newcommand{\IC}{\CC}
\newcommand{\e}[2]{e_{{#1}{#2}}}
\newcommand{\ee}[1]{e_{#1}}

\def\sqr#1#2{{\,\vcenter{\vbox{\hrule height.#2pt\hbox{\vrule width.#2pt
height#1pt \kern#1pt\vrule width.#2pt}\hrule height.#2pt}}\,}}
\def\bo{\sqr44\,}

\def\tp#1{\{#1\}}
\def\tpp#1#2#3{\{{#1}{#2}{#3}\}}

\def\tpc#1#2#3{\{#1\tp{#2}#3\}}

\begin{document}

\maketitle

\bigskip

This paper is an elaborated version of the material presented by the author in a  three hour minicourse at {\it V International Course of Mathematical Analysis in Andalusia}, at Almeria, Spain September 12-16, 2011. The author wishes to thank the scientific committee for the opportunity to present the course  and to the organizing committee for their hospitality.  The author also personally thanks Antonio Peralta  for his  collegiality and encouragement.\smallskip

The minicourse on which this paper is based had its genesis in a series of talks the author had given to undergraduates at Fullerton College in California.  I thank my former student Dana Clahane for his initiative in running the remarkable undergraduate research program at Fullerton College of which the seminar series is a part.  With their knowledge  only of  the product rule for differentiation as a starting point, these enthusiastic students were introduced to some aspects of the esoteric subject of non associative algebra, including triple systems as well as algebras.  Slides of these talks and of the minicourse lectures, as well as other related material,  can be found at the author's website (www.math.uci.edu/$\sim$brusso). Conversely, these undergraduate talks were motivated by the author's past and recent joint works on derivations of Jordan triples (\cite{HoMarPerRus02},\cite{HoPerRus11},\cite{PerRus10}), which are among the many results discussed here. \smallskip

Part I ({\it Derivations}) is devoted to an exposition of the properties of derivations on various algebras and triple systems in finite and infinite dimensions, the primary questions addressed being whether the derivation is automatically continuous and to what extent it is an inner derivation.  \smallskip

Part II ({\it Cohomology}) discusses cohomology theory of algebras and triple systems, in both finite and infinite dimensions.  Although the cohomology of associative and Lie algebras is substantially developed, in both finite and infinite dimensions (\cite{CarEil56},\cite{Fuks86},\cite{Knapp88}; see \cite{HochschildMR} for a review of \cite{CarEil56}), the same could not be said for Jordan algebras.  Moreover, the cohomology of triple systems has a rather sparse literature which is essentially non-existent in infinite dimensions.  Thus, one of the goals of this paper is to encourage the study of  continuous cohomology of some Banach triple systems. Occasionally, an idea for a research project is mentioned (especially in  subsection~\ref{7.2}). 
Readers are invited to contact the author, at brusso@uci.edu,  if they share this interest. \smallskip

 Part III ({\it Quantum Functional Analysis}) begins  with the subject of contractive projections, which  plays an important role in the structure theory of Jordan triples. The remainder of Part III discusses three topics, two very recent, which involve the interplay between Jordan theory and operator space theory (quantum functional analysis). The first one, a joint work of the author \cite{NeaRus03TAMS}, discusses the structure theory of contractively complemented Hilbertian operator spaces, and is instrumental to the third topic, which is concerned with some recent work on enveloping TROs and K-theory for JB*-triples \cite{BunFeeTim11},\cite{BunTim11},\cite{BohWer11}.  The second topic presents some very recent joint work by the author concerning quantum operator algebras which leans very heavily on contractive projection theory  \cite{NeaRus12}.\smallskip

Parts I and II, and the first section of Part III of this paper were covered in the three lectures making up the minicourse.  The rest of Part III of this paper, although not covered in the lectures, is included here as significant application of contractive projections and because there is some implicit reference to it in part I.\smallskip

A few proofs have been included but for the most part, theorems and other results are just stated and a reference for the proof and in some cases for a definition,  is provided.   Brief biographical information on some mathematicians, obtained primarily from Wikipedia,  is provided in footnotes.\smallskip

As part of the undergraduate lectures mentioned above, a set of exercises was developed, proceeding from elementary verifications to nontrivial theorems in the literature, with appropriate references. These occur at the end of sections 1 and 2.  \smallskip

Since this paper is meant to introduce a large number of topics, many of which are purely algebraical in nature, and others which the author has not fully studied yet, he has taken the liberty occasionally to  rely on Mathematical Reviews, as well as authors' introductions  for some of the summaries.

\tableofcontents

\part{Derivations}
\section{Algebras}\label{chpt1}

Let  ${\mathcal C}$ denote the algebra of continuous complex valued functions on a locally compact Hausdorff space.
 
\begin{definition}
A \underline{derivation} on $\mathcal C$
  is a linear mapping $\delta:{\mathcal C}\rightarrow {\mathcal C}$  satisfying the ``product''  rule:
$$\delta(fg)=\delta(f)g+f\delta(g).$$
 \end{definition}

 \begin{theorem}\label{thm:1112111}
There are no (non-zero) derivations on $\mathcal C$.  In other words, 
every derivation of $\mathcal C$ is identically zero.
 \end{theorem}

 The proof of this theorem can be broken up into two parts of independent interest.  The first theorem, due to Singer\footnote{{\bf Isadore Singer} (b. 1924).  Isadore Manuel Singer is an Institute Professor in the Department of Mathematics at the Massachusetts Institute of Technology. He is noted for his work with Michael Atiyah in 1962, which paved the way for new interactions between pure mathematics and theoretical physics} and Wermer\footnote{{\bf John Wermer} (b. 1925), retired from Brown University}, is valid for commutative semisimple Banach algebras. The second one is due to Sakai\footnote{{\bf Soichiro Sakai} (b. 1928), retired from University of Pennsylvania and returned to Japan}
 and is valid for all $C^*$-algebras.
   

\begin{theorem} [Singer and Wermer 1955 \cite{SinWer55}]
   Every \underline{continuous} derivation on $\mathcal C$ is zero.
  \end{theorem}
  
  \begin{theorem}[Sakai 1960 \cite{Sakai60}]
     Every derivation on $\mathcal C$ is continuous.
   \end{theorem}
 
 Singer and Wermer actually proved  that every continuous derivation on a commutative Banach algebra maps the algebra into the radical and they conjectured that the continuity assumption was unnecessary and could be removed.   Johnson showed this for commutative semisimple Banach algebras \cite{JohnsonAJM69} and  Thomas gave an affirmative answer to the Singer-Wermer conjecture for all commutative Banach algebras \cite{Thomas88}.  There are various noncommutative versions of the Singer-Wermer theorem for Banach algebras, for example \cite{MatMur91,MatRun92,BreMat95}. This is the first ripe area for exploration in the non associative contexts being considered in this paper.

\subsection{Derivations on finite dimensional algebras}\label{sub:1.1}
\subsubsection{Matrix multiplication---Associative algebras}

Let $M_n(\CC)$ denote the algebra of all $n$ by $n$ complex matrices, or more generally any finite dimensional semisimple  \underline{associative} algebra.

\begin{definition} 
 A \underline{derivation} on $M_n(\CC)$ with  respect to matrix multiplication is a linear mapping $\delta$ which satisfies the product rule
 
 $$\delta(AB)=\delta(A)B+A\delta(B).$$
\end{definition}

\begin{proposition} \label{prop:0202121}Fix a matrix  $A$ in $M_n(\CC)$ and define $$\delta_A(X)=AX-XA.$$ Then $\delta_A$ is a derivation with respect to matrix multiplication.
\end{proposition}

A more general form of the following theorem is due to Hochschild\footnote{
 {\bf Gerhard Hochschild} (1915--2010)
 Gerhard Paul Hochschild  was an American mathematician who worked on Lie groups, algebraic groups, homological algebra and algebraic number theory} \cite{Hochschild42}.   For an interesting summary of Hochschild's career  see \cite{SanMos11} and \cite{Santos11}. Theorem~\ref {thm:1112112}  is probably due to either Wedderburn\footnote{{\bf Joseph Henry Maclagan Wedderburn} (1882--1948).
 Scottish mathematician, who taught at Princeton University for most of his career. A significant algebraist, he proved that a finite division algebra is a field, and part of the Artin--Wedderburn theorem on simple algebras. He also worked on group theory and matrix algebra}
 or Noether\footnote{ {\bf Amalie Emmy Noether} (1882--1935)
  Amalie Emmy Noether was an influential German mathematician known for her groundbreaking contributions to abstract algebra and theoretical physics. Described as the most important woman in the history of mathematics, she revolutionized the theories of rings, fields, and algebras. In physics, Noether's theorem explains the fundamental connection between symmetry and conservation laws}.
The proof we present here is due to Jacobson \cite{Jacobson37}.

\begin{theorem} \label{thm:1112112} Every derivation on $M_n(\CC)$ with respect to matrix multiplication is of the form $\delta_A$ for some $A$ in $M_n(\CC)$.
 \end{theorem} 
 \pf\
 If $\delta$ is a derivation, consider the two representations of $M_n(\CC)$
\[
z\mapsto \matr{z}{0}{0}{z}\hbox{ and }z\mapsto  \matr{z}{0}{\delta(z)}{z}.
\]
The first is  a direct sum of two copies of the identity representation; but so is the second, since 
\[
 \matr{0}{0}{\delta(z)}{z}\hbox{ is equivalent to } \matr{0}{0}{0}{z}.
\]
So
$
\matr{z}{0}{\delta(z)}{0}\matr{a}{b}{c}{d}=\matr{a}{b}{c}{d}\matr{z}{0}{0}{z}.
$ Thus  $az=za, \ bz=zb$ and
\[
\delta(z)a=cz-zc\hbox{ and }\delta(z)b=dz-zd.
\]
Since $a$ and $b$ are multiples of $I$ they can't both be zero.

 \subsubsection{Bracket multiplication---Lie algebras}

 We next consider a multiplication on matrices (or associative algebras) which leads to the theory of Lie algebras\footnote{{\bf Sophus Lie}  (1842--1899)
Marius Sophus Lie was a Norwegian mathematician. He largely created the theory of continuous symmetry, and applied it to the study of geometry and differential equations}.
\begin{definition}
 A \underline{derivation} on $M_n(\CC)$ with respect to bracket multiplication $$[X,Y]=XY-YX$$
  is a linear mapping  $\delta$ which satisfies the product rule
 $$\delta([A,B])=[\delta(A),B]+[A,\delta(B)].$$
\end{definition}

 \begin{proposition}\label{prop:0202122}
Fix a matrix $A$ in $M_n(\CC)$ and define $$\delta_A(X)=[A,X]=AX-XA.$$ Then $\delta_A$ is a derivation with respect to bracket multiplication.
 \end{proposition}
 
 An algebra $\cl$ with multiplication $(x,y)\mapsto [x,y]$ is a Lie algebra if $$[x,x]=0$$ and $$[[x,y],z]+[[y,z],x]+[[z,x],y]=0.$$  Left multiplication in a Lie algebra is denoted by $\hbox{ad}\, x$: $\hbox{ad}\, x(y)=[x,y].$
 An associative algebra $A$ becomes a Lie algebra $A^-$ under the product, $[x,y]=xy-yx$.
 
The first axiom implies that $[x,y]=-[y,x]$ and the second (called the {\it Jacobi identity}) implies that 
$\hbox{ad}\,  x$ is a derivation, called an {\it inner derivation}, and $x\mapsto \hbox{ad}\,  x$ is a homomorphism of $\cl$ into the Lie algebra $(\hbox{End}\, \cl)^-$, that is, $\hbox{ad}\, [x,y]=[\hbox{ad}\, x,\hbox{ad}\, y]$.
 Assuming that $\cl$ is finite dimensional, the Killing form is defined by
$\lambda(x,y)=\tr{\hbox{ad}(x)\hbox{ad}(y})$. The proof of the following theorem is not given in Meyberg's notes \cite{Meyberg72},  but can be found in many places, for example \cite{Jacobson62}.

\begin{theorem}[CARTAN criterion \cite{Meyberg72}]
A finite dimensional Lie algebra $\cl$ over a field of characteristic 0 is semisimple if and only if the Killing form is nondegenerate.
\end{theorem}

 A more general form of the following theorem is due to Hochschild \cite{Hochschild42}. Theorem~\ref{thm:1112114}
has been attributed to Zassenhaus\footnote{  {\bf Hans Zassenhaus} (1912--1991)
 Hans Julius Zassenhaus  was a German mathematician, known for work in many parts of abstract algebra, and as a pioneer of computer algebra} in \cite{Meyberg72}, and to E. Cartan\footnote{{\bf Elie Cartan 1869--1951.} Elie Joseph Cartan was an influential French mathematician, who did fundamental work in the theory of Lie groups and their geometric applications. He also made significant contributions to mathematical physics, differential geometry, and group theory. He was the father of another influential mathematician, Henri Cartan.
} in \cite{Jacobson49}. We have taken the proof from \cite[p.42]{Meyberg72}.
  
\begin{theorem}\label{thm:1112114}
If the finite dimensional Lie algebra  $\cl$ over a field of characteristic 0 is semisimple (that is, its Killing form is nondegenerate), then every derivation is inner. 
\end{theorem}
\pf\
Let $D$ be a derivation of $\cl$.  Then $D\cdot \hbox{ad}(x)=\hbox{ad}(Dx)+\hbox{ad}(x)\cdot D$.  
Since $x\mapsto\tr{D\cdot \ad  x}$ is a linear form, there exists $d\in\cl$ such that $\tr{D\cdot \ad x}=\lambda(d,x)=\tr{\ad d\cdot \ad x}$.  Let $E$ be the derivation $E=D-\ad d$ so that 
\begin{equation}
\label{eq:1016101}
\tr{E\cdot\ad x}=0.
\end{equation}
Note next that 
\begin{eqnarray*}
E\cdot [\ad x,\ad y]&=&E\cdot\ad x\cdot\ad y-E\cdot\ad y\cdot \ad x\\
&=&(\ad x\cdot E+[E,\ad x])\cdot \ad y-E\cdot \ad y\cdot \ad x
\end{eqnarray*} so that 
\begin{eqnarray*}
[E,\ad x]\cdot \ad y&=&E\cdot [\ad x,\ad y]-\ad x\cdot E\cdot\ad y+E\cdot\ad y\cdot\ad x\\
&=&E\cdot [\ad x,\ad y]+[E\cdot\ad y,\ad x]
\end{eqnarray*}
and
$$\tr{[E,\ad x]\cdot \ad y}=\tr{E\cdot[\ad x,\ad y]}.$$
However,  since $E$ is a derivation 
\begin{eqnarray*}
[E,\ad x]\cdot \ad y&=&E\cdot \ad x\cdot \ad y-\ad x\cdot E\cdot \ad y\\
&=&(\ad{Ex}+\ad x\cdot E)\cdot \ad y-\ad x\cdot E\cdot \ad y\\
&=&\ad{Ex}\cdot \ad y.
\end{eqnarray*}
Thus
\begin{eqnarray*}
\lambda(Ex,y)&=&\tr{\ad{Ex}\cdot \ad y}\\
&=&\tr{[E,\ad x]\cdot \ad y}\\
&=&\tr{E\cdot [\ad x,\ad y]}=0 \hbox{ by (\ref{eq:1016101}))}.
\end{eqnarray*} 

Since $x$ and $y$ are arbitrary, $E=0$ and so $D-\ad d=0$.

\subsubsection{Circle multiplication---Jordan algebras}

We next consider a multiplication on matrices (or associative algebras) which leads to the theory of Jordan algebras\footnote{{\bf Pascual Jordan} (1902--1980)
Pascual Jordan was a German theoretical and mathematical physicist who made significant contributions to quantum mechanics and quantum field theory}

\begin{definition}
 A \underline{derivation} on $M_n(\CC)$ with respect to circle multiplication 
  $$
 X\circ Y=(XY+YX)/2
 $$
is a linear mapping  $\delta$ which satisfies the product rule
 $$\delta(A\circ B)=\delta(A)\circ B+A\circ \delta(B)$$
\end{definition}

 \begin{proposition} \label{prop:0201121}
Fix a matrix  $A$ in $M_n(\CC)$ and define $$\delta_A(X)=AX-XA.$$ Then $\delta_A$ is a derivation with respect to circle multiplication.
\end{proposition}

The first part of  the following theorem of Sinclair is valid for all semisimple Banach algebras (and continuous Jordan derivations).
 \begin{theorem}[Sinclair 1970 \cite{Sinclair70}] \label{thm:0201121} Every derivation on $M_n(\CC)$ with respect to circle multiplication  is also a derivation with respect to matrix multiplication, and is therefore of the form $\delta_A$ for some $A$ in $M_n(\CC)$.
 \end{theorem}
 
 According to the pioneering work of Jacobson\footnote{{\bf Nathan Jacobson} (1910--1999) Nathan Jacobson  was an American mathematician who was recognized as one of the leading algebraists of his generation, and he was also famous for writing more than a dozen standard monographs} the above proposition and theorem need to be modified for the real subalgebra (with respect to circle multiplication) of Hermitian matrices.

\begin{proposition}\label{prop:1.15}
Fix two hermitian matrices  $A,B$ in $M_n(\CC)$ and define $$\delta_{A,B}(X)=A\circ (B\circ X)-B\circ (A\circ X).$$ Then $\delta_{A,B}$ is a derivation of the real Jordan algebra of hermitian matrices with respect to circle multiplication.
\end{proposition}

 In preparation for the next theorem, we first note that for any algebra, $D$ is a derivation if and only if $[R_a,D]=R_{Da}$, where $R_a$ denotes the operator of multiplication (on the right) by $a$.
If you polarize the Jordan axiom $(a^2b)a=a^2(ba)$, you get
$[R_a,[R_b,R_c]]=R_{A(b,a,c)}$ where $A(b,a,c)=(ba)c-b(ac)$ is the ``associator''.
From the commutative law $ab=ba$, you get
$A(b,a,c)=[R_b,R_c]a$ and so $[R_b,R_c]$ is a derivation, sums of which are called {\bf inner},  
forming an  ideal in the Lie algebra of all derivations. 
The {\bf Lie multiplication algebra} $L$ of the Jordan algebra $A$ is the Lie algebra generated by the multiplication operators $R_a$.  It is given by
$$L=\{R_a+\sum_i [R_{b_i},R_{c_i}]:a,b_i,c_i\in A\}$$ so that $L$ is the sum of a Lie triple system (see subsection~\ref{3.2.2}) and the \underline{ideal} of inner derivations.

\begin{theorem}[1949 Jacobson \cite{Jacobson49}] Every derivation of a finite dimensional semisimple Jordan algebra  $J$ (in particular $M_n(\CC)$ with circle multiplication) is a sum of derivations of the form $\delta_{A,B}$ for some $A's$ and $B's$ elements of $J$.
 \end{theorem}

\pf\   Let $D$ be a derivation of a semisimple finite dimensional unital Jordan algebra $A$.  Then $\tilde D:X\mapsto [X,D]$ is a derivation of $L$.  
It is well known to algebraists that $L=L'+C$ where $L'$ (the derived algebra $[L,L]$) is semisimple and $C$ is the center of $L$.  Also $\tilde D$ maps $L'$ into itself and $C$ to zero. 

 By Theorem~\ref{thm:1112114}, $\tilde D$ is an inner derivation of $L'$ and hence also of $L$, so there exists $U\in L$ such that $[X,D]=[X,U]$ for all $X\in L$ and in particular
$[R_a,D]=[R_a,U]$.  
Then $Da=R_{Da}1=[R_a,D]1=[R_a,U]1=(R_aU-UR_a)1=a\cdot U1-Ua$ so that $D=R_{U1}-U\in L$.  
Thus, $D=R_a+\sum [R_{b_i},R_{c_i}]$ and so  $$0=D1=a+0=a.$$
\bigskip

 We summarize the previous three theorems, insofar as they concern $M_n(\CC)$, in the following table. Since $M_n(\CC)$ is not a semisimple Lie algebra, its  derivations must be taken modulo the center
(See \ref{sub:miers}).
 
 \begin{center}
{\bf Table 1}\nopagebreak

\smallskip

$M_n(\CC)$ (SEMISIMPLE ALGEBRAS)\nopagebreak

\smallskip

\begin{tabular}{|c|c|c|c|}\hline
matrix & bracket &circle\\\hline
$ab=a\times b$&$[a,b]=a b-b a$&$a\circ b=ab+ba$\\\hline
Th.~\ref{thm:1112112} &Th.~\ref{thm:1112114} &Th.~\ref{thm:0201121} \\\hline
 $\delta_a(x)$&$\delta_a(x)$&$ \delta_a(x)$\\
=&=&=\\
 $ax-xa$&$ax-xa$&$ax-xa$\\\hline
\end{tabular}

\end{center}
\bigskip

The following table  shows some properties of the three kinds of multiplication considered and thus provides the axioms for various kinds of algebras which we will consider.  These two tables also provide a model for the discussion of various kinds of triple systems which will be considered later.

\medskip

\begin{center}

\smallskip\smallskip\smallskip

{\bf Table 2}

\smallskip

ALGEBRAS
\smallskip\smallskip

{\bf commutative algebras}\\
$ab=ba$

\smallskip\smallskip

{\bf associative algebras}\\
$ a(bc)=(ab)c$\\

\smallskip\smallskip
{\bf Lie algebras}\\
$ a^2=0$
\\ $(ab)c+(bc)a+(ca)b=0$\\

\smallskip\smallskip
{\bf Jordan algebras}\\
$ ab=ba$\\
 $a(a^2b)=a^2(ab)$

\end{center}

 
 The algebra  $M_n(\CC)$, with matrix multiplication,  as well as the algebra  $\mathcal C$,  with ordinary multiplication, are examples of $C^*$-algebras  (finite dimensional and commutative respectively). 
The following theorem, which is due to Kadison\footnote{{\bf Richard Kadison} (b. 1925) Richard V. Kadison  is an American mathematician known for his contributions to the study of operator algebras}  and Sakai, and was preceded by some preliminary results of Kaplansky\footnote{ {\bf Irving Kaplansky } (1917--2006)
Kaplansky made major contributions to group theory, ring theory, the theory of operator algebras and field theory}
  thus explains Theorems~\ref{thm:1112111} and ~\ref{thm:1112112}.  We will have much more to say about C$^*$-algebras in subsection \ref{5.1}.
 
\begin{theorem}[Sakai, Kadison 1966 \cite{Sakai66},\cite{Kadison66}]\label{thm:1112113} Every derivation of a 
 $C^*$-algebra is of the form  $x\mapsto ax-xa$ for some  $a$ in the weak operator closure of the  $C^*$-algebra.
 \end{theorem}
  
 Why are derivations important?    One answer is given by cohomology, which is introduced here and revisited in later sections as the principal purpose of this paper.

Let $M$ be a Banach algebra and $X$ a Banach $M$-module.
For $n\ge 1$, let
$L^n(M,X)=$ all \underline{continuous} $n$-linear maps ($L^0(M,X)=X$).
The coboundary operator is  $\partial:L^n\rightarrow L^{n+1}$ (for $n\ge 1$), defined by
$$\partial\phi(a_1,\cdots,a_{n+1})=a_1\phi(a_2,\cdots,a_{n+1})
+\sum(-1)^j\phi(a_1,\cdots,a_{j-1},a_ja_{j+1},\cdots,a_{n+1})+(-1)^{n+1}\phi(a_1,\cdots,a_n)a_{n+1}$$
For $n=0$,
$\partial:X\rightarrow L(M,X)$\quad\quad $\partial x(a)=ax-xa$
so
$\hbox{Im}\partial=$ the space of inner derivations. Since $\partial\circ\partial=0$,
 $\hbox{Im}(\partial:L^{n-1}\rightarrow L^n)\subset \hbox{ker}(\partial:L^n\rightarrow L^{n+1})$,
 $H^n(M,X)=\hbox{ker}\partial/\hbox{Im}\partial$ is a vector space.
For $n=1$,
$\hbox{ker}\partial=
\{\phi:M\rightarrow X: a_1\phi(a_2)-\phi(a_1a_2)+\phi(a_1)a_2=0\}$
= the space of continuous derivations from $M$ to $X$
Thus, $$H^1(M,X)=\frac{\hbox{ derivations from }M\hbox{ to }X}{\hbox{inner derivations from }M\hbox{ to }X}$$ 
measures how close continuous derivations are to inner derivations.

Later we shall discuss what the spaces $H^2(M,X)$, $H^3(M,X)$,\ldots measure.  In anticipation of this, we present some miscellaneous remarks.

\begin{enumerate}

\item  One of the earliest papers on continuous cohomology is  the pioneering \cite{Kamowitz62} in which it is shown that
 $H^1(C(\Omega),E)=H^2(C(\Omega),E)=0$
\item ``The major open question in the theory of derivations on C*-algebras'' is whether $H^1(A,B(H))=0$ ($A\subset B(H)$)
\item A derivation from $A$ into $B(H)$ is inner if and only if it is completely bounded. (\cite[Theorem 3.1]{Christensen82})
\item The origin of continuous cohomology is the memoir of Barry Johnson \cite{Johnson72}\footnote{{\bf Barry Johnson} (1937--2002) see \cite{Sinclair04} for his fundamental contributions to cohomology of Banach algebras, among other things}

\end{enumerate}

\subsection{Exercises (Gradus ad Parnassum)---Algebras}

\begin{enumerate}

\item Prove the analog of Proposition~\ref{prop:0202121} for associative algebras: Fix an element $a$ in an associative algebra $A$ and define $$\delta_a(x)=[a,x]=ax-xa.$$ Then $\delta_a$ is a derivation with respect to the product of $A$.

\item Prove the analog of Proposition~\ref{prop:0202122} for associative algebras: Fix an element $a$ in an associative algebra $A$ and define $$\delta_a(x)=[a,x]=ax-xa.$$ Then $\delta_a$ is a derivation with respect to bracket multiplication $[x,y]=xy-yx$.

\item Prove the analog of Proposition~\ref{prop:0201121} for associative algebras: Fix an element $a$ in an associative algebra $A$ and define $$\delta_a(x)=[a,x]=ax-xa.$$ Then $\delta_a$ is a derivation with respect to circle multiplication $x\circ y=xy+yx$.

\item 
Let $a$ and $b$ be two fixed elements in an associative algebra $A$.  Show that the linear mapping
 $$\delta_{a,b}(x)=a\circ (b\circ x)-b\circ (a\circ x)$$ is a derivation of $A$ with respect to circle multiplication $x\circ y=xy+yx$.  

(cf. Proposition~\ref{prop:1.15})

\item Show that an associative algebra is a Lie algebra with respect to bracket multiplication.  

\item Show that an associative algebra  is a Jordan algebra with respect to circle multiplication.  

\item
Let us write $\delta_{a,b}$ for the linear mapping $\delta_{a,b}(x)=a(bx)-b(ax)$ in a Jordan algebra. Show that $\delta_{a,b}$ is a derivation of the Jordan algebra by following the outline below. (cf.  problem 4 above.)

\smallskip

(a) In the Jordan algebra axiom $$u(u^2v)=u^2(uv),$$ replace $u$ by $u+w$ to obtain the equation
\begin{equation}\label{eq:1.4}
2u((uw)v)+w(u^2v)=2(uw)(uv)+u^2(wv)\end{equation}

\smallskip

(b)  In (1), interchange $v$ and $w$ and subtract the resulting equation from  (1)
to obtain the equation
\begin{equation}\label{eq:1.5}
2u(\delta_{v,w}(u))=\delta_{v,w}(u^2).
\end{equation}

\smallskip

(c)  In (2), replace $u$ by $x+y$ to obtain the equation

$$
\delta_{v,w}(xy)=y\delta_{v,w}(x)+x\delta_{v,w}(y),
$$
which is the desired result.

\end{enumerate}

\section{Triple systems}

In this section, all algebras and triple systems will generally be finite dimensional.
 The first subsection parallels subsection~\ref{sub:1.1}.
 Infinite dimensional cases and references for the proofs  will be given in section~\ref{sec:5}.

\subsection{Derivations on finite dimensional triple systems}
\subsubsection{Triple matrix multiplication}

Let $M_{m,n}(\CC)$ denote the linear space of rectangular $m$ by $n$ complex matrices and note that it is closed under the operation $(A,B,C)\mapsto AB^*C$, which we will call triple matrix multiplication.
\begin{definition} 
 A \underline{derivation} on $M_{m,n}(\CC)$ with respect to triple matrix multiplication
is a linear mapping  $\delta$ which satisfies the triple product rule
 $\delta(AB^*C)= \delta(A) B^* C+A \delta(B)^* C+A B^* \delta(C)$
\end{definition}

\begin{proposition} \label{prop:0203121}
 For two matrices, $A\in M_m(\CC),B \in M_{n}(\CC)$, with $A^*=-A,B^*=-B$, 
define $\delta_{A,B}(X)=
AX+X B$.
Then  $\delta_{A,B}$ is a derivation with respect to triple matrix multiplication.
\end{proposition}
  
\begin{theorem}\label{thm:0202123} Every derivation on $M_{m,n}(\CC)$  with respect to triple matrix multiplication is of the form  $\delta_{A,B}$.
\end{theorem}

The proof of Theorem~\ref{thm:0202123} can be obtained by applying the result of 
Theorem~\ref{thm:0202125} below to the symmetrized product $(AB^*C+CB^*A)/2$.

\begin{remark}
 These results hold true and are of interest for the case $m=n$.
 \end{remark}
\subsubsection{Triple bracket multiplication}
  
Let's go back for a moment to square matrices and the bracket multiplication.  Motivated by the last remark, we define the triple bracket multiplication to be $(X,Y,Z)\mapsto [[X,Y],Z]$.

\begin{definition}
A \underline{derivation} on  $M_n(\CC)$ with respect to triple bracket multiplication is a linear mapping $\delta$ which satisfies the triple product rule
  $\delta([[A,B],C])=[[\delta(A),B],C]+[[A,\delta(B)],C]+[[A,B],\delta(C)]$ 
 \end{definition}
 
\begin{proposition}\label{prop:0203122}
Fix two matrices $A,B$ in $M_n(\CC)$ and define $\delta_{A,B}(X)=[[A,B],X]$.
Then  $\delta_{A,B}$ is a derivation with respect to triple bracket multiplication.
\end{proposition}

 $M_n(\CC)$ is not a semisimple Lie algebra or semisimple Lie triple system, so a modification needs to be made when considering Lie derivations  or Lie triple derivations on it, as in 
 Theorems~\ref{thm:miers} and \ref{thm:0202122}.  The proof of the following theorem is taken from \cite[Chapter 6]{Meyberg72}.  We do not define Lie triple system or semisimple here.

 \begin{theorem}\label{thm:0202124}
Every derivation of  a finite dimensional semisimple Lie triple system  $F$ is a sum of derivations of the form $\delta_{A,B}$, for some A's and B's in the triple system.  These derivations are called inner derivations and their set is denoted $\mbox{Inder}\, F$.
\end{theorem}
\pf\
Let $F$ be a finite dimensional semisimple Lie triple system (over a field of characteristic 0) and suppose that $D$ is a derivation of $F$. 
 Let $L$ be the Lie algebra $(\mbox{Inder}\, F)\oplus F$ with product
\[
[(H_1,x_1),(H_2,x_2)]=
([H_1,H_2]+L(x_1,x_2),H_1x_2-H_2x_1).
\]
A derivation of $L$ is defined by $\delta(H\oplus a)=[D,H]\oplus Da$.  
  Together with the definition of semisimple Lie triple system, it is proved in \cite{Meyberg72} that $F$ semisimple implies $L$ semisimple.
Thus there exists $U=H_1\oplus a_1\in L$ such that $\delta(X)=[U,X]$ for all $X\in L$. 
 Then
$0\oplus Da=\delta(0\oplus a)=[H_1+a_1,0\oplus a]=L(a_1,a)\oplus H_1a$ so $L(a_1,a)=0$ and $D=H_1\in \hbox{Inder}\, F$.
 \subsubsection{Triple circle multiplication}

Let's now return to rectangular matrices and form the triple circle multiplication
 $(A B^* C+C B^* A)/2$.
 For sanity's sake, let us  write this as 
 $$\{A,B,C\}=(A B^* C+C B^* A)/2$$

\begin{definition}
A \underline{derivation} on $M_{m,n}(\CC)$ with respect to triple circle multiplication is a linear mapping $\delta$ which satisfies the triple product rule
$\delta(\{A,B,C\})= \{\delta(A),B,C\}+\{A,\delta(B),C\}+\{B,A,\delta(C)\}$
\end{definition}

\begin{proposition}\label{prop:0203123}
Fix two matrices $A,B$ in $M_{m,n}(\CC)$ and define  $$\delta_{A,B}(X)=\{A,B,X\}-\{B,A,X\}.$$
Then  $\delta_{A,B}$ is a derivation with respect to triple circle multiplication.
\end{proposition}

The following theorem is   a special case of Theorem~\ref{thm:0601121} below.
 
 \begin{theorem}\label{thm:0202125}
Every derivation of  $M_{m,n}(\CC)$ with respect to triple circle multiplication is a sum of derivation of the form $\delta_{A,B}$.
\end{theorem}

\medskip

 We summarize Theorems~\ref{thm:0202123}, \ref{thm:0202124} and \ref{thm:0202125}, insofar as they concern $M_n(\CC)$, in the following table. Since $M_n(\CC)$ is not a semisimple Lie triple system, its  derivations must be taken modulo the center
(cf. \ref{sub:miers}).

\begin{center}

{\bf Table 3}

\smallskip

$M_{m,n}(\CC)$ (SS TRIPLE SYSTEMS)

\smallskip

\begin{tabular}{|c|c|c|c|}\hline
triple &triple &triple\\
matrix & bracket &circle\\\hline
$ab^*c$&$[[a,b],c]$&$ab^*c+cb^*a$\\\hline
Th. \ref{thm:0202123} &Th.~\ref{thm:0202124} &Th.~\ref{thm:0202125} \\\hline
 $\delta_{a,b}(x)$&$\delta_{a,b}(x)$& $\delta_{a,b}(x)$\\
=&=&=\\
 $ab^*x$&$abx$&$ab^*x$\\
 $+xb^*a$&$+xba$&$+xb^*a$\\
 $-ba^*x$&$-bax$&$-ba^*x$\\
 $-xa^*b$&$-xab$&$-xa^*b$\\\hline
 (sums)&(sums)&(sums)\\
             &($m=n$)&\\\hline
\end{tabular}

\end{center}

\subsection{Axiomatic approach for triple systems}

\begin{definition}
A  \underline{triple system} is defined to be a vector space with one ternary operation called triple multiplication.
Addition is denoted by $a+b$  and is commutative and associative:
$$a+b=b+a, \quad (a+b)+c=a+(b+c).$$
\end{definition}

Triple multiplication is denoted (simply and temporarily) by $abc$  and is required to be linear (conjugate-linear in some cases) in each of its three variables:
$$(a+b)cd=acd+bcd,\quad
 a(b+c)d=abd+acd,\quad
ab(c+d)=abc+abd.$$

Simple but important examples of triple systems can be formed from any algebra: if $ab$ denotes the algebra product, just define a triple multiplication to be
 $(ab)c$
\smallskip

Let's see how this works in the algebras we introduced in subsection~\ref{sub:1.1}.

\smallskip

$\mathcal C$; $fgh=(fg)h$, OR $fgh=(f\overline g)h$

\smallskip

$(M_n(\CC),\times)$; $abc=abc$ OR $abc=ab^*c$

\smallskip

$(M_n(\CC),[,])$; $abc=[[a, b], c] $

\smallskip

$(M_n(\CC),\circ)$; $abc=(a\circ b)\circ c $
\smallskip

All of these except the last one are useful.  It turns out that the appropriate form of the triple product for Jordan algebras is $abc=(a\circ b)\circ c+(c\circ b)\circ a-(a\circ c)\circ b$, since then you obtain a Jordan triple system (see subsection~\ref{3.2.3}). Also for the second example, by taking the symmetrized triple product,  in each case you also obtain a Jordan triple system.

\begin{definition}{\rm 
A triple system is said to be \underline{associative} if the triple multiplication is associative in the sense that $ab(cde)=(abc)de=a(bcd)e$ ({\bf first kind}) or  $ab(cde)=(abc)de=a(dcb)e$ ({\bf second kind}); and \underline{commutative} if  $abc=cba$.}
\end{definition}

\subsubsection{Associative triple systems}
The axiom which characterizes triple matrix multiplication is $$ (abc)de=ab(cde)=a(dcb)e.$$
The triple systems so defined are called {\bf associative triple systems} (of the second kind) or {\bf Hestenes\footnote{{\bf Magnus Hestenes (1906--1991)}
Magnus Rudolph Hestenes was an American mathematician. Together with Cornelius Lanczos and Eduard Stiefel, he invented the conjugate gradient method} algebras}

 \begin{theorem}[Lister 1971 \cite{Lister71}] Every derivation of a finite dimensional semisimple associative triple system (of the first kind) is inner.
\end{theorem}

\begin{theorem}[Carlsson 1976 \cite{Carlsson76}]
Every derivation of a finite dimensional semisimple associative triple system (first or second kind) into a module, is inner.
\end{theorem}

The following theorem really belongs in section~\ref{sec:5} since it concerns Banach associative triple systems. TROs play an important role in Part III, the definition is in subsection~\ref{TRO}.

\begin{theorem}[Zalar 1995 \cite{Zalar95}]
Let $W\subset B(H,K)$ be a TRO which contains all the compact operators.  If $D$ is a derivation of 
$W$ with respect to the associative triple product $ab^*c$ then there exist $a=-a^*\in B(K)$ and $b=-b^*\in B(H)$ such that $Dx=ax+xb$.
 \end{theorem}

This result has been extended to $B(X,Y)$ ($X,Y$ Banach spaces) in \cite{VelVil98}.

\subsubsection{Lie triple systems}\label{3.2.2}
The axioms which characterize triple bracket multiplication are
$$aab=0,\quad  
abc+bca+cab=0,$$ and
$$de(abc)=(dea)bc+a(deb)c+ab(dec).$$ The triple systems so defined are called {\bf Lie triple systems} and were developed by Jacobson and Koecher\footnote{{\bf Max Koecher (1924--1990)}
Max Koecher  was a German mathematician.
His main research area was the theory of Jordan algebras, where he introduced the Kantor-Koecher-Tits construction}.

\begin{theorem}[Lister 1952 \cite{Lister52}]
Every derivation of a finite dimensional semisimple Lie triple system is inner.
\end{theorem}

\subsubsection{Jordan triple systems} \label{3.2.3}
The axioms which characterize triple circle multiplication are $$abc=cba,$$ and 
$$de(abc)=(dea)bc-a(edb)c+ab(dec).$$ The triple systems so defined are called {\bf Jordan triple systems} and were developed principally by 
{\bf Kurt Meyberg, Ottmar Loos and Erhard Neher}, among many others.\smallskip

\begin{theorem}\label{thm:0601121} Every derivation of a finite dimensional  semisimple Jordan triple system is inner.
\end{theorem}

We outline a proof following a construction in  \cite{Chu12} (cf. \cite[Chapter 11]{Meyberg72}).
For simplicity, we assume non-degeneracy of a
Jordan triple system \cite[p.\ 25]{Chu12}.
\def\sqr#1#2{{\,\vcenter{\vbox{\hrule height.#2pt\hbox{\vrule width.#2pt
height#1pt \kern#1pt\vrule width.#2pt}\hrule height.#2pt}}\,}}
\def\bo{\sqr74\,}
Let $V$ be a Jordan triple and let $\mathcal L (V)$ be its TKK Lie
algebra ({\bf Tits-Kantor-Koecher}).
$\mathcal L (V) = V \oplus V_0 \oplus V$
and the Lie product is given by
$$[(x,h,y), (u,k,v)] = (hu - kx,\, [h,k]+x \bo v -u\bo y,\,
k^\natural y - h^\natural v).$$
Here, $a\bo b$ is the left multiplication operator $x\mapsto \{abx\}$ (also called the box operator),
 $V_0= {\rm span}\{ V\bo V\}$ is a Lie subalgebra of
$\mathcal L (V)$ and for $h= \sum_i a_i \bo b_i \in V_0$, the map
$h^\natural : V \rightarrow V$ is defined by
$$h^\natural = \sum_i b_i \bo a_i.$$

Let $\theta: \mathcal L(V) \rightarrow \mathcal
L(V)$ be the main involution 
$\theta(x\oplus h\oplus y)=y\oplus -h^\natural\oplus x$.
Let $\delta: V \rightarrow V$ be a derivation of a
Jordan triple $V$, with TKK Lie algebra $(\mathcal L(V), \theta)$.
Given $a,b\in V$, we define
\begin{eqnarray*}
D(a,0,0) &=& (\delta a,0,0)\\
D(0,0,b) &=& (0,0,\delta b)\\
D(0, a\bo b, 0) &=& (0,\, \delta a \bo b + a \bo \delta b,\, 0)
\end{eqnarray*}
and extend $D$ linearly on $\mathcal L(V)$. 
For the details of the proof of the following theorem and other results on Jordan triple cohomology, see the forthcoming paper \cite{ChuRus14}.
\begin{theorem}\label{3.19}
 Let $V$ be a Jordan triple with TKK Lie algebra $(\mathcal L(V), \theta)$.
There is a one-one correspondence between the triple derivations
of $V$ and the Lie derivations $D: \mathcal L(V) \rightarrow \mathcal
L(V)$ satisfying $D(V) \subset V$ and $D\theta = \theta D$.  Under this correspondence, $D$ is an inner Lie derivation if and only if $D|_V$ is an inner triple derivation.
\end{theorem}


Since the TKK Lie algebra $\mathcal L (V)$ of a semisimple Jordan triple system is semisimple,
Theorem~\ref{thm:0601121} follows from Theorem~\ref{thm:1112114} and Theorem~\ref{3.19}.
\smallskip

\medskip

 We summarize the previous three definitions in the following table
 
\medskip
 
\begin{center}

{\bf Table 4}

\smallskip

TRIPLE SYSTEMS

\smallskip

{\bf associative triple systems}\\
$ (abc)de=ab(cde)=a(dcb)e$\\

\smallskip\smallskip
{\bf Lie triple systems}\\
$aab=0$
\\ $abc+bca+cab=0$\\
$de(abc)=(dea)bc+a(deb)c+ab(dec)$

\smallskip\smallskip

{\bf Jordan triple systems} \\
$abc=cba$\\
$de(abc)=(dea)bc-a(edb)c+ab(dec)$

\end{center}

\subsection{Exercises (Gradus ad Parnassum)---Triple systems}

\begin{enumerate}
\item   Prove the analog of Proposition~\ref{prop:0203121} for associative algebras $A$ with involution:  For two elements $a,b\in A$ with $a^*=-a, b^*=-b$, define $\delta_{a,b}(x)=ax+xb$.
Then $\delta_{a,b}$ is a derivation with respect to the triple multiplication $ab^*c$.
  (Use the notation $\langle abc\rangle$ for $ab^*c$)
\item   Prove the analog of Proposition~\ref{prop:0203122} for associative algebras $A$:  Fix two
elements $a,b\in A$ and define $\delta_{a,b}(x)=[[a,b],x]$.  Then $\delta_{a,b}$ is a derivation with respect to the triple multiplication $[[a,b],c]$.
  (Use the notation $[abc]$ for $[[a,b],c]$)

\item  Prove Proposition~\ref{prop:0203123}:  Fix two matrices  $a,b$ in $M_{m,n}(\CC)$ and define $\delta_{a,b}(x)=\{a,b,x\}-\{b,a,x\}$.
Then  $\delta_{a,b}$ is a derivation with respect to triple circle multiplication.
(Use the notation $\{abc\}$ for $(ab^*c+cb^*a)/2$)

\item 
  Show that $M_n(\CC)$ is a Lie triple system with respect to triple bracket multiplication.  In other words, show that the three axioms for Lie triple systems in Table 4 are satisfied if $abc$ denotes $[[a,b],c]=(ab-ba)c-c(ab-ba)$ ($a,b$ and $c$ denote matrices).
  (Use the notation $[abc]$ for $[[a,b],c]$)

\item   
Show that $M_{m,n}(\CC)$ is a Jordan triple system with respect to triple circle multiplication.  In other words, show that the two axioms for Jordan triple systems in Table 4 are satisfied if $abc$ denotes $(ab^*c+cb^*a)/2$ ($a,b$ and $c$ denote matrices).
(Use the notation $\{abc\}$ for $(ab^*c+cb^*a)/2$)

\item 
Let us write $\delta_{a,b}$ for the linear process $$\delta_{a,b}(x)=abx$$ in a Lie triple system. Show that $\delta_{a,b}$ is a derivation of the Lie triple system by using the axioms for Lie triple systems in Table~4. 
  (Use the notation $[ abc]$ for the triple product in any Lie triple system, so that, for example, $\delta_{a,b}(x)$ is denoted by $[abx]$)

\item 
Let us write $\delta_{a,b}$ for the linear process $$\delta_{a,b}(x)=abx-bax$$ in a Jordan triple system. Show that $\delta_{a,b}$ is a derivation of the Jordan triple system by using the axioms for Jordan triple systems in Table 4. 
  (Use the notation $\{ abc\}$ for the triple product in any Jordan triple system, so that, for example, $\delta_{a,b}(x)=\{abx\}-\{bax\}$)

\item  
On the Jordan algebra $M_n(\CC)$ with the circle product $a\circ b=ab+ba$, define a triple product
$$
\{abc\}=(a\circ b)\circ c+(c\circ b)\circ a-(a\circ c)\circ b.
$$
Show that $M_n(\CC)$  is a Jordan triple system with this triple product.

Hint: show that $\{abc\}=2a bc+2c ba$

\item 
On the vector space $M_n(\CC)$, define a triple product $\langle abc\rangle=abc$ (matrix multiplication without the adjoint in the middle). Formulate the definition of a derivation of the resulting triple system, and state and prove a result corresponding to Proposition~\ref{prop:0203121}.
Is this triple system associative?

\item 
In an associative algebra, define a triple product $\langle abc\rangle$ to be $abc$.
Show that the algebra becomes an associative triple system with this triple product.

\item 
In an associative triple system with triple product denoted $\langle abc\rangle$, define a binary product $ab$ to be $\langle aub\rangle$, where $u$ is a fixed element.  Show that the triple system becomes an associative algebra with this product.  Suppose further that $\langle auu\rangle=\langle uua\rangle =a$ for all $a$. Show that we get a unital involutive algebra with involution $a^\sharp=\langle uau\rangle$.

\item 
 In a Lie algebra with product denoted by $[a,b]$, define a triple product $[abc]$ to be $[[a,b],c]$.  Show that the Lie algebra becomes a Lie triple system with this triple product.  (\cite[ch. 6, ex. 1, p. 43]{Meyberg72})

\item 
Let $A$ be an algebra (associative, Lie, or Jordan; it doesn't matter).  Show that the set ${\mathcal D}:=\hbox{Der}\, (A)$ of all derivations of $A$ is a Lie subalgebra of $\hbox{End}\, (A)$.  That is, $\mathcal D$ is a linear subspace of the vector space of linear transformations on $A$, and if $D_1,D_2\in \mathcal D$, then $D_1D_2-D_2D_1\in\mathcal D$.

\item 
Let $A$ be a triple system (associative, Lie, or Jordan; it doesn't matter).  Show that the set ${\mathcal D}:=\hbox{Der}\, (A)$ of derivations of $A$  is a Lie subalgebra of $\hbox{End}\, (A)$.  That is, $\mathcal D$ is a linear subspace of the vector space of linear transformations on $A$, and if $D_1,D_2\in \mathcal D$, then $D_1D_2-D_2D_1\in\mathcal D$.

\end{enumerate}

\subsection{Supplemental exercises (Gradus ad Parnassum)---Algebras and triple systems}

\begin{enumerate}
\item 
 In an arbitrary Jordan triple system, with triple product denoted by $\{abc\}$, define a triple product by
$$
[abc]=\{abc\}-\{bac\}.
$$
Show that the Jordan triple system becomes  a Lie triple system with this new triple product. \\
 (\cite[ch. 11, Th. 1, p. 108]{Meyberg72})

\item 
 In an arbitrary associative triple system, with triple product denoted by $\langle abc\rangle$, define a triple product by
$$
[xyz]=\langle xyz\rangle-\langle yxz\rangle-\langle zxy\rangle+\langle zyx\rangle.
$$
Show that the associative  triple system becomes  a Lie triple system with this new triple product. \\
 (\cite[ch. 6, ex. 3, p. 43]{Meyberg72})

\item 
 In an arbitrary Jordan algebra, with  product denoted by $xy$, define a triple product by
$
[xyz]=x(yz)-y(xz).
$
Show that the Jordan algebra becomes  a Lie triple system with this new triple product. \\
 (\cite[ch. 6, ex. 4, p. 43]{Meyberg72})

\item  
 In an arbitrary Jordan triple system, with triple product denoted by $\{abc\}$, fix an element $y$  and define a binary product by
$$
ab=\{ayb\}.
$$ 
Show that the Jordan triple system becomes a Jordan algebra with this (binary) product. \\
 (\cite[ch. 10, Th. 1, p. 94]{Meyberg72}---using different language;  see also \cite[Prop. 19.7, p. 317]{Upmeier85})

\item 
 In an arbitrary Jordan algebra with multiplication denoted by $ab$, define a triple product
$$
\{abc\}=(ab)c+(cb)a-(ac)b.
$$
Show that the Jordan algebra  becomes a  Jordan triple system with this triple product.  (\cite[ch. 10,  p. 93]{Meyberg72}---using different language; see also \cite[Cor. 19.10, p. 320]{Upmeier85})

\item   Show that every Lie triple system, with triple product denoted $[abc]$ is a subspace of some Lie algebra, with product denoted $[a,b]$, such that $[abc]=[[a,b],c]$. \\
 (\cite[ch. 6, Th. 1, p. 45]{Meyberg72})

\item Find out what a semisimple associative algebra is and prove that every derivation of a finite dimensional semisimple associative algebra is inner, that is, of the form $x\mapsto ax-xa$ for some fixed $a$ in the algebra. (\cite[Theorem 2.2]{Hochschild42})

\item Find out what a semisimple Lie algebra is and prove that every derivation of a finite dimensional semisimple Lie algebra is inner, that is, of the form $x\mapsto [a,x]$ for some fixed $a$ in the algebra.  \\
 (\cite[ch. 5, Th. 2, p. 42]{Meyberg72}; see also \cite[Theorem 2.1]{Hochschild42})

\item Find out what a semisimple Jordan algebra is and prove that every derivation of a finite dimensional semisimple Jordan algebra is inner, that is, of the form $x\mapsto \sum_{i=1}^n( a_i(b_ix)-b_i(a_ix))$ for some fixed elements $a_1,\ldots,a_n$ and $b_1,\ldots,b_n$ in the algebra. \cite[p. 320]{Jacobson68} and \cite[p. 270]{BraKoebook}

\item In an associative triple system with triple product $\langle xyz\rangle$, show that you get a Jordan triple system with the triple product $\{xyz\}=\langle xyz\rangle+\langle zyx\rangle$.   Then use Theorem~\ref{thm:0202125}
to prove Theorem~\ref{thm:0202123}.

\item Find out what a semisimple associative triple system is and prove that every derivation of a finite dimensional semisimple associative triple system is inner (also find out what inner means in this context).  (\cite{Carlsson76})

\item Find out what a semisimple Lie triple system is and prove that every derivation of a finite dimensional semisimple Lie triple system is inner, that is, of the form $x\mapsto \sum_{i=1}^n [a_ib_ix]$ for some fixed elements $a_1,\ldots,a_n$ and $b_1,\ldots,b_n$ in the Lie triple system. \\
 (\cite[ch. 6, Th.10, p. 57]{Meyberg72})

\item Find out what a semisimple Jordan triple system is and prove that every derivation of a finite dimensional semisimple Jordan triple system is inner, that is, of the form $x\mapsto \sum_{i=1}^n( \{a_ib_ix\}-\{b_ia_ix\})$ for some fixed elements $a_1,\ldots,a_n$ and $b_1,\ldots,b_n$ in the Jordan triple system. \\
 (\cite[ch. 11, Th. 8, p.123 and Cor. 2, p. 124]{Meyberg72})

\end{enumerate}

\section{Derivations on operator algebras and operator triple systems}\label{sec:5}
The theory of derivations on operator algebras is an important and well-investigated part of the general theory of operator algebras, with application in mathematical physics.  Investigation of unbounded derivations  were motivated mainly by the needs of mathematical physics, in particular by the problem of constructing the dynamics in quantum 
statistical mechanics. We do not consider unbounded derivations {\it per se} in this paper, cf. \cite{BraRob79,Sakai91}.

In this section we shall review the literature about derivations on associative and non associative operator algebras and operator triple systems, with the main focus being on some new results concerning the notions of ternary weak amenability and triple derivations on von Neumann algebras.

There are two basic questions concerning derivations of Banach algebras and Banach triple systems, namely, whether they are automatically continuous, and to what extent are they all inner, in the appropriate context.   The automatic continuity of various algebraic mappings plays
important roles in the general theory of Banach algebras and in
particular in operator algebra theory.  
These questions are significant in the case of a derivation of a space into itself, or into an appropriate module.\smallskip

We shall consider derivations primarily  in the following contexts: 

\begin{description}
\item[(i)]  C*-ALGEBRAS (more generally, associative Banach algebras)

\item[(ii)]  JC*-ALGEBRAS (more generally, Jordan Banach algebras)

\item[(iii)] JC*-TRIPLES (more generally, Banach Jordan triples)
\end{description}

One could and should also consider:
\begin{description}
 \item[(i$^\prime$)] Banach associative triple systems
\item[(ii$\prime$)] Banach Lie algebras  
 \item[(iii$^\prime$)] Banach Lie triple systems
 \end{description}

\subsection{C*-algebras}\label{5.1}

In the context of associative algebras, a derivation is a linear mapping satisfying the property $D(ab)=a\cdot Db+Da\cdot b$, and an inner derivation is a derivation of the form  $\hbox{ad x}(a)=x\cdot a-a\cdot x$.  

The main automatic continuity results for $C^*$-algebras are due to Kaplansky (commutative and type I cases), Sakai, and Ringrose\footnote{ {\bf John Ringrose (b. 1932)}
  John Ringrose is a leading world expert on non-self-adjoint operators and operator algebras. He has written  a number of influential texts including Compact non-self-adjoint operators (1971) and, with R V Kadison, Fundamentals of the theory of operator algebras in four volumes published in 1983, 1986, 1991 and 1992}.  
   The main inner derivation results are due primarily to Sakai, Kadison, Connes\footnote{ {\bf Alain Connes b. 1947}
 Alain Connes is  the leading specialist on operator algebras. 
 In his early work on von Neumann algebras in the 1970s, he succeeded in obtaining the almost complete classification of injective factors. 
  Following this he made contributions in operator K-theory and index theory, which culminated in the Baum-Connes conjecture.   
  He also introduced cyclic cohomology in the early 1980s as a first step in the study of noncommutative differential geometry.
Connes has applied his work in areas of mathematics and theoretical physics, including number theory, differential geometry and particle physics}
, and Haagerup\footnote{ {\bf Uffe Haagerup b. 1950}
Haagerup's research is in operator theory, and covers many subareas in the subject which are currently very active - random matrices, free probability, C*-algebras and applications to mathematical physics}.
We won't define the terms amenable and nuclear here.   We refer to the excellent notes by Runde \cite{RundeLNIM} for these definitions and much other information. Weak amenability will be defined shortly, as it is one of the main objects of our discussion.

\begin{theorem}[Sakai 1960 \cite{Sakai60}]
Every derivation from a
C$^*$-algebra into itself  is continuous.
\end{theorem}
\begin{theorem}[Ringrose 1972 \cite{Ringrose72}]
 Every derivation  from a
C$^*$-algebra  into a Banach $A$-{\bf bimodule}  is continuous.
\end{theorem}
\begin{theorem}[Sakai \cite{Sakai66}, Kadison \cite{Kadison66}]
Every derivation of a $C^*$-algebra is of the form 
 $x\mapsto ax-xa$ for some $a$ in the weak closere of the $C^*$-algebra
 \end{theorem}

\begin{theorem}[Connes 1976 \cite{Connes76}]
Every amenable $C^*$-algebra is nuclear.
 \end{theorem}

 \begin{theorem} [Haagerup 1983 \cite{Haagerup83}]
 Every nuclear $C^*$-algebra is amenable.
 \end{theorem}

 \begin{theorem}[Haagerup 1983 \cite{Haagerup83}] Every $C^*$-algebra is weakly amenable.
 \end{theorem}

\subsection{Algebras of unbounded operators}

\subsubsection{Murray-von Neumann algebras}
Physics considerations demand that
the  Hamiltonian of a quantum system will, in general, correspond to an unbounded
operator on a Hilbert space H. These unbounded operators will not lie in a von Neumann algebra, but they may be ÒaffiliatedÓ with
the von Neumann algebra corresponding to the quantum system.  In general, unbounded operators do not behave well with respect to addition and multiplication. However, as noted by Murray and von Neumann \cite{MurvN36}, for the finite von Neumann algebras, their families of affiliated operators do form a *-algebra. Thus, it natural to study derivations of algebras that
include such unbounded operators.

A closed densely defined operator
T on a Hilbert space H is affiliated with a von Neumann algebra R when $UT = TU$ for each unitary operator U in $R'$, the commutant
of R. If operators S and T are affiliated with R, then S + T and ST are densely defined, preclosed and their closures are affiliated with R.  Such algebras are referred to as {\bf Murray-von Neumann algebras}.

If R is a finite von Neumann algebra, we denote by $A_f(R)$ its associated Murray-von Neumann algebra. It is natural to conjecture that every derivation of  $A_f(R)$ should be inner.  In \cite[Theorem 5]{KadLiu14} it is proved that extended derivations
of $A_f(R)$ (those that map R into R) are inner. In \cite[Theorem 12]{KadLiu14}, it is  proved that each derivation of $A_f(R)$  with $R$ a factor
of type $II_1$ that maps $A_f(R)$  into $R$ is 0. In other words, by requiring that the range of the derivation is in $R$, the ÒboundedÓ
part of $A_f(R)$, allows a noncommutative unbounded version of the well-known  Singer-Wermer conjecture.  This is extended to the general von Neumann algebra of type $II_1$ in \cite{KadLiuMS}.
\smallskip

We turn our attention now to  commutators in $A_f(N)$, where $N$ is a finite von Neumann algebra.  The following theorem is related to the work of Kadison and Liu described above.
\begin{theorem}[Theorem 4.13 of \cite{Liu11}]
Let $N$ be $II_1$-factor. If the element $b$ of $N$ is  a commutator of self-adjoint elements in $A_f(N)$, then $b$ has trace zero.  In particular, the scalar operator $i1 \in  N$
is not a commutator of self-adjoint elements in $A_f(N)$.
\end{theorem}
This has been complemented in \cite{Thom14} as follows.
\begin{theorem}[Theorem 4 of \cite{Thom14}]
Let $a,b\in A_f(N)$ where $N$ is a $II_1$-factor.  Assume that either $a,b$ belong to the Haagerup-Schultz algebra {\rm (}a *-subalgebra of $A_f(N)$, see {\rm \cite{HagSch07})} or if $ab\in A_f(N)$ is conjugate to a self-adjoint element. If $[a,b]=\lambda 1$, then $\lambda=0$.
\end{theorem}

As for sums of commutators, we have the following theorem.

\begin{theorem}[Theorems 5 and 6 of \cite{Thom14}]
 Let $N$ be a $II_1$-factor. There exist $a,b,c,d\in A_f(N)$  such that
$1=[a,b]+[c,d]$. Every element of $A_f(N)$ is the sum of two commutators.
\end{theorem}

\subsubsection{Algebras of measurable operators}

Noncommutative integration theory was initiated by Segal (\cite{Segal53}), who considered new classes of algebras of unbounded operators, in particular the algebra $S(M)$ of all measurable operators  affiliated with a von Neumann algebra $M$.   A study of the derivations on the algebra $S(M)$ was initiated by Ayupov\footnote{{\bf Shavkat Ayupov  b. 1952}  is the Director of the Institute of Mathematics, National University of Uzbekistan, well-known for his work in operator algebras and quantum probability, as well as in the theory of Leibniz algebras} in 2000.  For example, in the commutative case $M=L^\infty(\Omega,\Sigma,\mu)$, $S(M)$ is isomorphic to the algebra $L^0(\Omega)$ of all complex measurable functions and it is shown in \cite{BerChiSuk06} that $L^0(0,1)$ admits nonzero derivations which are discontinuous in the measure topology. The study of derivations on various subalgebras of the algebra $LS(M)$ of all locally measurable operators in the general semifinite case was initiated in \cite{AyuAlbKud07,AyuAlbKud08}, with the most complete results obtained in the type I case, \cite{AyuAlbKud09}.  See the survey \cite{AyuKud10} for more detail up to 2010.

It was shown in \cite{AyuAlbKud09} (and in \cite{BerPagSuk11} for separable predual) that if $M$ is properly infinite and of type I, then every derivation of $LS(M)$ is continuous  in the local measure topology $t(M)$ on $LS(M)$. The same holds for $M$ of type III as shown in \cite{AyuKud10} 
\ and for type $II_\infty$ in \cite{BerChiSuk13}.  

The following is the main theorem of \cite{BerChiSuk14} and illustrates the state of the art.
\begin{theorem} Let $M$ be any von Neumann algebra.
Every derivation on the *-algebra $LS(M)$ continuous with respect to the topology $t(M)$ is inner.
\end{theorem}
\begin{corollary}
If $M$ is a properly infinite von Neumann algebra, then every derivation on $LS(M)$ is inner.
\end{corollary}

\subsection{A bridge to Jordan algebras}

A {\it Jordan derivation} from a Banach algebra $A$ into a Banach
$A$-module is a linear map $D$ satisfying $D(a^2) = a D(a) + D(a)
a,$ ($a\in A$), or equivalently, $D(ab+ba)=aD(b)+D(b)a + D(a)b
+bD(a),$ ($a,b\in A$). 
Sinclair proved in 1970 \cite{Sinclair70} that a bounded Jordan
derivation from a semisimple Banach algebra to itself is a
derivation, although this result fails for derivations of
semisimple Banach algebras into a Banach bi-module.

Nevertheless, a celebrated result of B.E.
Johnson in 1996 \cite{Johnson96} states that every bounded Jordan derivation from a
C$^*$-algebra $A$ to a Banach $A$-bimodule is an associative
derivation.
In view of the intense interest in automatic continuity problems
in the past half century, it was natural to consider the following problem.

\begin{problem}
 {\rm  Is every Jordan derivation from a C$^*$-algebra
$A$ to a Banach $A$-bimodule automatically continuous (and hence a derivation, by Johnson's theorem)?}
\end{problem}

In 2004, J.
Alaminos\hyphenation{Alaminos}, M. Bre\v{s}ar and A.R. Villena \cite{AlaBreVil04}
gave a positive answer to the above problem for some classes of
C$^*$-algebras including the
class of abelian C$^*$-algebras.
Combining a theorem of Cuntz from 1976 \cite{Cuntz76} with the theorem just quoted  yields

\begin{theorem}
Every Jordan derivation from a C$^*$-algebra $A$ to a Banach $A$-module  is continuous.
\end{theorem}

In the same way, using  the solution in 1996 by Hejazian-Niknam  \cite{HejNik96} in the commutative case we have

\begin{theorem}
Every Jordan derivation from a C$^*$-algebra $A$ to a \underline{\it Jordan} Banach $A$-module  is continuous.
\end{theorem}

(Jordan module will be defined below)
\smallskip

These last two results are among the consequences of some work of Peralta and Russo \cite{PerRus10}, described below,  on automatic continuity of derivations into Jordan triple modules.

\subsection{JC*-algebras}

In the context of Jordan algebras and Jordan modules, a derivation is a linear mapping satisfying the property $D(a\circ b)=a\circ Db+Da\circ b$, and an inner derivation is a derivation of the form  $\sum_i [L(x_i)L(a_i)-L(a_i)L(x_i)]$,  ($x_i\in M, a_i\in A$), that is,
$b\mapsto \sum_i [x_i\circ(a_i\circ b)-a_i\circ(x_i\circ b)]$

The main automatic continuity results for $JC^*$-algebras are due to Upmeier\footnote{  {\bf Harald Upmeier  (b. 1950)}, Professor at University of Marburg, Germany, pioneer of topological Jordan structures}, Hejazian-Niknam, and Alaminos-Bresar-Villena.  The main inner derivation results are due to Jacobson and Upmeier.
We have already mentioned the results of Hejazian-Niknam and Alaminos-Bresar-Villena in this context.
\begin{theorem}
[Upmeier 1980 \cite{Upmeier80}]
Every derivation of a reversible $JC^*$-algebra extends to a derivation of its enveloping $C^*$-algebra.
 \end{theorem}
In particular, this theorem improves Theorem~\ref{thm:0201121} in the case of a $C^*$-algebra.

\begin{theorem} [Jacobson 1951 \cite{Jacobson68}]
\label{thm:0208121}
 Every derivation of a finite dimensional semisimple Jordan algebra into a (Jordan) module is inner.
\end{theorem}
 The proof of this theorem uses the theory of Lie algebras and Lie triple systems.\smallskip

 Because of the structure theory for real Jordan operator algebras, the following theorem completely answers the question of when all derivations are inner on a JBW-algebra.
 
 \begin{theorem}[Upmeier 1980 (Theorem 3.10 in \cite{Upmeier80})]\label{4.14} 
 1.  Purely exceptional JBW-algebras have the inner derivation property, that is, every derivation is inner\\
2.   Reversible JBW-algebras have the inner derivation property\\  
3.  $\oplus L^\infty(S_j,U_j)$ has the inner derivation property if and only if $\sup_j\dim U_j<\infty$ ($U_j$ spin factors).
\end{theorem}
 
 \subsection{Lie derivations of operator algebras and Banach algebras}\label{sub:miers}

  In this subsection, we record some results on derivations of the Lie algebra structure of operator algebras, including  triple structure, which is one theme of this survey paper.  

\begin{theorem}[Miers 1973 \cite{Miers73}]
If $M$ is a von Neumann algebra, $[M,M]$ the Lie algebra linearly generated by $\{[X, Y ] = XY-YX:X,Y\in M\}$
 and $ L: [M,M]\rightarrow M$ a Lie derivation, i.e., $L$ is linear and $L[X, Y ] = [LX, Y ]+
[X,LY ]$, then  $L$ has an extension $D:M\rightarrow M $ that is a derivation of the
associative algebra. 
\end{theorem}
The proof involves matrix-like computations. 
Using Theorem~\ref{thm:1112113} leads to 
\begin{theorem} [Miers 1973 \cite{Miers73}] \label{thm:miers}If $L:M \rightarrow M$ is a Lie derivation, then $L = D+\lambda$,  where $D$ is an associative derivation and $\lambda$  is a linear map into the center of $M$ vanishing on
$[M,M]$. 
\end{theorem}

A Lie triple derivation of an associative algebra $M$ is a linear map $L:M \rightarrow
M$ such that $$L[[X, Y ],Z] = [[L(X), Y ],Z]+
[[X,L(Y )],Z]+[[X, Y ],L(Z)]$$ for all $X, Y,Z\in M$.

\begin{theorem}[Miers 1978 \cite{Miers78}]\label{thm:0202122}
If $M$ is a von Neumann algebra with no central Abelian summands and $L$ is a Lie triple derivation,  then there exists
an operator $A\in M$ such that $L(X) = [A,X]+ \lambda(X)$ where $\lambda :M \rightarrow Z_M$ is a linear map which
annihilates brackets of operators in $M$.
\end{theorem}
 
Notice that no assumption of continuity was made in these theorems.  Continuous Lie derivations of associative Banach algebras have also been studied outside the operator algebra context,  as evidenced in the following.

\begin{theorem}[Johnson 1996 \cite{Johnson96}]
Every continuous Lie derivation of a symmetrically amenable Banach algebra $A$ into a Banach  $A$-module $X$ is the sum of an associative derivation and a ``trivial'' derivation, that is, a linear map which vanishes on commutators and maps into the ``center'' of the module. The same holds if $A$ is a $C^*$-algebra.
\end{theorem}
 
As shown in \cite{MatVil03} and \cite{AlaMatVil04}, the continuity assumption can be dropped if $X=A$ and $A$ is a C*-algebra or a semisimple symmetrically amenable Banach algebra.\smallskip

According to a review by Garth Dales \cite{Dalesrev}: ``It remains an open question whether an analogous result (automatic continuity)  for Lie derivations from A into a Banach A-bimodule holds when A is an arbitrary C*-algebra and when A is an arbitrary symmetrically amenable Banach algebra.
     It is also an interesting open question whether or not every Lie derivation on a semisimple Banach algebra  to itself has this form.''\smallskip
  
As pointed out to the author  by Martin Mathieu, a negative answer has been given, by Read \cite{Read03}, to the last question raised by Dales in his review.
\smallskip




A search for ``Lie derivation'' in the ``review text''  option in MathSciNet turns up 172  entries. A good number of these concern Lie derivations on nonselfadjoint operator algebras, for example nest algebras, reflexive algebras, $CSL$ algebras, $UHF$-algebras\footnote{It is worth mentioning in connection with section~\ref{chpt6}, that there are corresponding results for the cohomology of nonselfadjoint operator algebras.  See references [14]-[17],[38],[39] in \cite{Ringrose96} as well as the recent papers  \cite{Hou10} and \cite{HouWei07}}.  Several others are concerned with von Neumann algebras.    We quote a few of them.

\begin{theorem}[Ji and Qi 2011 Corollary 2.1 in \cite{JiQi11}]
Let $N$ be an arbitrary nest on a Hilbert space $H$ of dimension greater than 2, and $\hbox{Alg} N$ be the associated nest algebra.  Let $\delta$ be a linear mapping on $\hbox{Alg} N$  which satisfies the Lie derivation property on commutators $[a,b]$ with $ab=0$.  Then there exists an operator $r\in \hbox{Alg} N$ and a linear center valued map $\tau$ vainshing on such commutators such that $\delta(a)=ra-ar+\tau(a)$.
\end{theorem}

\begin{theorem}[Zhuraev 2011 \cite{Zhuraev11}]
For a von Neumann algebra $M$ on a Hilbert space $H$ denote by $S(M)$ the unital $^*$-algebra of all measurable operators affiliated with $M$. If $M$ is a homogeneous von Neumann algebra of type $I_n$, then very Lie derivation on $S(M)$ can be uniquely decomposed into the sum of a derivation and a centre-valued trace.
\end{theorem}

In any algebra $A$  set $p_1(x)=x$ and for $n\ge 2$,
$
p_n(x_1,\ldots,x_n))=[p_{n-1}(x_1,\ldots,x_{n-1},x_n]
$.
Thus $p_2(x_1,x_2)=[x_1,x_2]$, $p_3(x_1,x_2,x_3)=[[x_1,x_2],x_3]$. A map $\varphi:A\rightarrow A$ is called a multiplicative Lie $n$-derivation if
\[
\varphi(p_n(x_1,\ldots,x_n))=\sum_i^n p_n(x_1,\ldots,\varphi(x_i),x_{i+1},\ldots,x_n).
\]

\begin{theorem}[Abdullaev 1992 \cite{Abdullaev92}]
Every $n$-Lie derivation $L$ on some von Neumann algebra $M$ (or on its skew-adjoint part) can be decomposed as $L=D+E$, where $D$ is an ordinary derivation on $M$ and $E$ is a linear map from $M$ into its center which annihilates the Lie products. 
\end{theorem}

\subsection{JC*-triples}

In the context of Jordan triple systems, a derivation is a linear  or conjugate linear mapping satisfying the property $D\{a,b,c\}=\{Da.b,c\}+\{a,Db,c\}+\{a,b,Dc\}$, and an inner derivation is a derivation of the form  $\sum_i [L(x_i,a_i)-L(a_i,x_i)]$  ($x_i\in M, a_i\in A$), that is, $b\mapsto \sum_i [\{x_i,a_i,b\}-\{a_i,x_i,b\}]$.  \smallskip

For reasons related to the definition of Jordan triple module in the complex case, and which are discussed in both \cite{HoPerRus11} and \cite{PerRus10}, although a derivation from a Jordan triple system to itself is  linear, nevertheless, in the complex case, a derivation of a Jordan triple system into a Jordan triple module (other than  the triple itself) is declared to be a conjugate linear map.\smallskip

The main automatic continuity results for $JC^*$-triples, where the triple product is given by
$\{x,y,z\}=(xy^*z+zy^*x)/2$
 are due to Barton\footnote{ {\bf Tom Barton (b. 1955)}
  Tom Barton is Senior Director for Architecture, Integration and CISO at the University of Chicago. He had similar assignments  at the University of Memphis, where he was a member of the mathematics faculty before turning to administration}
  and Friedman\footnote{{\bf Yaakov Friedman (b. 1948)}
   Yaakov Friedman is director of research at Jerusalem College of Technology};
   and Peralta\footnote{ {\bf Antonio Peralta (b. 1974)}, Professor at University of Granada, Spain} and Russo\footnote{ {\bf Bernard Russo (b. 1939)}, Professor Emeritus at University of California, Irvine}, and the  main inner derivation results are due to Ho-Martinez-Peralta-Russo, Meyberg, K\"uhn-Rosendahl, and Ho-Peralta-Russo.

Before stating these results, we must acknowledge the pioneering efforts of Harris\footnote{{\bf Lawrence A. Harris} (PhD 1969), Professor at University of Kentucky} to the study of Jordan operator triple systems, especially for the papers \cite{Harris74} and \cite{Harris81} devoted respectively to infinite dimensional holomorphy 
and
spectral and ideal theory.

\begin{theorem}[Barton-Friedman 1990 \cite{BarFri90}]
Every derivation of a $JB^*$-triple is continuous.
 \end{theorem}

\begin{theorem}[Peralta-Russo 2010 \cite{PerRus10}]
Necessary and sufficient conditions under which a derivation of a $JB^*$-triple into a Jordan triple module is continuous.
 \end{theorem}

(JB$^*$-triple and Jordan triple module are defined below)\smallskip

 The proofs of the following two theorems rely heavily on the theories of Lie algebras and Lie triple systems.

\begin{theorem}[Meyberg 1972 \cite{Meyberg72}]  Every derivation of a finite dimensional semisimple Jordan triple system is inner
\end{theorem}

\begin{theorem}[K\"uhn-Rosendahl 1978 \cite{KuhRos78}]   Every derivation of a finite dimensional semisimple Jordan triple system into a Jordan triple module is inner.
 \end{theorem}
 
  The proofs of the following two theorems rely on operator theory and functional analysis.
 
\begin{theorem}[(Ho-Martinez-Peralta-Russo 2002 (Theorem 2 in \cite{HoMarPerRus02})]\label{4.19}
 Cartan factors of type  $I_{n,n}$, II (even or $\infty$), and III have the inner derivation property. 
\end{theorem}
  
  \begin{theorem}[(Ho-Martinez-Peralta-Russo 2002 \cite{HoMarPerRus02}]  Infinite dimensional Cartan factors of type $I_{m,n}, m\ne n$, and IV do not have the inner derivation property.\end{theorem}
 


 

\subsubsection{
Ternary weak amenability}

Recall that a Banach algebra $A$ is said to be \emph{amenable} if every bounded
derivation of $A$ into a dual $A$-module is inner, and  \emph{weakly amenable} if
every (bounded) derivation from $A$ to $A^*$ is inner.
A Jordan Banach triple $E$ is said to be \emph{weakly amenable} or
\emph{ternary weakly amenable} if every continuous triple derivation from
$E$ into its dual space is necessarily inner. 
The main results of \cite{HoPerRus11} are:

 \begin{enumerate}
 \item  Commutative $C^*$-algebras are ternary weakly amenable.
\item Commutative $JB^*$-triples are approximately weakly amenable, that is, every derivation into its dual is approximated in norm by inner derivations.
\item  $B(H),K(H)$  are ternary weakly amenable if and only if they are finite dimensional.
\item  Cartan factors of type $I_{m,n}$ of finite rank with  $m\ne n$, and of type IV are ternary weakly amenable if and only if they are finite dimensional.
\end{enumerate}

We provide a taste of the details of one of the results of \cite{HoPerRus11}.

 \begin{lemma}
 The C$^*$-algebra $A=K(H)$ of all compact operators on an
infinite dimensional Hilbert space $H$  is not Jordan weakly amenable, that is, not every Jordan derivation into the dual is inner.
\end{lemma}
{\it Proof:}
We shall identify $A^*$ with the trace-class operators on $H$.
Supposing that $A$ were Jordan weakly amenable, let $\psi\in A^*$ be arbitrary.
Then $D_\psi$ ($=\hbox{ad}\, \psi$)  is an associative derivation and hence a Jordan derivation, so by assumption would be an inner Jordan derivation.  Thus  there would exist $\varphi_j\in A^*$
and $b_j\in A$ such that $$D_\psi(x)=\sum_{j=1}^n[\varphi_j\circ(b_j\circ x)-b_j\circ(\varphi_j\circ x)]$$
 for all $x\in A$.
 For $x,y\in A$, a direct calculation yields
$$ \psi(xy-yx)=-\frac{1}{4}\left(\sum_{j=1}^n b_j \varphi_j-\varphi_j b_j\right)(xy-yx). $$\smallskip

It is known (Pearcy-Topping 1971 \cite{PeaTop71}) 
that every compact operator on a separable (which we may assume without loss of generality) infinite dimensional Hilbert space  is
a finite sum of commutators of compact operators, and it follows that the trace-class operator
$$\psi = -\frac{1}{4}\left(\sum_{j=1}^n b_j \varphi_j-\varphi_j b_j\right)$$ is a finite sum of
commutators of compact and trace-class operators, and hence has trace zero.
This is a contradiction, since $\psi$ was arbitrary, completing the proof.

\subsubsection{Automatic continuity}\label{5.5}

We shall now give some details about the paper \cite{PerRus10} which has been mentioned several times already.  This involves Jordan triples and Jordan triple modules.


A complex (resp., real) \emph{Jordan triple} is a complex (resp., real)
vector space $E$ equipped with a non-trivial  triple
product $$ E \times E \times E \rightarrow E$$
$$(x,y,z) \mapsto \tp{x}{y}{z} $$
which is bilinear and symmetric in the outer variables and
conjugate linear (resp., linear) in the middle one satisfying the
so-called \emph{``Jordan Identity''}:
\begin{equation}\label{eq:0221121}
L(a,b) L(x,y) -  L(x,y) L(a,b) =
 L(L(a,b)x,y) - L(x,L(b,a)y),\end{equation}
for all $a,b,x,y$ in $E$, where $L(x,y) z := \tp{x}{y}{z}$.

A JB$^*$-algebra is a complex Jordan Banach algebra $A$ equipped
with an algebra involution $^*$ satisfying  $\|\{a,a^*,a\} \|= \|a\|^3$, $a\in
A$.  (Recall that $\{a,a^*,a\}  =
 2 (a\circ a^*) \circ a - a^2 \circ a^*$).
A  (complex) \emph{JB$^*$-triple} is a complex Jordan Banach triple
${E}$ satisfying the following axioms: 
\smallskip

(a) For each $a$ in ${E}$ the map $L(a,a)$ is an hermitian
operator on $E$ with non negative spectrum.
\smallskip

(b)  $\left\|
\{a,a,a\}\right\| =\left\| a\right\| ^3$ for all $a$ in ${A}.$
\medskip

Every C$^*$-algebra (resp., every JB$^*$-algebra) is a JB$^*$-triple with respect to the product
$\{abc\} = \frac12 \ ( a b^* c + cb^* a) $ (resp., $\{abc\} := (a\circ b^*) \circ c + (c\circ b^*) \circ a - (a\circ c) \circ b^*$).\smallskip

In order to motivate the definition of triple module, we first recall two basic definitions.
If $A$ is an associative algebra, an
\emph{$A$-bimodule} is a vector space $X$, equipped with two
bilinear products $(a,x)\mapsto a x$ and $(a,x)\mapsto x a$ from
$A\times X$ to $X$ satisfying the following axioms: $$a (b x) = (a
b) x ,\ \ a (x b) = (a x) b, \hbox{ and, }(xa) b = x (a b),$$ for
every $a,b\in A$ and $x\in X$.
If $J$ is a Jordan algebra, a \emph{Jordan $J$-module} is a
vector space $X$, equipped with two bilinear products
$(a,x)\mapsto a \circ x$ and $(x,a)\mapsto x \circ a$ from
$J\times X$ to $X$, satisfying: $$a \circ x = x\circ a,\ \ a^2
\circ (x \circ a) = (a^2\circ  x)\circ a, \hbox{ and, }$$ $$2(
(x\circ a)\circ  b) \circ a + x\circ (a^2 \circ b) = 2 (x\circ
a)\circ  (a\circ b) + (x\circ b)\circ a^2,$$ for every $a,b\in J$
and $x\in X$
\smallskip

 If $E$ is a complex  Jordan triple, a \emph{Jordan
triple $E$-module}  (also called \emph{triple $E$-module}) is a
vector space $X$ equipped with three mappings\medskip

$\{.,.,.\}_1 :
X\times E\times E \to X$\quad , \quad
$\{.,.,.\}_2 : E\times X\times E \to
X$\quad , \quad
$ \{.,.,.\}_3: E\times E\times X \to X$\smallskip
\\ satisfying:

\begin{enumerate}
\item $\{ x,a,b \}_1$ is linear in $a$ and $x$ and conjugate
linear in $b$, $\{ abx \}_3$ is linear in $b$
and $x$ and conjugate linear in $a$ and
$\{a,x,b\}_2$ is conjugate linear in $a,b,x$

\item  $\{ x,b,a \}_1 = \{ a,b,x \}_3$, and $\{ a,x,b \}_2  = \{
b,x,a \}_2$  for every $a,b\in E$ and $x\in X$.

 \item Denoting by
$\J ...$ any of the products $\{ .,.,. \}_1$, $\{ .,.,. \}_2$ and
$\{ .,.,. \}_3$, the identity $\J {a}{b}{\J cde} = \J{\J abc}de $
$- \J c{\J bad}e +\J cd{\J abe},$ holds whenever one of the
elements  $a,b,c,d,e$ is in $X$ and the rest are in $E$.
\end{enumerate}

It is a little bit  laborious to check that the dual space,
$E^*$, of a complex (resp., real) Jordan Banach triple $E$ is  a
complex (resp., real) triple $E$-module with respect to the
products: \begin{equation}\label{eq module product dual 1}   \J
ab{\varphi} (x) = \J {\varphi}ba (x) := \varphi \J bax
\end{equation} and \begin{equation}\label{eq module product dual
2} \J a{\varphi}b (x) := \overline{ \varphi \J axb },  \forall \varphi\in
E^*, a,b,x\in E. \end{equation}
For each submodule $S$ of a  triple $E$-module $X$, we define its
\emph{quadratic annihilator}, Ann$_{E} (S)$,  as the set 
$\{ a\in
E : Q (a) (S) = \J aSa = 0\}$.

\begin{theorem} Let $E$ be a complex
JB$^*$-triple, $X$ a Banach triple $E$-module, and let $\delta:
E\to X$ be a triple derivation. Then $\delta$ is continuous if and
only if $\hbox{Ann}_{E}(\sigma_{_X} (\delta))$ is a (norm-closed)
linear subspace of $E$ and $$\J {\hbox{Ann}_{_{E}} (\sigma_{_X}
(\delta))}{\hbox{Ann}_{_{E}} (\sigma_{_X} (\delta))}{\sigma_{_X}
(\delta)} = 0.$$ 
\end{theorem}

\begin{corollary}
 Let $E$ be a real or complex
JB$^*$-triple.   Then \\  
(a) Every derivation
$\delta : E \to E$ is continuous. {\rm (Barton-Friedman)}\\
(b) Every derivation $\delta :
E \to E^*$ is continuous {\rm (motivated ternary weak amenability)}.\smallskip

In particular,\\
(c) Automatic continuity of derivation of $JB^*$-algebra into a Jordan module {\rm (Hejazian-Niknam)}\\
(d) Automatic continuity of derivation of  $C^*$-algebra into a module {\rm (Ringrose)}\\
(e) Automatic continuity of Jordan derivation of $C^*$-algebra into a Jordan module {\rm (Hejazian-Niknam)}
\end{corollary}

 

\subsection{Triple derivations on von Neumann algebras}
As noted in subsection~\ref{5.1}, building on earlier work of Kadison \cite{Kadison66}, Sakai \cite{Sakai66} proved that
very derivation of a von Neumann algebra into itself is 
 inner.   
 Building on earlier work of Bunce and Paschke \cite{BunPas80}, Haagerup, on his way to proving that every C$^*$-algebra is weakly amenable,  showed in \cite{Haagerup83}
 that every derivation of a von Neumann algebra into its predual is inner.
 Thus the first  Hochschild cohomology groups $H^1(M,M)$ and $H^1(M,M_*)$ vanish for any von Neumann algebra $M$.

  It is also known that every triple derivation of a von Neumann algebra into itself is 
 an inner triple derivation (Theorem~\ref{4.19}), and every Jordan derivation of a von Neumann algebra into itself is an inner Jordan derivation (Theorem~\ref{4.14}).  
 Triple derivations and inner triple derivations into a triple module are defined and discussed in subsection~\ref{5.5}.\footnote{In this subsection, we use the terms `triple' and `ternary' interchangeably, while mindful that in some quarters `triple' means  `Jordan triple' and `ternary' refers to an associative triple setting, such as a TRO (ternary ring of operators)}  However, in this subsection we shall only be concerned with these concepts for von Neumann algebras, with their associated Jordan structure.

 Two consequences of \cite{HoPerRus11} are that  finite dimensional von Neumann algebras and abelian von Neumann algebras  have the  property that every triple derivation into the predual is an inner triple derivation, analogous to the Haagerup result.  This property was called  {\bf normal ternary weak amenability} in \cite{PluRus14} where it is shown that it rarely holds in a general von Neumann algebra, but that it comes close. 
 The main results of \cite{PluRus14} are the following two theorems. The first one gives a cohomological characterization of finite factors  and the second gives a ``zero-one'' law for factors.
 
\begin{theorem}\label{thm:0411131}
Let $M$ be a von Neumann algebra.
\begin{description}
\item[(a)] If every triple derivation of $M$ into $M_*$ is approximated in norm by inner triple derivations, then $M$ is finite.
\item[(b)] If $M$ is a finite von Neumann algebra acting on a separable Hilbert space or if $M$ is a finite factor, then every triple derivation of $M$ into $M_*$ is approximated in norm by inner triple derivations.
\end{description}
\end{theorem}

\begin{theorem}\label{thm:0416131}  
If $M$ is a properly infinite factor, then the real vector space of triple derivations of $M$ into $M_*$, modulo the norm closure of the inner triple derivations,  has dimension~1.
\end{theorem}

 Inner triple 
derivations  on a von Neumann algebra $M$  into its predual $M_*$ are closely related to the span of commutators of normal functionals with elements of $M$, denoted by $[M_*,M]$. 
 As noted earlier, for a finite factor of  type $I_n$, all triple derivations into the predual are inner triple derivations.
\ A consequence of \cite[Proposition 4.1]{PluRus14} and the  work of Dykema and Kalton \cite{DykKal05} is that no factor of type $II_1$ is normally ternary weakly amenable (\cite[Corollary 4.3]{PluRus14}, even though each triple derivation into the predual is a norm limit of inner triple derivations by Theorem~\ref{thm:0411131}.  

It is also shown in \cite{PluRus14} that a finite countably decomposable von Neumann algebra $M$
of type $I_n$ is normally ternary weakly amenable   if and only if every element of $M_*$ of central trace zero belongs to $[M_*,M]$. In particular, for such an algebra $M$ with a faithful normal finite trace,  if every element of $M_*$ of  trace zero belongs to $[M_*,M]$ then $M$ is normally ternary weakly amenable.

\subsection{Local derivations on $JB^*$-triples}

Linear maps which agree with a derivation at each point are called local derivations.  These have been studied in the Banach setting by
Kadison \cite{Kadison90}, Johnson \cite{Johnson00}, Ajupov et. al  \cite{AyuKudNur09},\cite{AyuKudNur10},\cite{AlbAyuKudNur11}, among others. Kadison proved that a continuous local derivation from a von Neumann algebra into a Banach module is a derivation.\smallskip

Let $X$ and $Y$ be Banach spaces. According to the terminology employed in the literature (see, for example, \cite{BattMol}), a subset $\mathcal{D}$ of the Banach space $B(X,Y)$, of all bounded linear operators from $X$ into $Y$, is called \emph{algebraically reflexive} in $B(X,Y)$ when it satisfies the property:
\begin{equation}\label{eq reflexivity} T\in B(X,Y) \hbox{ with } T(x)\in \mathcal{D} (x), \ \forall x\in X \Rightarrow T\in \mathcal{D}.
\end{equation}
\emph{Algebraic reflexivity} of ${\mathcal D}$ in the space $L(X,Y)$, of all linear mappings from $X$ into $Y$, a stronger version of the above property not requiring continuity of $T$, is defined by:
\begin{equation}\label{eq reflexivity2} T\in L(X,Y) \hbox{ with } T(x)\in \mathcal{D} (x), \ \forall x\in X \Rightarrow T\in \mathcal{D}.
\end{equation}

In 1990, Kadison proved that (\ref{eq reflexivity}) holds if ${\mathcal D}$ is the set Der$(M,X)$ of all (associative) derivations on a von Neumann algebra $M$ into a dual $M$-bimodule $X$ \cite{Kadison90}.   Johnson extended Kadison's result by establishing that the set ${\mathcal D}=$ Der$(A,X),$ of all (associative) derivations from a C$^*$-algebra $A$ into a Banach $A$-bimodule $X$ satisfies (\ref{eq reflexivity2}) \cite{Johnson00}.

Michael Mackey gave a talk on this topic  at the conference in honor of Cho-Ho Chu's 65th birthday in May 2012 in Hong Kong. He proved that every continuous local derivation on a $JBW^*$-triple is a derivation (see \cite{Mack}), and he suggested some problems, among them whether every local derivation on a $JB^*$-triple into itself or into a Banach module is automatically continuous, and if so, whether it is a derivation.  There are other problems in this area, some involving nonlinear maps. Many of these have now been answered.

 Algebraic reflexivity of the set of local triple derivations on a C$^*$-algebra and on a JB$^*$-triple have been studied in \cite{Mack, BurFerGarPe2014, BurFerPe2013} and \cite{FerMolPe}. More precisely, Mackey proves in \cite{Mack} that the set ${\mathcal D}=\hbox{Der}_{t} (M),$ of all triple derivations on a JBW$^*$-triple $M$ satisfies (\ref{eq reflexivity}). The result has been supplemented in \cite{BurFerPe2013}, where Burgos, Fernández-Polo and Peralta  prove that for each JB$^*$-triple $E$, the set ${\mathcal D}=\hbox{Der}_{t} (E)$ of all triple derivations on $E$ satisfies (\ref{eq reflexivity2}). In what follows, {\it algebraic reflexivity} will refer to the stronger version  (\ref{eq reflexivity2}) which does not assume the continuity of $T$.\smallskip

In \cite{BreSemrl95}, Bre\v{s}ar and \v{S}emrl proved that the set of all (algebra) automorphisms of $B(H)$ is algebraically reflexive  whenever $H$ is a separable, infinite-dimensional Hilbert space. Given a Banach space $X$. A linear mapping $T: X\to X$ satisfying the hypothesis at \ref{eq reflexivity2} for $\mathcal{D} = \hbox{Aut} (X)$, the set of automorphisms on $X$, is called a \emph{local automorphism}.
Larson and Sourour showed in \cite{LarSou} that for every infinite dimensional Banach space $X$, every surjective local automorphism $T$ on the Banach algebra $B(X),$  of all bounded linear operators on $X$, is an automorphism.
\smallskip

Motivated by the results of \v{S}emrl in \cite{Semrl97}, references witness a growing interest in a subtle version of algebraic reflexivity called \emph{algebraic 2-reflexivity} (cf. \cite{AyuKuday2012,AyuKuday2014, BurFerGarPe2014preprint, BurFerGarPe2014preprintTripleHom,KimKim05,LiuWong06,Mol2002,Mol2003} and \cite{Pop}). A subset $\mathcal{D}$ of the set $\mathcal{M}(X,Y) = Y^{X},$ of all mappings from $X$ into $Y$, is called \emph{algebraically 2-reflexive} when the following property holds: for each mapping $T$ in $\mathcal{M}(X,Y)$ such that for each $a,b\in X,$ there exists $S= S_{a,b}\in  \mathcal{D}$ (depending on $a$ and $b$), with $T(a) =S_{a,b} (a)$ and $T(b) = S_{a,b} (b)$, then $T$ lies in $\mathcal{D}.$ A mapping $T: X\to Y$ satisfying that for each $a,b\in X,$ there exists $S= S_{a,b}\in  \mathcal{D}$ (depending on $a$ and $b$), with $T(a) =S_{a,b} (a)$ and $T(b) = S_{a,b} (b)$ will be called a 2-local $\mathcal{D}$-mapping. \smallskip

\v{S}emrl establishes in \cite{Semrl97} that for every infinite-dimensional separable Hilbert space $H$, the sets Aut$(B(H))$ and Der$(B(H))$, of all  (algebra) automorphisms and associative derivations on $B(H)$, respectively, are algebraically 2-reflexive in $\mathcal{M}(B(H)) = \mathcal{M}(B(H),B(H)).$ Ayupov and Kudaybergenov proved in \cite{AyuKuday2012} that the same statement remains true for general Hilbert spaces (see \cite{KimKim04} for the finite dimensional case). Actually, the set Hom$(A)$, of all homomorphisms on a general C$^*$-algebra $A,$ is algebraically 2-reflexive in the Banach algebra $B(A)$, of all bounded linear operators on $A$, and the set $^*$-Hom$(A)$, of all $^*$-homomorphisms on $A,$ is algebraically 2-reflexive in the space $L(A)$, of all linear operators on $A$ (cf. \cite{Pe2014}).\smallskip

In recent contributions,  Burgos, Fernández-Polo and Peralta prove that the set Hom$(M)$ (respectively, Hom$_{t} (M)$), of all homomorphisms (respectively, triple homomorphisms) on a von Neumann algebra (respectively, on a JBW$^*$-triple) $M$, is an algebraically 2-reflexive subset of $\mathcal{M}(M)$ (cf. \cite{BurFerGarPe2014preprint}, \cite{BurFerGarPe2014preprintTripleHom}, respectively), while Ayupov and Kudaybergenov establish that set Der$(M)$ of all derivations on $M$ is algebraically 2-reflexive in $\mathcal{M}(M)$ (see \cite{AyuKuday2014}).\smallskip

In \cite{KudOikPerRus14} the set Der$_{t} (A)$ of all triple derivations on a C$^*$-algebra $A$ is considered. We recall that every C$^*$-algebra $A$ can be equipped with a ternary product of the form $$\{a,b,c\} = \frac12 (a b^* c + c b^* a). $$ When $A$ is equipped with this product it becomes a JB$^*$-triple in the sense of \cite{Kaup83}. A linear mapping $\delta: A\to A$ is said to be a \emph{triple derivation} when it satisfies the (triple) Leibnitz rule: $$\delta\{a,b,c\} = \{\delta(a),b,c\} + \{a,\delta(b),c\}+ \{a,b,\delta(c)\}.$$
According to the standard notation, 2-local Der$_{t} (A)$-mappings from $A$ into $A$ are called \emph{2-local triple derivations}. \smallskip

The main result of \cite{KudOikPerRus14} is the following theorem.

\begin{theorem}
 Every {\rm(}not necessarily linear nor continuous{\rm)} 2-local triple derivation on an arbitrary von Neumann algebra $M$ is a triple derivation {\rm(}hence linear and continuous{\rm)}, equivalently, Der$_{t} (M)$ is algebraically 2-reflexive in $\mathcal{M}(M)$.
\end{theorem}

\part{Cohomology}

\section{Cohomology of finite dimensional algebras}\label{sec:cohomology}

\subsection{Associative algebras}

The starting point for the cohomology theory of algebras is the paper of Hochschild from 1945 \cite{Hochschild45}. The standard reference of the theory is \cite{CarEil56}.  Two other useful references are due to Weibel (\cite{Weibel94},\cite{Weibel99}). We review here the definitions of the cohomology groups $H^n(M,X)$ and the interpretation of them in the cases $n=1,2$.

Let $M$ be an associative algebra and $X$ a two-sided $M$-module.
For $n\ge 1$, let
$L^n(M,X)$ be the linear space of all  $n$-linear maps from $M$ to $X$  ($L^0(M,X)=X$).  The
coboundary operator is the linear mapping  $\partial:L^n\rightarrow L^{n+1}$ (for $n\ge 1$) defined by
$$\partial\phi(a_1,\cdots,a_{n+1})=a_1\phi(a_2,\cdots,a_{n+1})
+\sum(-1)^j\phi(a_1,\cdots,a_{j-1},a_ja_{j+1},\cdots,a_{n+1})
+(-1)^{n+1}\phi(a_1,\cdots,a_n)a_{n+1}$$

For $n=0$,
$\partial:X\rightarrow L(M,X)$ is defined by $\partial x(a)=ax-xa$.
Since $\partial\circ\partial=0$,
 $\hbox{Im}(\partial:L^{n-1}\rightarrow L^n)\subset \hbox{ker}(\partial:L^n\rightarrow L^{n+1})$, and 
$H^n(M,X)=\hbox{ker}\partial/\hbox{Im}\partial$ is a vector space.

For $n=1$,
$\hbox{ker}\, \partial=$
$\{\phi:M\rightarrow X: a_1\phi(a_2)-\phi(a_1a_2)+\phi(a_1)a_2=0\}$
which is  the space of derivations from $M$ to $X$; and 
$\partial:X\rightarrow L(M,X)$ is given by $\partial x(a)=ax-xa$ so that 
$\hbox{Im}\, \partial$ is the space of inner derivations
\smallskip

Thus  $H^1(M,X)$ measures how close derivations are to inner derivations.
\smallskip

An associative algebra $B$ is an extension of associative algebra $A$ if there is a homomorphism $\sigma$ of $B$ onto $A$.  The extension splits if $B=\ker \sigma\oplus A^*$ where $A^*$ is an algebra isomorphic to $A$, and is singular if $(\ker\sigma)^2=0$.

\begin{proposition}
There is a one to one correspondence between isomorphism  classes of singular extensions of $A$ and $H^2(A,A)$
\end{proposition}

\subsection{Lie algebras}
 
Shortly after the introduction of the cohomology for associative algebras was introduced, there appeared a corresponding theory for Lie algebras.   We follow \cite{Jacobson62} for the definitions and initial results. Applications can be found in \cite{Fuks86} and \cite{Knapp88}.

If $L$ is a Lie algebra, then an $L$-module is a vector space $M$ and a mapping of $M\times L$ into $M$, $(m,x)\mapsto mx$, satisfying
$(m_1+m_2)x=m_1x+m_2x$,
$\alpha(mx)=(\alpha m)x=m(\alpha x)$, and
$m[x_1,x_2]=(mx_1)x_2-(mx_2)x_1$.
\smallskip

Let $L$ be a Lie algebra, $M$ an $L$-module.  If $i\ge 1$,  an {\it $i$-dimensional $M$-cochain} for $L$ is a skew symmetric $i$-linear mapping $f$ of $L\times L\times\cdots\times L$ into $M$. Skew symmetric means that if two arguments in $f(x_1,\cdots,x_i)$ are interchanged, the value of $f$ changes sign. A 0-dimensional cochain is a constant function from $L$ to $M$. \smallskip

The coboundary operator $\delta$ (for $i\ge 1$) is:
\[
\delta(f)(x_1,\cdots,x_{i+1})
=\sum_{q=1}^{i+1}(-1)^{i+1}f(x_1,\cdots,\hat x_q,\cdots,x_{i+1})x_q
+
\sum_{q<r=1}^{i+1}(-1)^{r+q}f(x_1,\cdots,\hat x_q,\cdots,\hat x_r,\cdots,x_{i+1},[x_q,x_r]).
\]
and for $i=0$, $\delta(f)(x)=ux $ (module action), if $f$ is the constant $u\in M$.\smallskip

One verifies that $\delta^2=0$ giving rise to cohomology groups 
\[
H^i(L,M)=Z^i(L,M)/B^i(L,M)
\]
If $i=0$ we take $B^i=0$ and $H^0(L,M)=Z^0(L,M)=\{u\in M:ux=0, \forall x\in L\}$.

\begin{theorem}[Whitehead lemmas]
If $L$ is a finite dimensional semisimple Lie algebra over a field of characteristic 0, then
$$H^1(L,M)=H^2(L,M)=0$$
 for every finite dimensional module $M$ of $L$.
 \end{theorem}

\begin{theorem}[Whitehead]
If $L$ is a finite dimensional semisimple Lie algebra over a field of characteristic 0, then
$$H^i(L,M)=0\ (\forall i\ge 0)$$ 
 for every finite dimensional \underline{irreducible} module $M$ of $L$ such that $ML\ne 0$.
\end{theorem}

\subsection{Jordan algebras}\label{sub:hint}

The cohomology theory for Jordan algebras is less well developed than for associative and Lie algebras.  A starting point would seem to be the papers of Gerstenhaber in 1964 \cite{Gerstenhaber64} and Glassman in 1970 \cite{Glassman70PJM}, which concern arbitrary 
nonassociative algebras.  A study focussed primarily on Jordan algebras is \cite{Glassman70JAlg}.

Let $A$ be an algebra defined by a set of identities and let $M$ be an $A$-module. A singular extension of length 2 is, by definition, a null extension of $A$ by $M$.
A null extension is a short exact sequence
\[
0\rightarrow M\stackrel{\alpha}{\rightarrow}E\stackrel{\beta}{\rightarrow} A\rightarrow 0
\]
where, provisionally, $M$ is an algebra (rather than an $A$-module) with $M^2=0$.
 If $n>2$, a singular extension of length $n$ is an exact sequence of bimodules
\[
0\rightarrow M\rightarrow M_{n-1}\rightarrow \cdots\rightarrow M_2\rightarrow E\rightarrow A\rightarrow 0
\]

Morphisms, equivalences, addition, and scalar multiplication of equivalence classes of singular extensions can be defined.  
Then for $n\ge 2$,
$H^n(A,M):=$ equivalence classes of singular extensions of length $n$.
These definitions are equivalent to the classical ones in the associative and Lie cases.\smallskip

It is shown in \cite{Gerstenhaber64} (using generalized projective resolutions) showed that if $0\rightarrow M'\rightarrow M\rightarrow M''\rightarrow 0$ is an exact sequence of $A$-bimodules, then there are natural homomorphisms $\delta^n$ so that the long sequence
\[
0\rightarrow \hbox{Der}\, (A,M')\rightarrow \hbox{Der}\, (A,M)\rightarrow \hbox{Der}\, (A,M'')\stackrel{\delta^1}{\rightarrow}
\]
\[
 H^2(A,M')\rightarrow H^2(A,M)\rightarrow H^2(A,M'')
\stackrel{\delta^2}{\rightarrow}
\]
\[ H^3(A,M')\rightarrow H^3(A,M)\rightarrow H^3(A,M'')\rightarrow
\]
\[
\cdots\rightarrow H^n(A,M'')\stackrel{\delta^n}{\rightarrow} H^{n+1}(A,M')\rightarrow\cdots
\]
is exact.  In particular,
\[
H^n(A,M)=\ker \delta^n/\hbox{im}\, \delta^{n-1}\quad (n\ge 2)
\]

What about $H^0(A,M)$ and $H^1(A,M)$ and Jordan algebras?
  For this we turn first to Glassman's thesis of 1968, embodied in \cite{Glassman70JAlg} and \cite{Glassman70PJM}.
Given an algebra $A$,
consider the functor $\mathcal T$ from the category $\mathcal C$ of $A$-bimodules and $A$-homomorphisms to the category $\mathcal V$ of vector spaces and linear maps:

\[
M\in {\mathcal C}\mapsto \hbox{Der}\, (A,M)\in{\mathcal V}\quad, \quad \eta\mapsto \tilde\eta
\]
where $\eta\in\hbox{Hom}_A(M_1,M_2)$ and
 $\tilde\eta\in \hbox{Hom}(\hbox{Der}\, (A,M_1),\hbox{Der}\, (A,M_2))$ is given by
  $\tilde\eta=\eta\circ D$. 
$\left(A\stackrel{D}{\rightarrow}M_1\stackrel{\eta}{\rightarrow}M_2\right)$
\smallskip

An \underline{inner derivation functor} is a subfunctor $\mathcal J$ which respects epimorphisms, that is, 
\[
M\in {\mathcal C}\mapsto {\mathcal J}(A,M)\subset \hbox{Der}\, (A,M)\in{\mathcal V}
\]
and if  $\eta\in\hbox{Hom}_A(M_1,M_2)$ is onto, then so is ${\mathcal J}(A,\eta):=\tilde\eta|{\mathcal J}(A,M)$
\smallskip

Relative to the choice of $\mathcal J$ one defines
$
H^1(A,M)=\hbox{Der}\, (A,M)/{\mathcal J}(A,M).
$
(The definition of $H^0(A,M)$ is more involved so is omitted here).
Glassman then proves that
\[
0\rightarrow H^0(A,M')\rightarrow H^0 (A,M)\rightarrow H^0(A,M'')\stackrel{\delta^0}{\rightarrow}
\]
\[
 H^1(A,M')\rightarrow  H^1(A,M)\rightarrow  H^1(A,M'')
\stackrel{\delta^1}{\rightarrow}
\]
\[
 H^2(A,M')\rightarrow H^2(A,M)\rightarrow H^2(A,M'')
\stackrel{\delta^2}{\rightarrow}
\]
\[ H^3(A,M')\rightarrow H^3(A,M)\rightarrow H^3(A,M'')\rightarrow
\]
\[
\cdots\rightarrow H^n(A,M'')\stackrel{\delta^n}{\rightarrow} H^{n+1}(A,M')\rightarrow\cdots
\]
is exact. 

 We mention just one other result from \cite{Glassman70JAlg}.
In a Jordan algebra, recall that $\{xby\}=(xb)y+(by)x-(xy)b$. 
 The $b$-homotope of $J$, written $J^{(b)}$ is the Jordan algebra structure on the vector space $J$ given by the multiplication $x\cdot_b y=\{xby\}$. 
  If $M$ is a bimodule for $J$ and $b$ is invertible, the corresponding bimodule $M^{(b)}$ for $J^{(b)}$ is the vector space $M$ with action
$a\cdot_b m=m\cdot_b a=\{abm\}$.

\begin{lemma}
$H^n(J,M)\sim H^n(J^{(b)},M^{(b)})$
\end{lemma}

We next record the Jordan analogs of the first and second Whitehead lemmas as described in 
\cite{Jacobson57}.  The first one has already been mentioned in Theorem~\ref{thm:0208121}.

\begin{theorem}[Jordan analog of first Whitehead lemma \cite{Jacobson51}]
Let $J$ be a finite dimensional semisimple Jordan algebra over a field of characteristic 0 and let $M$ be a $J$-module.  Let $f$ be a linear mapping of $J$ into $M$ such that 
\[
f(ab)=f(a)b+af(b).
\]
Then there exist $v_i\in M, b_i\in J$ such that
\[
f(a)=\sum_i((v_ia)b-v_i(ab_i)).
\]
\end{theorem}

\begin{theorem}[Jordan analog of second Whitehead lemma \cite{Penico51}]\label{5.6}
Let $J$ be a finite dimensional  separable\footnote{See the footnote to Theorem~\ref{thm:0131141}}  Jordan algebra  and let $M$ be a $J$-module.  Let $f$ be a bilinear mapping of $J\times J$ into $M$ such that 
\[
f(a,b)=f(b,a)
\]
and
\[
f(a^2,ab)+f(a,b)a^2+f(a,a)ab
=f(a^2b,a)+f(a^2,b)a+(f(a,a)b)a
\]
Then there exist a linear mapping $g$ from $J$ into $M$ such that
\[
f(a,b)=g(ab)-g(b)a-g(a)b
\]
\end{theorem}

Two proofs of Theorem~\ref{5.6} are given in \cite{Jacobson68}.  One of them, which uses the classification of finite dimensional Jordan algebras, is outlined in the next subsection (see Corollary~\ref{5.23}).   The other proof,  uses Lie algebras and is contained in \cite[pp. 324--336]{Jacobson68} and will not be discussed here.

 \subsubsection{Jordan classification approach}\label{5.3.1}

In this subsection, as just indicated, we shall give a detailed outline of one proof of Theorem~\ref{5.6}, following \cite{Jacobson68}. This will serve as a model, or road map, for a proposed infinite dimensional generalization.

\begin{center} {\bf Step 1} Extensions of algebras and factor sets (\cite[pp.91-92]{Jacobson68})\end{center}

Let $A$ and $M$ be 
 algebras  (associative, Jordan, Lie, alternative). 
  An {\bf extension of $A$ by $M$} is a short exact sequence
$
0\rightarrow M\stackrel{\alpha}{\rightarrow} E\stackrel{\beta}{\rightarrow} A\rightarrow 0.
$
Thus $\alpha$ is an injective homomorphism and $\beta$ is a surjective homorphism.
Extensions 
$
0\rightarrow M\stackrel{\alpha}{\rightarrow} E\stackrel{\beta}{\rightarrow} A\rightarrow 0.
$
and 
$
0\rightarrow M\stackrel{\alpha'}{\rightarrow} E'\stackrel{\beta'}{\rightarrow} A\rightarrow 0.
$
are {\bf equivalent} if there exists a homomorphism $\gamma:E\rightarrow E'$ such that $\alpha'=\gamma\circ \alpha$ and $\beta=\beta'\circ\gamma$.   Thus $\gamma$ is an isomorphism of $E$ onto $E'$.\smallskip

An extension 
$
0\rightarrow M\stackrel{\alpha}{\rightarrow} E\stackrel{\beta}{\rightarrow} A\rightarrow 0.
$
is {\bf split} (or {\bf inessential}) is there exists a homomorphism $\delta:A\rightarrow E$ with $\beta\circ \delta=1_A$.  Thus $E=\delta(A)\oplus \alpha(M)$ as vector spaces and $\delta(A)$ is a subalgebra of $E$ isomorphic to $A$.
An extension 
$
0\rightarrow M\stackrel{\alpha}{\rightarrow} E\stackrel{\beta}{\rightarrow} A\rightarrow 0.
$
is {\bf null} if $M^2=0$, that is, all products in $M$ are zero. 

Let 
$
0\rightarrow M\stackrel{\alpha}{\rightarrow} E\stackrel{\beta}{\rightarrow} A\rightarrow 0.
$
be any extension and identify $M$ with $\alpha(M)\subset E$.   Then we may write $E=M\oplus \delta(A)$, a vector space direct sum, for some linear map $\delta:A\rightarrow E$ such that $\beta\circ \delta=1_A$.   $E$ is an $E$-module and since $M$ is an ideal in $E$, it is a submodule.  If $M^2=0$ then 
 $M$ is an $E/M$-module via $(e+M)\cdot u=eu,\ u\cdot(e+M)=ue$.
Now  the isomorphism $\beta(e)\rightarrow e+M$ of $A$ with $E/M$ completes the proof of the following lemma.

\begin{lemma}\label{exercise:1219131}
$M$ is a 
 $A$-module under the module actions $a\cdot u=\delta(a)u$, $u\cdot a=u\delta(a)$ for $a\in A$ and $u\in M\subset E$.  {\rm (Multiplication in $E$)} 
\end{lemma}

Although the next theorem is proved in \cite{Jacobson68} simultaneously for associative, Lie. and Jordan algebras, we shall now restrict to Jordan algebras.


\begin{definition}
Let $M$ be a Jordan $A$-module.  A bilinear map $h:A\times A\rightarrow M$ is a  (Jordan) {\bf 2-cocycle} if it is symmetric and satisfies
\[
(h(a,a)\cdot b)\cdot a+h(a^2,b)\cdot a+h(a^2b,a)=a^2\cdot h(b,a)+h(a,a)\cdot (ba)+h(a^2,ba).
\] 
A  map $h:A\times A\rightarrow M$ is a  (Jordan) {\bf 2-coboundary} if it is of the form
\[
h(a,b)=\mu(ab)-a\cdot \mu(b)-\mu(a)\cdot b
\] 
for some linear map $\mu:A\rightarrow M$.
The vector space of all 2-cocycles modulo 2-coboundaries is denoted $H^2(A,M)$.
\end{definition}

\begin{theorem}[THEOREM 12, p.94]\label{thm:0207141}
Let $M$ be a Jordan $A$-module.  Then there is a bijection of $H^2(A,M)$ and the set of equivalence classes of null extensions of $A$ by $M$ such that the associated bimodule structure on $M$ given by Lemma~\ref{exercise:1219131} is the given one.   In this correspondence the equivalence class of 0 in $H^2(A,M)$ corresponds to the isomorphism class of inessential extensions.
\end{theorem}

\begin{center} {\bf Step 2} Peirce decomposition; connected idempotents (\cite[pp.117--124]{Jacobson68})\end{center}

Let $e$ be an idempotent
in a Jordan algebra $J$. The {\bf Peirce decomposition} of $J$ with respect to $e$ is 
\[
J=J_0(e)\oplus J_1(e)\oplus J_{\frac{1}{2}}(e),
\]
where $J_i(e)=\{x_i\in J:x_ie=ix_i\}$.

\begin{lemma}[Lemma 1, p.119]
If $e$ is an idempotent, then
\begin{itemize}
\item $J_0^2\subset J_0$, $J_1^2\subset J_1$, $J_0J_1=0$
\item $J_{\frac{1}{2}}^2\subset J_0+J_1$, $J_{\frac{1}{2}}(J_0+J_1)\subset J_{\frac{1}{2}}$
\end{itemize}
\end{lemma}

Suppose $J$ has an identity element and 
$1=e_1+\ldots+e_n$ where $e_i^2=e_i$ and $e_ie_j=0$ ($i\ne j$). The {\bf Peirce decomposition} of $J$ with respect to the idempotents $e_i$ (orthogonal with sum 1)  is 
\[
J=\sum_{i\le j}\oplus J_{ij},
\]
where $$J_{ii}=\{x_{ii}\in J:x_{ii}e_i=x_{ii}\}=J_1(e_i)$$ and $$J_{ij}=J_{\frac{1}{2}}(e_i)\cap 
J_{\frac{1}{2}}(e_j)=\{x_{ij}:x_{ij}e_i=\frac{1}{2}x_{ij}=x_{ij}e_j\}.$$

\begin{lemma}[Lemma 2, p.120]
If $e_i$ are orthogonal idempotents with sum 1,  then
\begin{itemize}
\item $J_{ii}^2\subset J_{ii}$,\quad  $J_{ij}J_{ii}\subset J_{ij}$, \quad $J_{ij}^2\subset J_{ii}+J_{jj}$, \quad $J_{ii}J_{jj}=0\ (i\ne j)$
\item $J_{ij}J_{jk}\subset J_{ik}$, \quad $J_{ij}J_{kk}=0$, \quad $J_{ij}J_{kl}=0$ if $i,j,k,l$ are distinct.
\end{itemize}
\end{lemma}

Let $e_1$ and $e_2$ be orthogonal idempotents so that we have $$J_1(e_1+e_2)=J_{11}+J_{22}+J_{12}=J_1(e_1)+J_1(e_2)+J_{12}.$$
The orthogonal idempotents $e_1$ and $e_2$ are {\bf connected} if there exists $u_{12}\in J_{12}$ which is invertible in $J_1(e_1+e_2)$.  They are {\bf strongly connected} if here exists $u_{12}\in J_{12}$ with $u_{12}^2=e_1+e_2$.

\begin{center} {\bf Step 3} Jordan Matrix Algebras  (\cite[pp.125--131]{Jacobson68})\end{center}

Let $D$ be a unital algebra with involution $j(d)=\overline{d}$. Then $D_n$ denotes the algebra of  $n$ by $n$ matrices with entries from $D$, with involution $X\mapsto X^J=\overline{X}^t$, and $H(D_n)\subset D_n^+$ denotes the subalgebra of symmetric elements.  Algebras of the form $H(D_n)$ which are Jordan algebras are called {\bf Jordan matrix algebras}.\smallskip

Let $e_{ij}$ be the usual matrix units in $D_n$ and if $x\in D$, we identify $x$ with the diagonal matrix  in $D_n$ all of whose diagonal entries are $x$.  Then $xe_{ij}$ is the matrix with $x$ in the $(i,j)$-position and zeros elsewhere.  We set
\[
x[ij]=xe_{ij}+(xe_{ij})^J.
\]

Theorems~\ref{thm:0201141} and \ref{thm:0201142} are not used in  this paper, but are stated here for their intrinsic interest.

\begin{theorem}[THEOREM 1, p.127]\label{thm:0201141}  For $n\ge 3$,
$H(D_n)$ is a Jordan algebra if and only if $D$ is associative, or $n=3$ and $D$ is alternative with symmetric elements in the nucleus.
\end{theorem}

\begin{theorem}[THEOREM 2, p.129]\label{thm:0201142}  For $n\ge 3$, let $H=H(D_n)$ be a Jordan matrix algebra.
There are one to one correspondences
\[
\{\hbox{*-subalgebras } E \hbox{ of } D\}\leftrightarrow \{\hbox{*-subalgebras } H\cap E_n\hbox{ containing }1[ij]\}
\]
and 
\[
\{\hbox{*-ideals } E \hbox{ of } D\}\leftrightarrow \{\hbox{*-ideals } H\cap E_n\}
\]
In each case, the map $E\mapsto H\cap E_n$ is a lattice isomorphism, and in the second case, $(H\cap E_n)^2=0$ if and only if $E^2=0$,

\end{theorem}

\begin{theorem}[THEOREM 3, p.130]  For $n\ge 3$, let $H(D_n)$ and $H(E_n)$  be  Jordan matrix algebras, and suppose $\eta:D\rightarrow E$ is a  unital *-homomorphism (into).
Then the restriction $\sigma$ to $H(D_n)$ of the mapping $(d_{ij})\mapsto (\eta(d_{ij}))$ of $D_n$ is a *-homomorphism of $H(D_n)$ into $H(E_n)$ such that $\sigma(1[ij]])=1[ij]$.

Conversely, if $\sigma$ is a *-homomorphism of $H(D_n)$ into $H(E_n)$ such that $\sigma(1[ij]])=1[ij]$, then there exists a unital *-homomorphism $\eta$ of $D$ into $E$ such that $\sigma$ is the restriction of the mapping $(d_{ij})\mapsto (\eta(d_{ij}))$.
\end{theorem}

\begin{center} {\bf Step 4} Coordinatization Theorems  (\cite[pp.132--138]{Jacobson68})\end{center}

\begin{remark}
If $J=H(D_n)$ and we write $x[ij]=xe_{ij}+(xe_{ij})^{J}$, $e_i=\frac{1}{2}[ii]=e_{ii}$, then the $e_i$ are orthogonal idempotents with sum 1, and $e_1$ and $e_j$ are strongly connected.  
\end{remark}

\begin{theorem}[THEOREM 5, p.133; strong coordinatization theorem]  
Let $J$ be a Jordan algebra whose identity is the sum of $n\ge 3$ strongly connected idempotents.
Then $J$ is isomorphic to a Jordan matrix algebra $H(D_n)$.

More precisely, if $1=\sum_1^n e_i$ and we have $u_{1j}$ with $e_1u_{1j}=\frac{1}{2}u_{1j}=e_ju_{1j}$ and $u_{1j}^2=e_1+e_j$, then there is a Jordan matrix algebra $H(D_n)$ and an isomorphism $\eta$ of $J$ onto $H(D_n)$ such that $\eta(e_i)=\frac{1}{2}[ii]$ and $\eta(u_{1j})=1[1j]$.
\end{theorem}

\begin{center} {\bf Step 5} Lifting of idempotents (\cite[pp.148--151]{Jacobson68})\end{center}

An algebra is {\bf nil} if every element is nilpotent.

\begin{lemma}[Lemma 5, p.150]
Let $J$ be a Jordan algebra with 1, $N$ a nil ideal, $e_1$ and $e_2$
nonzero orthogonal idempotents in $J$.   Then $e_1$ and $e_2$ are connected (strongly connected) if and only if $\overline{e}_1=e_1+N$ and $\overline{e}_2=e_2+N$ are connected (strongly connected) in $\overline{J}=J/N$.
\end{lemma}

\begin{theorem}[THEOREM 10, p.151]  \label{thm:0208141}
Let $J$ be a Jordan algebra with 1, $N$ a nil ideal such that $\overline{J}=J/N$ is isomorphic to a Jordan matrix algebra $H(\overline{D}_n)$ of order $n\ge 3$.  Then $J$ is isomorphic to a Jordan matrix algebra $H(D_n)$ where the ideal in $H(D_n)$ corresponding to $N$ has the form 
$M_n\cap H(D_n)$ where $M$ is an ideal in $D$ and $\overline{D}$ is isomorphic to $D/M$ as algebras with involution.
\end{theorem}

\begin{center} {\bf Step 6} Solvable ideals and the radical (\cite[pp.190--196]{Jacobson68})\end{center}

A subspace $B$ if a Jordan algebra $J$ is an {\bf associator ideal} if any associator $[a_1,a_2,a_3]:=(a_1a_2)a_3-a_1(a_2a_3)\in B$ whenever any one of $a_i\in B$.

\begin{theorem}[THEOREM 1, p.191]  
Let $J$ be a Jordan algebra, $B$ an ideal of finite codimension $m$.  Define 
\begin{enumerate}
\item $B^{[0]}=B,\quad B^{[k]}=B^{[k-1]}B+(B^{[k-1]}B)J$
\item $B^{(0)}=B,\quad B^{(k)}=B^{(k-1)}B^{(k-1)}+(B^{(k-1)}B^{(k-1)})J$.
\end{enumerate}
Then $B^{[k]}$ and $B^{(k)}$ are ideals  and $B^{[m+1]}\cup B^{(m+1)}\subset B^2$.
\end{theorem}

The powers $J^{2^k}$ of a Jordan algebra $J$ are $J^{2^0}=J$ and $J^{2^k}=(J^{2^{k-1}})^2$.
$J^2$ is an ideal and $J^{2^k}$ is a subalgebra.  $J$ is {\bf solvable } if there exists an integer $N$ with $J^{2^N}=0$.  A solvable Jordan algebra is nil.
If $J$ is finite dimensional (more generally, if $J$ satisfies the maximum condition for ideals), then
it contains a solvable ideal $R$, the {\bf radical} ($=\hbox{rad}\, J$), such that $R$ contains every solvable ideal of $J$. $J$ is {\bf semisimple} if $R=0$ ($J$ has no nonzero solvable ideals).  
Note that $J/\hbox{rad}\, J$ is semisimple.

The next two theorems, due to Albert,  are included for their intrinsic interest.
An algebra is {\bf nilpotent} if there exists an integer $N$ such that every product (in any association) of $N$ (distinct) elements is zero.   A nilpotent Jordan algebra is solvable.
\begin{theorem}[COROLLARY 1, p.195]  
A finite dimensional solvable Jordan algebra is nilpotent.
\end{theorem}

\begin{theorem}[THEOREM 3, p.196]  
Any finite dimensional nil algebra is solvable.
\end{theorem}

\begin{center} {\bf Step 7} The theorem of Albert-Penico-Taft (\cite[pp.287--292]{Jacobson68})\end{center}

\begin{theorem}[Theorem 13, p. 292]\label{thm:0131141}
Let $E$ be a finite dimensional Jordan algebra, $M$ an ideal in $E$ such that $E/M$ is separable\footnote{Separable, in this context, means that the algebra remains semisimple with respect to all extensions of the ground field.  For algebraically closed fields, this is the same as being semisimple}. Then there is a subalgebra $R$ of $E$ such that $E=R\oplus M$ as vector spaces. In other words, every finite dimensional extension of a separable Jordan algebra splits.
\end{theorem}

By Theorem~\ref{thm:0207141}, we have vanishing of the second cohomology group and a proof of Theorem~\ref{5.6}.
\begin{corollary}[Corollary, p 292]\label{5.23}
If  $M$ is a finite dimensional bimodule for $J$, a separable Jordan algebra,  then  $H^2(J,M)=0$.
\end{corollary}

\begin{description}
\item[Reduction I, p.288]  If Theorem \ref{thm:0131141} is true for solvable $M$, then it is true for arbitrary $M$.
\item[Reduction II, p.288] If Theorem \ref{thm:0131141} is true for solvable $M$ with $M^2=0$, then it is true for arbitrary solvable $M$.
\item[Reduction III, p.288] If Theorem \ref{thm:0131141} is true for algebras $E$ with a unit, then it is true for arbitrary $E$.
\item[Reduction IV, p.288-289] If Theorem \ref{thm:0131141} is true for algebras $E$ over an algebraically closed field, then it is true for arbitrary $E$

\item[Reduction V, p.289-290] If Theorem \ref{thm:0131141} is true for algebras $E$ with $E/M$ simple, then it is true for arbitrary $E$
\end{description}

The proof is now completed by appealing to the classification of unital simple Jordan algebras over an algebraically closed field in order to prove

\begin{proposition}[Lemmas 2-4, pp. 290--291]
Let $E$ be  a finite dimensional Jordan algebra with unit over an algebraically closed field and let $M$ be an ideal in $E$ such that $M^2=0$ and $E/M$ is separable. Then 
$E$ contains a subalgebra isomorphic to $E/M$.
\end{proposition}

For later use, we indicate the proof of Reduction IV when the underlying field is the real numbers.
If   we complexify the short exact sequence
\begin{equation}\label{eq:0207143}
0\rightarrow M\rightarrow E\rightarrow J=E/M\rightarrow 0,
\end{equation}
 and this complexified sequence splits, then so does the original one.
To show this,  one uses the equivalence given by Theorem~\ref{thm:0207141}. as follows

Let $h:J\times J\rightarrow M$ be a Jordan 2-cocycle (over the real field).  We need to show that there is a linear map $\mu:J\rightarrow M$ such that 
\begin{equation}\label{eq:0207141}
h(a,b)=\mu(ab)-\mu(a)\cdot b-\mu(b)\cdot a.
\end{equation}
Let $u_1,\ldots,u_n$ be a basis for $J$ over ${\bf R}$ and $v_1,\ldots, v_r$ a basis for $M$ over ${\bf R}$.  Then (\ref{eq:0207141}) holds if and only if it holds for $a,b\in
\{u_1,\ldots,u_n\}$.  With $\mu$ provisionally defined by $\mu(u_i)=\sum_p \mu_{ip}v_p$, define the quantities $  \eta_{ijq},\ \    \gamma_{ijk},\ \  \delta_{piq}$ by the formulas
 $$h(u_i,u_j)=\sum_q\eta_{ijq}v_q,\ \
u_iu_j=\sum_k\gamma_{ijk}u_k;\ \ 
v_p\cdot u_i=\sum_q\delta_{piq}v_q.$$
Then  (\ref{eq:0207141}) for $a=u_i$ and $b=u_j$ is equivalent to the set of linear equations
\begin{equation}\label{eq:0207142}
\eta_{ijk}=\sum_k\gamma_{ijk}\mu_{kq}-\sum_p\mu_{ip}\delta_{pjq}-\sum_p\mu_{jp}\delta_{piq}.
\end{equation}
for the $\mu$'s in ${\bf R}$.

The solvability  of  (\ref{eq:0207142}) for the  $\mu$'s in ${\bf R}$ is a necessary and sufficient condition that the extension splits.  If we complexity the sequence (\ref{eq:0207143}), and extend the maps $\alpha,\beta$, and $h$, we still have $(M^{{\bf C}})^2=0$, $J^{{\bf C}}$ is separable and $h^{{\bf C}}$ is a 2-cocycle.  The bases $\{u_i\}$ and $\{v_j\}$ remain bases over the complex field.  Assuming the theorem holds in the algebraically closed case,  the equations 
(\ref{eq:0207142}) have a solution $\{\mu_{ip}\}$ in ${\bf C}$.  Since these are linear equations with coefficients in ${\bf R}$, it follows that they have a solution in ${\bf R}$.\smallskip

  \subsubsection{Cohomology groups of Jordan algebras and triples, after McCrimmon}
\newcommand{\cm}{\hbox{$\mathcal M$}}

A study of low dimensional cohomology for quadratic Jordan algebras is given in \cite{McCrimmon71}.  Since quadratic Jordan algebras (which we do not define here as they coincide with ``linear'' Jordan algebras over characteristic 0 fields) can be considered a bridge from Jordan algebras to Jordan triple systems, this would seem to be a good place to look for exploring cohomology theory for Jordan triples.   Indeed, this is hinted at in \cite{McCrimmon82}.  Let us review some aspects of \cite{McCrimmon71,McCrimmon82}.  Although \cite{McCrimmon71} is about Jordan algebras, the concepts are phrased in terms of the associated triple product $\{abc\}=(ab)c+(cb)a-(ac)b$.

Let $J$ be a unital (quadratic) Jordan algebra and $\cm$ a unital $J$-module.  A linear map $d:J\rightarrow \cm$ is a {\it derivation} if $d(1)=0$ and $d\tp{x}{y}{x}=\tpc{x}{d(y)}{x}+2\tpc{d(x)}{y}{x}$.   An {\it inner derivation} is a sum of derivations of the form $d_{x,m}(y)=\tp{x}{m}{y}-\tp{m}{x}{y}$.  
We shall use the notation  $\hbox{Inder}\, (J,\cm)\subset \hbox{Der}\, (J,\cm)$.  $H^1(J,\cm)= \hbox{Der}\, (J,\cm)/\hbox{Inder}\, (J,\cm)$.

We can express $H^1$ (and define $H^2$) in terms of cocycles and coboundaries.  

\begin{definition}{\rm  ($n$-cochains, $n=0,1,2,3$)}
\begin{itemize}
\item $C^0(J,\cm)=J\otimes\cm$
\item $C^1(J,\cm)=\hbox{ 1-cochains }=\{ f:J\rightarrow \cm: f \hbox{ linear}, f(1)=0\}$
\item $C^2(J,\cm)=\{q:J\times J:\rightarrow \cm: q(x;y)\hbox{ quadratic in }x; \hbox{linear in }y, q(1;y)\equiv 0\}$
\item $C^3(J,\cm)=C^3_1(J,\cm)\oplus C^3_2(J,\cm)$, where
$$C^3_1(J,\cm)=\{q:J\times J\times J\rightarrow \cm: q(x;y;z)\hbox{ quardic in }x;\hbox{ quadratic in }y;\hbox{ linear in }z\}$$
$$C^3_2(J,\cm)=\{q:J\times J\times J\rightarrow \cm: q(x;y;z)\hbox{ cubic in }x;\hbox{ linear in }y,z\}$$
\end{itemize}
\end{definition}

\begin{definition}(coboundary maps)
\begin{itemize}
\item $\partial^0:C^0\rightarrow C^1$: $\partial^0(x\otimes m)(y)=\tp{x}{m}{y}-\tp{m}{x}{y}$
\item $\partial^1:C^1\rightarrow C^2$: $\partial^1f(x,y)=f\tp{x}{y}{x}-2\tpc{x}{y}{f(x)}-\tpc{x}{f(y)}{x}$
\item $\partial^2=\partial^2_1\oplus \partial^2_2:C^2\rightarrow C^3$, where
$\partial^2_1:C^2\rightarrow C^3_1$ and $\partial^2_2:C^2\rightarrow C^3_2$ are given by  
\begin{eqnarray*}
\partial^2_1q(x;y;z)&=&
q(\tp{x}{y}{x};z)+\tpc{\tp{x}{y}{x}}{z}{q(x;y)}\\
&-&q(x;\tpc{y}{\tp{x}{z}{x}}{y})-\tpc{x}{q(y;\tp{x}{z}{x})}{x}\\
&-&\tpc{x}{\tpc{y}{q(x;z)}{y}}{x},
\end{eqnarray*}
\begin{eqnarray*}
\partial^2_2q(x;y;z)&=&
q(x,\tp{x}{z}{x};y)
+2\tpc{x}{y}{q(x;z)}\\
&-&2q(x;\tp{y}{x}{z})-\tpc{x}{q(y,z;x)}{x},
\end{eqnarray*}
and where $q(x,y;z)=(q(x+y,z)-q(x,z)-q(y,z))/2$ denotes the bilinearization of $q(x;z)$ in the first variable.
\end{itemize}
\end{definition}
 
 \begin{definition} (cocycles and coboundaries)
 \begin{itemize}
 \item   $B^1(J,\cm)=\hbox{Im}\, \partial^0=\{\sum(\tp{x_i}{m_i}{\cdot}-\tp{m_i}{x_i}{\cdot})\}$ (=$\hbox{Inder}\, (J,\cm)$)
 \item $Z^1(J,\cm)=\ker\partial^1 (=\hbox{Der}\, (J,\cm))$
 \item $B^2(J,\cm)=\hbox{Im}\, \partial^1$
  \item $Z^2(J,\cm)=\ker \partial^2$ ($=\ker\partial^2_1\oplus \ker\partial^2_2$).
 \end{itemize}
\end{definition}
 Then $H^n(J,\cm)=Z^n(J,\cm)/B^n(J,\cm)$, $n=1,2$.\smallskip
 
 The paper \cite{McCrimmon71}, which is mostly concerned with representation theory, proves a version of Theorem~\ref{thm:0207141}  for the cohomology groups just defined,  and the linearity of the functor $H^n$:
  \[
 H^n(J,\oplus_i M_i)=\oplus_i H^n(J,M_i),\quad n=1,2.
 \]
The paper \cite{McCrimmon82}, which is mostly concerned with compatibility of tripotents in Jordan triple systems, proves versions of the linearity of the functor $H^n$ corresponding to the Jordan triple structure.
 
\section{Cohomology of operator algebras}\label{chpt6}

We are now going to summarize the main points of this theory using the  two survey articles \cite{Ringrose96}, \cite{SinSmi04} as our primary guide, and using the notation found there.

\subsection{Continuous Hochschild cohomology}
 
As we saw in section~\ref{sec:cohomology}, 
Hochschild cohomology involves an associative algebra $A$ and $A$-bimodules $X$ and gives rise to

\begin{itemize}
\item $n$-cochains $L^n(A,X)$, 
 \item coboundary operators $\Delta_n$, 
 \item $n$-coboundaries $B^n$, 
 \item $n$-cocycles $Z^n$ and 
 \item
 cohomology groups $H^n(A,X)$.
\end{itemize}

If $A$ is a Banach algebra and $X$ is a Banach $A$-bimodule ( that is, a Banach space with module actions jointly continuous) we have the continuous versions of the above concepts 
$L^n_c\ ,\ B_c^n\ ,\ Z_c^n\ ,\  H_c^n(A,X)$.
Warning: $B_c^n$ is not always closed, so $H_c^n$ is still only a vector space.
\smallskip

Let $A$ be a C*-algebra of operators acting on a Hilbert space $H$ and let $X$ be a dual normal $A$-module ($X$ is a dual space and the module actions are separately ultra weakly-weak*-continuous). Further, we now have 
\begin{itemize}
\item normal $n$-cochains $L^n_w(A,X)$, that is,  bounded and separately weakly continuous $n$-cochains 
 \item coboundary operators $\Delta_n$, 
 \item normal $n$-coboundaries $B^n_w$, 
 \item normal $n$-cocycles $Z^n_w$ and 
 \item
 normal cohomology groups $H^n_w(A,X)$.
\end{itemize}

Thus, for a C*-algebra acting on a Hilbert space we  have three possible cohomology theories:
\begin{itemize}
\item the purely algebraic Hochschild theory $H^n$
\item the bounded theory $H^n_c$
\item the normal theory $H^n_w$
\end{itemize}

\begin{theorem}[Johnson,Kadison,Ringrose 1972 \cite{JohKadRin72}] \label{thm:1114111}
$H^n_w(A,X)\sim H_w^n(R,X)$\\
 ($R=$ultraweak closure of $A$)
\end{theorem}

\begin{theorem}[Johnson,Kadison,Ringrose 1972 \cite{JohKadRin72}]\label{thm:1114112}
 $H^n_w(A,X)\sim H_w^c(A,X)$
\end{theorem}

By Theorems~\ref{thm:1114111} and ~\ref{thm:1114112}, all four cohomology groups
\[
H^n_w(A,X)\ , \  H_w^n(R,X)\ ,\ H^n_c(R,X)\ ,\  H_w^n(R,X)
\]
are pairwise isomorphic.

\begin{theorem}[Johnson,Kadison,Ringrose 1972 \cite{JohKadRin72}]
$H^n_c(R,X)=0\ \forall n\ge 1$\\ 
 ($R=$hyperfinite von Neumann algebra)
\end{theorem}

\begin{theorem}[Connes 1978 \cite{Connes78}]
 If $R$ is a von Neumann algebra with a separable predual, and $H^1_c(R,X)=0$ for 
 every dual normal $R$-bimodule $X$, then  $R$ is hyperfinite.
\end{theorem}

At this point, there were two outstanding problems of special interest;

\begin{problem}\label{prob:1114111}
  $H_c^n(R,R)=0\ \forall n\ge 1$
 for every von Neumann algebra $R$?
 \end{problem}

\begin{problem}\label{prob:1114112}
$H_c^n(R,B(H))=0\  \forall n\ge 1$
 for every von Neumann algebra $R$ acting on a Hilbert space $H$?
 \end{problem}
 
 In the words of John Ringrose, 
 ``The main obstacle to advance was a paucity of information about the general bounded linear (or multilinear) mapping between operator algebras.  The major breakthrough, leading to most of the recent advances, came through the development of a rather detailed theory of completely bounded mappings.''
  This development occurred about a decade later and will be described in subsection~\ref{6.3}

 \subsection{Completely bounded maps}
 
 In this subsection, we shall give an introduction to the category of operator spaces from a point of view advantageous to the goals of this paper.  The immediate purpose is to give context to completely bounded maps for use when we return to cohomology in the next subsection.  Besides the standard monoraphs which appeared after the beginning of the 21st century \cite{EffRua00}, \cite{Paulsen02},\cite{Pisier03},\cite{BleLeM04},\cite{Helemskii10}, there is quite a good bit of the basic theory in the lecture notes of Runde \cite{RundeLNIM}.  Proofs of all the results in this subsection  can be found in one or more of these references.

\subsubsection{Banach spaces}

Why are normed spaces important?  One answer is that 
$\RR^n$ is a vector space with a norm, so you can take derivatives
and integrals.  Thus, normed spaces are important because you can do calculus.
\smallskip

Why is  completeness important?  
Three basic principles of functional analysis are based in some
way on completeness:
\begin{itemize}
\item Hahn-Banach theorem  (which depends on Zorn's lemma)
\item Open mapping theorem (which depends on Baire category)
\item Uniform boundedness theorem  (which depends on Baire category)
\end{itemize}

\noindent Some examples which will be of interest to us are:
\begin{itemize}
\item $C[0,1],\ L^p, \ \ell^p,\ c_0$, $C(\Omega),\ L^p(\Omega,\mu),\ C_0(\Omega)$
\item $B(X,Y)$, Banach algebras, Operator algebras ${\mathcal A}\subset B(X)$
\item Operator spaces; C*-algebras, Ternary rings of operators
\end{itemize}

The Hahn-Banach theorem states that every bounded linear functional
on a subspace of a Banach space has a bounded extension to the larger space 
with the same norm.

\begin{definition}  A Banach space is {\bf injective} if every bounded linear 
\underline{operator} from a subspace of a Banach space into it has a bounded 
extension to the larger space with the same norm.
\end{definition}

\begin{theorem}[Kelley 1952 \cite{Kelley52}]
A Banach space is injective if and only if it is isometric to
$C(\Omega)$, where $\Omega$ is extremely disconnected.
\end{theorem}

Injectivity in a category of Banach spaces is intimately related to contractive projections. Our first example of this is the following easily verified fact.
\begin{remark} A Banach space is injective if and only if there is an injective
Banach space containing it and a contractive projection of that space onto it.
\end{remark}

\subsubsection{Operator spaces}

Every Banach space $X$  is commutative, by which we mean that by the Hahn-Banach theorem, $X\subset C(X_1^*)$ (where $X_1^*$ denotes the unit ball of $X^*$), that is, $X$ is a linear subspace
of a {\bf commutative} C*-algebra. Hence: $C(\Omega)$ is the ``mother of all Banach spaces.''

\begin{definition}
An {\bf operator space} (or noncommutative Banach space, or
quantized Banach space) is a linear subspace of a C*-algebra (or $B(H)$).
\end{definition}

Hence: $B(H)$ is the ``mother of all operator spaces.''

An operator space may be viewed as a C*-algebra for which the
multiplication and involution have been ignored.  
You must replace
these by some other structure. 
What should this structure be and
what can you do with it?
The answer lies in the morphisms. 
In the category of Banach spaces,
the morphisms are the bounded linear maps. 
In the category of operator
spaces, the morphisms are the {\bf completely bounded maps}.
\smallskip

Operator spaces are intermediate between Banach spaces and C*-algebras.
The advantage of operator spaces over C*-algebras is that they allow
the use of finite dimensional tools and isomorphic invariants (``local
theory'').  C*-algebras are too rigid: morphisms are contractive, norms
are unique.
Operator space theory opens the door to a massive transfer of technology 
coming from  Banach space theory. 
This quantization process has benefitted 
operator algebra theory mainly (as opposed to Banach space theory).

\begin{definition}  An operator space is {\bf injective} if every completely contractive  linear 
operator from a subspace of an operator space into it has a completely contractive 
extension to the larger space. (See the beginning of section~\ref{sec:9} for the definition of completely bounded map)
\end{definition}

\begin{theorem} [Arveson-Wittstock-Paulsen Hahn-Banach] $B(H)$ is injective in the category of operator spaces and
completely contractive maps.
\end{theorem}

Our second example of the connection of contractive projections to injectivity is the following early verified fact.

\begin{remark} $X\subset B(H)$ is injective if and only if there is 
a completely contractive projection $P:B(H)\rightarrow B(H)$ with
$P(B(H))=X$.
\end{remark}

\subsubsection{Injective and mixed injective  operator spaces}

\begin{theorem} An operator space is injective if and only if it
is completely isometric to a corner of a C*-algebra. In particular,
it is completely isometric to a ternary ring of operators (TRO).
\end{theorem}

\begin{theorem} Every operator space has an injective envelope.
\end{theorem}
There are also the related notions of ternary envelope and Shilov boundary, which are useful but are not needed here, so will not be defined.
The following concept will appear later in part 3.

\begin{definition} An operator space is a {\bf mixed injective} if every completely
contractive  linear operator from a subspace of an operator space
into it has a bounded contractive
extension to the larger space.
\end{definition}

\begin{remark}
$X\subset B(H)$ is a mixed injective if and only if there is 
a contractive projection $P:B(H)\rightarrow B(H)$ with
$P(B(H))=X$.\end{remark}

\subsubsection{Applications of operator space theory}

The following theorems highlight the profound applications that operator space theory has had  to operator algebra theory.

\begin{itemize}
\item {\bf Similarity problems}. The Halmos problem is solved in the following theorem. For progress on and a description of the  Kadison problem, see \cite{PisierLNIM1618}.

\begin{theorem} [Pisier 1997 \cite{Pisier97}]
There 
exists a polynomially bounded operator which is not
similar to a contraction.
\end{theorem}

\item  {\bf Tensor products}

\begin{theorem}[Junge-Pisier 1995 \cite{JunPis95}]
$B(H)$ is not a nuclear pair.  That is,
$B(H)\otimes B(H)$ does not have a unique C*-norm.
\end{theorem}

To quote Pisier \cite{Pisier98}: ``This is a good case study. The proof is based on the solution to a problem
which would be studied for its own sake: 
whether the set of finite
dimensional operator spaces is separable (it is not).''

\item {\bf Operator amenable groups}. A well-known theorem of Johnson, which is false for the Fourier algebra, states that a group is amenable if and only if its group algebra is
an amenable Banach algebra under convolution.

\begin{theorem}[Ruan  1995 \cite{Ruan95}] A group is amenable if and only if its Fourier algebra is
amenable in the operator space formulation. 
\end{theorem}

\item Operator local reflexivity. A well known theorem of Lindenstrauss-Rosenthatl with many applications states that every Banach space is locally reflexive (Finite dimensional
subspaces of the second dual can be approximated by finite dimensional
subspaces of the space itself).  One application occurs in the proof that the second dual of a JB*-triple is a JB*-triple, with the consequence being the Gelfand-Naimark for JB*-triples.
Not all operator spaces are locally reflexive (in the appropriate
operator space sense).  However, there is the following theorem.

\begin{theorem}[Effros-Junge-Ruan 2000 \cite{EffJunRua00}] The dual of a C*-algebra is a locally reflexive operator space.
\end{theorem}

\end{itemize}
\subsection{Completely bounded cohomology}\label{6.3}

After this brief excursion into the world of operator spaces, we return to cohomology.

Let $A$ be a C*-algebra and let $S$ be a von Neumann algebra, both acting on the same Hilbert space $H$ with $A\subset S$.  We can view $S$ as  a dual normal $A$-module with $A$ acting on $S$ by left and right multiplication. We now have 
\begin{itemize}
\item completely bounded $n$-cochains $L^n_{cb}(A,S)$
 \item coboundary operators $\Delta_n$, 
 \item completely bounded $n$-coboundaries $B^n_{cb}$, 
 \item completely bounded  $n$-cocycles $Z^n_{bc}$ 
 \item
completely bounded cohomology groups $H^n_{cb}(A,S)$.
\end{itemize}

For a C*-algebra $A$ and a von Neumann algebra $S$ with $A\subset S\subset B(H)$ we thus have two new  cohomology theories:
\begin{itemize}
\item the completely bounded theory $H^n_{cb}$
\item the completely bounded normal theory $H^n_{cbw}$
\end{itemize}

By straightforward analogues of Theorems~\ref{thm:1114111} and ~\ref{thm:1114112}, the four cohomology groups
\[
H^n_{cb}(A,S)\ , \  H_{cbw}^n(A,S)\ ,\ H^n_{cb}(R,S)\ ,\  H_{cbw}^n(R,S)
\]
are mutually isomorphic, where $R$ is the ultraweak closure of $A$.

\begin{theorem}[Christensen-Effros-Sinclair 1987 \cite{ChrEffSin87}]
$H^n_{cb}(R,B(H))=0\ \forall n\ge 1$
 ($R=$any von Neumann algebra acting on $H$)
\end{theorem}

\begin{theorem}[Christensen-Sinclair 1987 \cite{ChrSin87}]
$H^n_{cb}(R,R)=0\ \forall n\ge 1$
 ($R=$any von Neumann algebra)
 \end{theorem}

 \subsubsection{Advances using complete boundedness}

A noteworthy quote from \cite{SinSmi04} is the following:
``Cohomology and complete boundedness have enjoyed a symbiotic relationship where advances in one have triggered progress in the other.''

\smallskip

Theorems~\ref{thm:1114113} and ~\ref{thm:1114114} are also due to Christensen-Effros-Sinclair \cite{ChrEffSin87}.

\begin{theorem}\label{thm:1114113}
$H^n_{c}(R,R)=0\ \forall n\ge 1$
 ($R=$ von Neumann algebra of type $I$, $II_\infty$, $III$, or of type $II_1$ and stable under tensoring with the hyperfinite factor)\end{theorem}

\begin{theorem}\label{thm:1114114}
$H^n_{c}(R,B(H))=0\ \forall n\ge 1$
 ($R=$ von Neumann algebra  of type $I$, $II_\infty$, $III$, or of type $II_1$ and stable under tensoring with the hyperfinite factor, acting on a Hilbert space $H$)\end{theorem}

 The following theorem is the culmination of earlier work in \cite{PopSmi94} and \cite{ChrPopSinSmi97}.

\begin{theorem}[Sinclair-Smith \cite{SinSmi98}]
$H^n_{c}(R,R)=0\ \forall n\ge 1$
 ($R=$ von Neumann algebra  of type  $II_1$ with a Cartan subalgebra and a separable\footnote{The separability assumption was removed in 2009 \cite{Cameron09}} predual)
\end{theorem}

The following theorem was proved in \cite{ChrPopSinSmi03AM} and \cite{ChrPopSinSmi03PNAS} and was new only for $n\ge 3$.  For $n=1$, the two cases are in \cite{Kadison66},\cite{Sakai66} and \cite{ChrPopSinSmi97}, while for $n=2$ one can refer to \cite{ChrSin87} and \cite{Christensen01}.

\begin{theorem}[Christensen-Pop-Sinclair-Smith $n\ge 3$ 2003] $H^n_{c}(R,R)=H^n_{c}(R,B(H))=0\ \forall n\ge 1$
 ($R=$ von Neumann algebra \underline{factor}  of type  $II_1$ with property $\Gamma$, acting on a  Hilbert space $H$)\end{theorem}
 
 A more recent result is the main result in \cite{PopSmi2010JFA}
\begin{theorem}
Let $M$ and $N$ be von Neumann algebras of type $II_1$.   Then $H_c^2(M\overline\otimes N,M\overline\otimes N)=0$.
\end{theorem}

We can now add a third problem to our previous two (Problems~\ref{prob:1114111},\ref{prob:1114112}). A candidate for a counterexample is the factor arising from the free group on 2 generators.

\begin{problem}
$H_c^n(R,R))=0\  \forall n\ge 2$?
($R$ is a von Neumann algebra of type $II_1$) 
\end{problem}

\subsubsection{Another approach (Paulsen)}
 
A different approach, which is popular among algebraists is taken in the paper \cite{Paulsen98}.  We quote from the review of this paper in Mathematical Reviews \cite{Rozenblumrev}.

``Let A be a subalgebra in the algebra B(H) of bounded operators in a Hilbert space H, and X be a linear subspace in B(H) such that AXA$\subset$X. Then X can be considered as an A-bimodule, and the Hochschild cohomology groups for the pair (A,X) can be constructed in the usual way, but with the cocycles being n-linear maps from A to X satisfying certain extra `complete boundedness' conditions.''

``Continuing the study of A. Ya. Khelemskii \cite{Helemskii89}, himself and others, Paulsen presents two new alternative presentations of these completely bounded Hochschild cohomologies, one of them as a relative Yoneda cohomology, i.e. as equivalence classes of relatively split resolutions, and the second as a derived functor, making it similar to EXT. 
These presentations make clear the importance of the notions of relative injectivity, projectivity and amenability, which are introduced and studied.''

 ``Using the relative injectivity of von Neumann algebras, the author proves the triviality of completely bounded Hochschild cohomologies for all von Neumann algebras. 
 The Yoneda representation makes the proofs of a number of classical results more transparent.''
 
\subsection{Perturbation of Banach algebras}

One of the most striking applications of the cohomology theory of Banach algebras lies in the area of perturbation theory.

\subsubsection{Theorem of Johnson and  Raeburn-Taylor}
\begin{definition}
A Banach algebra  is {\bf stable} if any other product making it into a Banach algebra, which is sufficiently  close to the original product gives rise to a Banach algebra which is topologically algebraically isomorphic to the original Banach algebra..
\end{definition}
Let's make this more precise. 
If $m$ is a Banach algebra multiplication on $A$, then its norm is defined by
$\|m(x,y)\|\le \|m\|\|x\|\|y\|$.  The following theorem was proved simultaneously and independently by Johnson \cite{Johnson77} and the team of Raeburn and Taylor \cite{RaeTay77}.

\begin{theorem} \label{thm:0215121} For any Banach algebra $A$ with multiplication $m$,
  if $H^2(A,A)=H^3(A,A)=0$, then 
there exists $\epsilon>0$ such that if $m_1$ is another Banach algebra multiplication on $A$ such that $\|m_1-m\|<\epsilon$ then $(A,m_1)$ and $(A,m)$ are topologically algebraically isomorphic.
\end{theorem}

For a stimulating discussion of this topic and many others, we recommend \cite{Helemskii89}.
For a similar problem on operator algebras one can look into \cite{KadKas72}.

\subsubsection{Perturbation for nonassociative Banach algebras}

The origin of perturbation theory for general non associative algebras is deformation theory.
Let $c_{ij}^k$ be the structure constants of a finite dimensional Lie algebra $L$. 
Let  $c_{ij}^k(\epsilon)\rightarrow c_{ij}^k$ 
Stability in this context means $(L,c_{ij}^k(\epsilon))$ is isomorphic to $(L,c_{ij}^k)$ if $\epsilon $ is sufficiently small.

\begin{theorem}[Gerstenhaber 1964 \cite{Gerstenhaber64AM}]
Finite dimensional semisimple Lie algebras are stable.
\end{theorem}

Perturbation results similar to Theorem~\ref{thm:0215121} have been proved recently for some non associative Banach algebras in \cite{Dosi09}.

\begin{theorem}[Dosi 2009]
If $L$ is a Banach Lie algebra and 
 $H^2(L,L)=H^3(L,L)=0$, then $L$ is a stable Banach Lie algebra.
\end{theorem}

An entirely similar result for Banach Jordan algebras is also proved in \cite{Dosi09}.  Of necessity, the definition of cohomology for Jordan algebras in \cite{Dosi09} is made only in dimensions 3 or less.

\section{Cohomology of   triple systems}

In this section we shall describe the results which the author learned about in preparing this mini course.   Since the author of this paper is not (yet) an expert in the purely algebraic side of cohomology theory, he is going to rely on information obtained from reviews of these papers in Mathematical Reviews. 

\subsection{Cohomology of finite dimensional triple systems}
\subsubsection{Cohomology of Lie triple systems}

The earliest work on cohomology of triple systems seems to be \cite{Harris61}, of which the following description is taken from its review in Mathematical Reviews \cite{Legerrev}. Four decades later, the second paper on the subject appeared \cite{HodPar02}.

A Lie triple system $T$ is a subspace of a Lie algebra $L$ closed under the ternary operation $[xyz]=[x,[y,z]]$ or, equivalently, it is the subspace of $L$ consisting of those elements $x$ such that $\sigma (x)=-x$, where $\sigma$ is an involution of $L$. 
A $T$-module $M$ is a vector space such that the vector-space direct sum $T\oplus M$ is itself a Lie triple system in such a way that

\begin{enumerate}
\item $T$ is a subsystem
\item $ [xyz]\in M$ if any of $x,y,z$ is in $M$
\item $ [xyz]=0$ if two of $x,y,z$ are in $M$. 
\end{enumerate}

 A universal Lie algebra $L_u(T)$ and an $L_u(T)$-module $M_s$ can be constructed in such a way that both are operated on by an involution $\sigma$ and so that $T$ and $M$ consist of those elements of $L_u(T)$ and $M_s$ which are mapped into their negatives by $\sigma$. 
\smallskip

   Now suppose $L$ is a Lie algebra with involution $\sigma$ and $N$ is an $L$-$\sigma$ module. Then $\sigma$ operates on $H^n(L,N)$ so that $$H^n(L,N)=H_+^n(L,N)\oplus H_-^n(L,N)$$ with both summands invariant under $\sigma$. 
      The cohomology of the Lie triple system is defined by $H^n(T,M)=H_+^n(L_u(T),M_s)$.  Harris 
   investigates these groups for $n=0,1,2$. 
      
\begin{enumerate}
\item   $H^0(T,M)=0$ for all $T$ and  $M$
   \item $H^1(T,M)=$ derivations of $T$ into $M$ modulo inner derivations    
\item    $H^2(T,M)=$ factor sets of $T$ into $M$ modulo trivial factor sets. 
   \end{enumerate}

   Turning to the case of finite-dimensional simple $T$ and ground field of characteristic 0, one has the Whitehead lemmas   
    $H^1(T,M)=0=H^2(T,M)$ and    
     Weyl's theorem: Every finite-dimensional module is semi-simple. 
     Also, if in addition, the ground field $\Phi$ is algebraically closed, then $H^3(T,\Phi)$ is 0 or not 0, according as $L_u(T)$ is simple or not.

The following is the verbatim review \cite{Koshrev} of \cite{HodPar02}.

``This is a study of representations of Lie triple systems, both ordinary
and restricted. 
The theory is based on the connection between Lie algebras and Lie triple systems.
In addition, the authors begin the study of the cohomology theory for Lie triple systems and their
restricted versions. 
They also sketch some future applications and developments of the theory.''

\subsubsection{Cohomology of associative triple systems}

The following is the verbatim review \cite{Seibtrev} of \cite{Carlsson76}.

``A cohomology for associative triple systems is defined, with the main purpose to get quickly the cohomological triviality of finite-dimensional separable objects over fields of characteristic $\ne  2$, i.e., in particular the Whitehead lemmas and the Wedderburn principal theorem. 

This is achieved by embedding an associative triple system $A$ in an associative algebra $U(A)$ and associating with every trimodule $M$ for $A$ a bimodule $M_u$ for $U(A)$ such that the cohomology groups $H^n(A,M)$ are subgroups of the classical cohomology groups $H^n(U(A),M_u)$.

 Since $U(A)$ is chosen sufficiently close to $A$, in order to inherit separability, the cohomological triviality of separable $A$ is an immediate consequence of the associative algebra theory. 
 
 The paper does not deal with functorialities, nor with the existence of a long exact cohomology sequence.''

\subsubsection{Wedderburn decomposition}
The classical Wedderburn decomposition is the following theorem, which we state in its Banach algebra form from \cite{BadCur60Wed}. 
It  has been generalized, not only to non associative algebras, see for example \cite{Jacobson62} and \cite{Jacobson68},  but   to some classes of Banach algebras, see for example \cite{Dales00}.

\begin{theorem}
Let $A$ be a Banach algebra with radical $R$.  If the dimension of $A$ is finite, there is  a subalgebra $B$ of $A$ such that $A=B+R$ and $B\cap R=0$.  Moreover, if $A$ is commutative, then $B$ is necessarily unique.
\end{theorem}

The following definition is the model to use in extending this concept to infinite dimensional triples.

\begin{definition}
If $A$ is a Banach algebra with radical $R$, then $A=B\oplus R$ is a strong Wedderburn decomposition of $A$ if $B$ is a closed subalgebra of $A$.
\end{definition}

The author knows of only two papers dealing with this concept in the context of triple systems, namely, \cite{Carlsson77} and \cite{KuhRos78}.
The following, including Theorem~\ref{thm:0215122} is the verbatim review \cite{Boersrev} of \cite{Carlsson77}.

The Wedderburn principal theorem, known for Lie triple systems, is proved for alternative triple systems and pairs.
If $i$ is an involution of an alternative algebra $ B$, then $\langle xyz\rangle:=(x\cdot  i(y))\cdot z$ is an alternative triple ($x,y,z\in B$). 
A polarisation of an alternative triple $A$ is a direct sum of two submodules $A^1\oplus A^{-1}$ with $$\langle A^1A^{-1}A^1\rangle\ \subset A^1, \langle A^{-1}A^1A^{-1}\rangle\ \subset A^{-1}$$ and

$$
\langle A^1A^1A^1\rangle\ =\langle A^{-1}A^{-1}A^{-1}\rangle\ =\langle A^1A^1A^{-1}\rangle\ =\langle A^{-1}A^{-1}A^1\rangle\ =\langle A^{-1}A^1A^1\rangle\ =\langle A^1A^{-1}A^{-1}\rangle =\{0\}.$$

 An alternative pair is an alternative triple with a polarisation. 

\begin{theorem}\label{thm:0215122}
 If $A$ is a finite-dimensional alternative triple system (or an alternative pair) over a field $K, R$ the radical and $A/R$ separable, then $A=B\oplus R$, where $B$ is a semisimple subtriple (subpair) of $A$ with $B=A/R$.
\end{theorem}

The following description of \cite{KuhRos78} is taken from the abstract and the review \cite{Thedyrev}.

This paper summarizes
some properties of Jordan pairs, states some results about some groups defined by Jordan pairs, and 
constructs a Lie algebra to a Jordan pair. 
This construction is a generalization of the well-known Koecher-Tits-construction. 
The radical of this Lie algebra is calculated in terms of the given Jordan pair and a Wedderburn decomposition theorem for Jordan pairs (and triples) in the characteristic zero case is proved. \smallskip

More precisely, an observation of Koecher that the theorem of Levi for Lie algebras of characteristic 0 implies the Wedderburn principal theorem for Jordan algebras is extended to Jordan pairs (and Jordan triples) V over a field of characteristic 0. In addition, the authors show that any two Wedderburn splittings of V are conjugate under a certain normal subgroup of the automorphism group of V.

\subsubsection{Cohomology of algebras and triple systems}
The reader may be wondering about the cohomology of Jordan triple systems.  As noted in \ref{sub:hint}, there is a hint of this in the paper of McCrimmon \cite{McCrimmon82}.  Building on the 1 and 2 dimensional cohomology for quadratic Jordan algebras, \cite{McCrimmon82}  broaches the analogous construction for Jordan triple systems.\smallskip

There is also a more general approach in the paper of Seibt \cite{Seibt75}.  
The following is from the introduction to \cite{Seibt75}. 
\smallskip

``The classical cohomologies of unital associative algebras and of Lie algebras have both a double algebraic character:
They are embedded in all of the machinery of derived functors, and they allow full extension theoretic interpretations ({\bf Yoneda}-interpretation) of the higher cohomology groups---which seems natural since the coefficient category for cohomology is actually definable in terms of singular extension theory.

 If one wants to define a uniform cohomology theory for (linear) nonassociative algebras and triple systems which "generalizes'' these two classical cohomologies one may proceed either via derived functors \cite{BarRin66}
    or via singular extension theory \cite{Gerstenhaber64}, \cite{Glassman70PJM}
    
 The purpose of this paper which adopts the first point of view is to discuss compatibility questions with the second one.''

\subsection{Cohomology of Banach triple systems---Prospectus}\label{7.2}

In this subsection we propose, albeit quite vaguely and briefly, some places to look for extending some of the preceding material to infinite dimensional Banach triples of various kinds as well as infinite dimensional non associative Banach algebras\footnote{See the forthcoming monograph {\it 
Non-Associative Normed Algebras: Volume 1, The Vidav-Palmer and Gelfand-Naimark Theorems} (Encyclopedia of Mathematics... by Miguel Cabrera Garc'a and Angel Rodr'guez Palacios)}. We begin with the following quotation from \cite{Kadison00}.
\begin{quotation}
"A veritable army of researchers took the theory of derivations of operator algebras to dizzying heights---producing a theory of cohomology of operator algebras as well as much information about automorphisms of operator algebras." ---Dick Kadison ({\it Which Singer is that?} 2000)
\end{quotation}

In addition to associative algebras, cohomology groups are defined for Lie algebras and, to some extent, for Jordan algebras.  Since the structures of Jordan derivations and Lie derivations on von Neumann algebras are well understood, isn't it time to study the higher dimensional non associative cohomology of a von Neumann algebra?  This subsection includes an introduction to the second Jordan cohomology groups of a von Neumann algebra.

Also included, with or without comment,  are some papers related to topics discussed in this survey which the author has downloaded and  thinks might be worth exploring. In some cases possible connections between papers are suggested.

\subsubsection{Jordan 2-cocycles on von Neumann algebras}

Let $M$ be a von Neumann algebra.   A {\bf Hochschild 2-cocycle} is a bilinear map $f:M\times M\rightarrow M$ satisfying
\begin{equation}\label{eq:0504131}
af(b,c)-f(ab,c)+f(a,bc)-f(a,b)c=0
\end{equation} An example is a {\bf Hochschild 2-coboundary}:  $f(a,b)=a\mu(b)-\mu(ab)+\mu(a)b$ for some $\mu:M\rightarrow M$ linear.
A {\bf Jordan 2-cocycle} is a  bilinear map $f:M\times M\rightarrow M$ satisfying
\[
f(a,b)=f(b,a) \hbox{ (symmetric) }
\] 
\begin{equation}\label{eq:0520131}
f(a^2,a\circ b)+f(a,b)\circ a^2+f(a,a)\circ(a\circ b)
\end{equation}
\[-f(a^2\circ b,a)-f(a^2,b)\circ a-(f(a,a)\circ b)\circ a=0
\]
An example is a {\bf Jordan 2-coboundary}:$f(a,b)=a\circ \mu(b)-\mu(a\circ b)+\mu(a)\circ b$ for some 
$\mu:M\rightarrow M$ linear.
We have the following cohomology groups and their interpretations:
$$H^1(M,M)=\frac{1\hbox{-cocycles}}{1\hbox{-coboundaries}}=\frac{\hbox{derivations}}{\hbox{inner derivations}}$$
$$H^1_J(M,M)=\frac{\hbox{Jordan }1\hbox{-cocycles}}{\hbox{Jordan }1\hbox{-coboundaries}}=\frac{\hbox{Jordan derivations}}{\hbox{inner Jordan derivations}}$$
$$H^2(M,M)=\frac{2\hbox{-cocycles}}{2\hbox{-coboundaries}}
=\frac{\hbox{null extensions}}{\hbox{split null extensions}}$$
$$H^2_J(M,M)=\frac{\hbox{Jordan }2\hbox{-cocycles}}{\hbox{Jordan }2\hbox{-coboundaries}}
=\frac{\hbox{Jordan null extensions}}{\hbox{Jordan split null extensions}}$$

For ``almost all'' von Neumann algebras, $H^2(M,M)=0$. For finite dimensional $M$, the Albert-Penico-Taft theorem (step 7 in subsection~\ref{5.3.1}) shows that $H^2_J(M,M)=0$? This is still unknown for infinite dimensional $M$. 

 There are at least two elegant approaches to this problem: Jordan classification (subsection \ref{5.3.1}) and Lie algebras (not discussed in this survey).  We present now a pedestrian approach based on solving linear equations and motivated by the proof of reduction IV at the end of subsection~\ref{5.3.1}.  We first illustrate with two known associative cases.\smallskip
 
 {\bf Test case 1}

Let $h$ be a Hochschild 1-cocycle, that is, a linear map $h:\mn\rightarrow \mn$ satisfying
$
h(ab)-ah(b)-h(a)b=0.
$
To show that there is an element $x\in \mn$ such
that 
$
h(a)=xa-ax,
$
it is enough to prove this with $a\in\{    \e{i}{j}    \}$.
With 
\begin{equation}\label{eq:0515144}
x=\sum_{p,q}x_{pq}\ee{pq}.
\end{equation}
and  $\gamma_{ijpq}$ defined by
\begin{equation}\label{eq:0515145}
h(\ee{ij})=\sum_{p,q}\gamma_{ijpq}\ee{pq},
\end{equation}
we arrive at the system of  linear vector equations
\begin{equation}\label{eq:0515146}
\sum_{p,q}\gamma_{ijpq}\ee{pq}=\sum_{p,q}\delta_{qi}x_{pq}\ee{pj}-\sum_{p,q}\delta_{jp}x_{pq}\ee{iq}.
\end{equation}
with $n^2$ unknowns $x_{ij}$. Then any solution of (\ref{eq:0515146}) proves the result.\smallskip

{\bf Test case 2}

Let $h$ be a Hochschild 2-cocycle, that is, a bilinear map $h:\mn\times\mn\rightarrow \mn$ satisfying
$
ah(b,c)-h(ab,c)+h(a,bc)-h(a,b)c=0.
$
 To show that there is a linear transformation $\mu:\mn\rightarrow \mn$ such
that 
$
h(a,b)=\mu(ab)-a\mu(b)-\mu(a)b,
$ 
it is enough to prove that this holds with $a,b\in\{    \e{i}{j}    \}$, that is
\begin{equation}\label{eq:0419142}
h(\e{i}{j},\e{k}{l})=\delta_{jk}\mu(\e{i}{l})-\e{i}{j}\mu(\e{k}{l})-\mu(\e{i}{j})\e{k}{l}.
\end{equation}
With
$
\mu(\e{i}{j})=\sum_{k,l}\mu_{ijkl}\e{k}{l}
$ and
  $\gamma_{ijklpq}$ defined by
$
h(\e{i}{j},\e{k}{l})=\sum_{p,q}\gamma_{ijklpq}\e{p}{q},
$ \\
we arrive at  the system of  $n^6$ linear equations
\begin{equation}\label{eq:0419143}
\sum_{p,q}\gamma_{ijklpq}\e{p}{q}=\sum_{p,q}\delta_{jk}\mu_{ilpq}-\sum_{p,q}\delta_{jp}\mu_{klpq}\e{i}{q}-\sum_{p,q}\delta_{qk}\mu_{ijpq}\e{p}{l},
\end{equation}
with $n^4$ unknowns $\mu_{ijkl}$. Then any solution of (\ref{eq:0419143}) proves (\ref{eq:0419142}).
\smallskip

We propose to use this method to prove that $H^2_J(M,M)=0$ if $M$ is a finite von Neumann algebra of type I.  (A separate proof would be needed for commutative $M$.)
Let $M=M_2(L^\infty(\Omega))$ be a finite von Neumann algebra of type $I_n$ with $n=2$.
Let $f$ be a Jordan 2-cocycle, that is, a symmetric bilinear map $f:M\times M\rightarrow M$ with
\begin{equation}\label{eq:0510141}
f(a^2,ab)+f(a,b)a^2+f(a,a)(ab)-f(a^2b,a)-f(a^2,b)a-(f(a,a)b)a=0.
\end{equation}
(To save space, $ab$ temporarily denotes the Jordan product in the associative algebra $M$)

\medskip

To show that there is a linear transformation $\mu:M\rightarrow M$ such
that 
\begin{equation}\label{eq:0509148}
f(a,b)=\mu(ab)-a\mu(b)-\mu(a)b,
\end{equation}
it is enough to prove, for $a,b\in Z(M)$,
\begin{equation}\label{eq:0509149}
f(a\e{i}{j},b\e{k}{l})=\delta_{jk}\mu(ab\e{i}{l})-a\e{i}{j}\mu(b\e{k}{l})-\mu(a\e{i}{j})b\e{k}{l}.
\end{equation}

With
$
\mu(a\e{i}{j})=\sum_{k,l}\mu_{ijkl}(a)\e{k}{l}
$
and $\gamma_{ijklpq}(a,b)\in Z(M)$ defined by
\begin{equation}\label{eq:05091410}
f(a\e{i}{j},b\e{k}{l})=\sum_{p,q}\gamma_{ijklpq}(a,b)\e{p}{q},
\end{equation}
we arrive at the system of  $n^6$ linear vector equations 
with $3n^4$ unknowns $\mu_{ijkl}(ab),\mu_{ijkl}(a),\mu_{ijkl}(b)$
$$
2  \sum_{p,q}\gamma_{ijklpq}(a,b)\e{p}{q}=  \delta_{jk}\sum_{p,q} \mu_{ilpq}(ab)\ee{pq}
$$
\begin{equation}\label{eq:0510141bis}
-  \sum_{p,q} a(\delta_{jp}\mu_{klpq}(b)\e{i}{q}+\delta_{iq}\mu_{klpq}(b)\e{p}{j})
\end{equation}
$$
-  \sum_{p,q} b(\delta_{qk}\mu_{ijpq}(a)\e{p}{l}+\delta_{lp}\mu_{ijpq}(a)\e{k}{q}).
$$
 Then any solution of (\ref{eq:0510141bis}) proves (\ref{eq:0509149}) and hence (\ref{eq:0509148}).
(Warning: it is not clear that this approach will be successful!)\smallskip

We shall now describe some properties of Jordan 2-cocycles.
Just as every derivation is a Jordan derivation, it is also true that every symmetric Hochschild 2-cocycle is a Jordan 2-cocyle. Recall that every Jordan derivation on a semisimple Banach algebra is a derivation (Sinclair). If every Jordan 2-cocycle was a Hochschild 2-cocycle, we would have $H^2(M,M)=H^2_J(M,M)$.  Fortunately, we have the following examples.

\begin{proposition}\label{prop:7.4} Let $M$ be a von Neumann algebra.
\begin{description}
\item[(a)] Let $f:M\times M\rightarrow M$ be defined by $f(a,b)=a\circ b$.
Then $f$ is a Jordan 2-cocycle with values in $M$, which is not a Hochschild 2-cocycle unless $M$ is commutative.
\item[(b)] If $M$ is finite with a faithful normal finite trace $\tr$, then $f:M\times M\rightarrow M_*$ defined by $f(a,b)(x)=\tr((a\circ b)x)$ is a Jordan 2-cocycle with values in $M_*$ which is not a Hochschild 2-cocycle unless $M$ is commutative.
\end{description}
\end{proposition}

It can also be shown that if $f$ is a Jordan 2-cocycle on the von Neumann algebra $M$, then $f(1,x)=xf(1,1)$ for every $x\in M$ and $f(1,1)$ belongs to the center of $M$.  

Let us show the proof of (a) of Proposition~\ref{prop:7.4}. 
Recall the definition of Jordan 2-cocycle
\begin{equation}\label{eq:0530131}
f(a^2,a\circ b)+f(a,b)\circ a^2+f(a,a)\circ(a\circ b)
\end{equation}
\[-f(a^2\circ b,a)-f(a^2,b)\circ a-(f(a,a)\circ b)\circ a=0
\]
Let $f(a,b)=a\circ b$.   The equation (\ref{eq:0530131}) reduces to 
\[
a^2\circ(a\circ b)+(a\circ b)\circ a^2+a^2\circ(a\circ b)-(a^2\circ b)\circ a-(a^2\circ b)\circ a-(a^2\circ b)\circ a,
\]
which is zero by the Jordan axiom, so $f$ is a Jordan 2-cocycle.
If this $f$ were a Hochschild 2-cocycle, we would have
\[
c(a\circ b)-(ca)\circ b+c\circ(ab)-(c\circ a)b=0,
\]
which reduces to $[[c,b],a]=0$ and therefore $[M,M]\subset Z(M)$ (the center of $M$). Since
$M=Z(M)+[M,M]$, $M$ is commutative.  This proves (a).

While we are at it, let us prove that every symmetric Hochschild 2-cocycle is a Jordan 2-cocycle.
The Jordan 2-cocycle identity reduces to 
\[
f(a^2,ab)+f(a^2,ba)+f(a,b)a^2+a^2f(a,b)
\]
\begin{equation}\label{eq:0504132}
+f(a,a)ab/2+baf(a,a)/2-f(a^2b,a)-f(ba^2,a)
\end{equation}
\[
-f(a^2,b)a-af(a^2,b)-bf(a,a)a/2-af(a,a)b/2=0
\]
In (\ref{eq:0504131}), replace $a$ by $a^2$ and $c$ by $a$ to obtain
\[
a^2f(b,a)-f(a^2b,a)+f(a^2,ba)-f(a^2,b)a=0.
\]
Thus, taking into account that $f$ is symmetric, four of the twelve terms in (\ref{eq:0504132}) sum to zero.
In (\ref{eq:0504131}), replace $c$ by $a^2$  to obtain
\[
af(b,a^2)-f(ab,a^2)+f(a,ba^2)-f(a,b)a^2=0.
\]
Thus, four more of the twelve terms in (\ref{eq:0504132}) sum to zero.
It remains to show that
\begin{equation}\label{eq:0504133}
f(a,a)ab+baf(a,a)-bf(a,a)a-af(a,a)b=0.
\end{equation}
In (\ref{eq:0504131}), replace $b$ and $c$  by $a$  to obtain\quad
$
af(a,a)=f(a,a)a
$\\  from which (\ref{eq:0504133}) follows.

\subsubsection{Ternary weak amenability for group algebras and Fourier algebras}

The purpose of the summary in this subsection is to suggest that the concepts of triple derivation and ternary weak amenability be explored in the contexts of group algebras and Fourier algebras.

Group algebras $L^1(G)$ and Fourier algebras $A(G)$ of a locally compact group $G$ are closely related and share some techniques with operator algebras, as do the measure algebra $M(G)$ and the Fourier-Stieltjes algebra $B(G)$. The Fourier algebra and the Fourier-Stieltjes algebras were introduced in the thesis of Pierre Eymard \cite{Eymard64} and have provided fertile ground for the extension of many classical harmonic analysis results to the noncommutative setting.

Johnson showed very early that $G$ is an amenable group if and only if $L^1(G)$ is an amenable Banach algebra  \cite{Johnson72}. Two decades later he showed that $L^1(G)$ is weakly amenable \cite{Johnson91}. A decade after that it was shown that $M(G)$ is weakly amenable if and only if $G$ is discrete (and therefore $M(G)=\ell^1(G)$) with the consequence that $M(G)$ is amenable if and only if $G$ is discrete \cite{DalGhaHel02}. 

The Fourier algebra $A(G)$ is amenable if and only if $G$ has an abelian group of finite index. The ``if'' was proved in \cite{LauLoyWil96} and the ``only if'' in \cite{ForRun05}. The following was also  proved in the latter reference: $B(G)$ is amenable if and only if $G$ has a compact abelian subgroup of finite index.   

A characterization for weak amenability for the Fourier algebras is not yet known, but there are several partial results which can be found in the survey \cite{Spronk2010}, which contains much other information, some of which has been summarized here.  We mention only one such result: If $G$ is a compact group, then $A(G)$ is weakly amenable if and only if the connected component of the identity in $G$ is abelian. The  ``only if'' part was proved in \cite{ForSamSpr09} and the ``if'' part holds for all locally compact groups, as shown in  \cite{ForRun05}.

\subsubsection{Cohomology of commutative JB*-triples and TROs} 

The first three cohomology groups of a commutative $C^*$-algebra with values in a Banach module were studied in the pioneering paper \cite{Kamowitz62}.   The cohomology groups of a finite dimensional associative triple system with values in the relevant module were studied in the \cite{Carlsson76}. By some combination of the techniques in these two papers, it should be possible to develop the cohomology of a commutative $JB^*$-triple with values in a module.\smallskip

More generally, the cohomology of an arbitrary TRO (ternary ring of operators) should be pursued. Once this has been worked out for TROs, the author feels that the enveloping TRO of a $JC^*$-triple, discussed in section 11 may be helpful in working out the theory for $JC^*$-triples, as it seems to be for working out  the $K$-theory of $JB^*$-triples, studied in \cite{BohWer11bis}.

\subsubsection{Some other avenues to pursue}

\begin{itemize}

\item {Wedderburn decompositions for JB*-triples} (K\"uhn-Rosendahl \cite{KuhRos78})
\item {Low dimensional cohomology for JBW*-triples and algebras-perturbation} (Dosi \cite{Dosi09}, McCrimmon \cite{McCrimmon82})
\item {Structure group of JB*-triple} (McCrimmon \cite{McCrimmon82}, Meyberg \cite{Meyberg72})
\item {Alternative Banach triples} (Carlsson \cite{Carlsson77}, Braun \cite{Braun84})
\item {Completely bounded triple cohomology} (Bunce-Feely-Timoney \cite{BunFeeTim11}, Bunce-Timoney \cite{BunTim11}, Christensen-Effros-Sinclair \cite{ChrEffSin87})
\item {Lie algebraic techniques. Applications of the Koecher-Kantor-Tits construction}(\cite{Chu12})
\end{itemize}

As an illustration of the last bullet, we shall now give a preview of the forthcoming paper \cite{ChuRus14}, which develops a cohomology theory for Jordan triples based on the corresponding theory for Lie algebras (cf. Theorem~\ref{3.19}).

Let $V$ be a non-degenerate Jordan
triple and let $ L (V)$ be its TKK Lie algebra as defined in
\cite[p.45]{Chu12}, where
$$ L (V) = V \oplus V_0 \oplus V$$
and the Lie product is given by
$$[(x,h,y), (u,k,v)] = (hu - kx,\, [h,k]+x \bo v -u\bo y,\,
k^\natural y - h^\natural v).$$
We note that $V_0= {\rm Span}\{ V\bo V\}$ is a Lie subalgebra of
$ L (V)$ and for $h= \sum_i a_i \bo b_i \in V_0$, the map
$h^\natural : V \rightarrow V$ is defined by
$$h^\natural = \sum_i b_i \bo a_i.$$
Given a Lie algebra $ L$ and an $  L$-module $X$, we
denote the action of $  L$ on $X$ by
$$(\ell, x) \in   L \times X \mapsto \ell . x \in X$$
so that
$$ [\ell, \ell '].x = \ell.(\ell'. x) - \ell' .(\ell. x).$$

Let $V$ be a Jordan triple and $M$ be a triple $V$-module as
defined in subsection~\ref{5.5}, but for convenience, we omit the numeric
subscripts of the triple actions on $M$. Given $a,b\in V$, the box operator $a\bo b:V\rightarrow V$ can also  be considered  as a mapping from $M$ to $M$.  Similarly, for $u\in V$ and $m \in
M$, the ``{\it box operators}'' $$u \bo m,\, m\bo u : V
\longrightarrow M$$ are defined in a natural way as $v\mapsto \{u,m,v\}$ and $v\mapsto \{m,u,v\}$ respectively.  In view of the main identity for Jordan triples, for
$a\bo b \in V_0$, we define naturally\footnote{Note that the left sides of these two definitions are not  commutators, but just a convenient notation}
$$[a\bo b, u\bo m] = \{a,b,u\}\bo m - u\bo \{m,a,b\}= - [u\bo m, a\bo b] $$
and
$$[a\bo b, m\bo u] = \{a,b,m\}\bo u - m\bo \{u,a,b\} = -[m \bo u, a\bo b].$$
for $u \in V$ and $m\in M$.  Let $M_0$ be the linear span of
$$\{u\bo m, n\bo v: u,v \in V, m,n \in M\}$$
in the vector space $L(V,M)$ of linear maps from $V$ to $M$.
Extending the above product by linearity, we can define an action
of $V_0$ on $M_0$ by
$$ (h, \varphi) \in V_0 \times M_0 \mapsto  [h,\varphi]\in M_0.$$

\begin{lemma} $M_0$ is a $V_0$-module of the Lie algebra $V_0$.
\end{lemma}
\pf\
We are required to show that
\begin{equation}\label{eq:0710141}
[[h,k],\varphi]=[h,[k,\varphi]]-[k,[h,\varphi]].
\end{equation}
We can assume that $h=a\bo b$, $k=c\bo d$ and $\varphi=w\bo m$ or $m\bo w$.  We assume $\varphi=w\bo m$, the other case being similar.
For the left side of (\ref{eq:0710141}), we have
\begin{eqnarray*}
[[a\bo b,u\bo v],w\bo m]&=&[\tp{abu}\bo v-u\bo\tp{vab},w\bo m]\\
&=&\tpp{\tp{abu}}{v}{w}\bo m-w\bo \tpp{m}{\tp{abu}}{v}\\
&-&\tpp{u}{\tp{vab}}{w}\bo m+w\bo \tpp{m}{u}{\tp{vab}}\\
&=&(\tpp{\tp{abu}}{v}{w}-\tpp{u}{\tp{vab}}{w})\bo m\\
&-&w\bo (\tpp{m}{\tp{abu}}{v}-\tpp{m}{u}{\tp{vab}}.
\end{eqnarray*}
For the right side of (\ref{eq:0710141}), we have
\begin{eqnarray*}
\lefteqn{[a\bo b,[u\bo v,w\bo m]]-[u\bo v,[a\bo b,w\bo m]]=}\\
&&[a\bo b,\tp{uvw}\bo m-w\bo \tp{muv}]-[u\bo v,\tp{abw}\bo m-w\bo\tp{mab}]\\
&=&\tpp{a}{b}{\tp{uvw}}\bo w-\tp{uvw}\bo\tp{mab}-\tp{abw}\bo\tp{muv}+w\bo\tpp{\tp{muv}}{a}{b}\\
&-&\tpp{u}{v}{\tp{abw}}\bo m+\tp{abw}\bo \tp{muv}+\tp{uvw}\bo\tp{mab}-w\bo\tpp{\tp{mab}}{u}{v}\\
&=&(\tpp{a}{b}{\tp{uvw}}-\tpp{u}{v}{\tp{abw}} )\bo m-w\bo(\tpp{\tp{mab}}{u}{v}-\tpp{\tp{muv}}{a}{b}).
\end{eqnarray*}
(\ref{eq:0710141}) now follows from the main identity for Jordan triples.
\smallskip

Let $V$ be  Jordan triple and $  L (V)$ its TKK Lie algebra.
Given a triple $V$-module $M$, we now construct a corresponding
Lie module $\mathcal{L} (M)$ of the Lie algebra $  L (V)$ as
follows.
Define $$\mathcal{L} (M) = M \oplus M_0 \oplus M$$ and define the
action
$$ ((a,h,b), (m,\varphi,n)) \in   L(V) \times \mathcal{L} (M)
\mapsto (a,h,b). (m,\varphi,n) \in \mathcal{L} (M)$$ by
$$(a,h,b). (m,\varphi,n) = (hm-\varphi a,~ [h,\varphi] +a\bo
n - m\bo b,~ \varphi^\natural b - h^\natural(n)\,)$$ where, for $h
= \sum_i a_i\bo b_i$ and $\varphi = \sum_i u_i \bo m_i + \sum_j
n_j \bo v_j$, we have the following natural definitions
$$hm= \sum_i \{a_i,b_i,m\},\quad \varphi a= \sum_i \{u_i,  m_i,a\} +
\sum_j\{n_j,
 v_j, a\}, \quad  \varphi^\natural = \sum_i m_i \bo u_i +\sum_j v_j \bo
 n_j.$$

The following theorem is the starting point of the paper \cite{ChuRus14}.

 \begin{theorem}\label{m} Let $V$ be a Jordan triple and let $  L (V)$
 be its TKK Lie algebra. Let $M$ be a triple $V$-module. Then
 $\mathcal{L} (M)$ is a Lie $  L(V)$-module.
 \end{theorem}

\part{Quantum functional analysis}

\section{Contractive projections}

The purpose of this section is to describe the role which contractive projections have played in the theory of $JB^*$-triples, which is the primary algebraic structure of interest in this survey.   
Contractive projections on Banach spaces have been the focus of much study. For a survey see \cite{Rand01}.   We are concerned here  with the interplay between  algebraic structure and contractive projections.

\subsection{Projective stability}\label{proj1}

A well-known and useful result in the structure theory of operator triple
systems is the ``contractive projection principle,'' that is, the fact
that the range of a contractive projection on a \jbst\  is linearly
isometric in a natural way to another \jbst.
The genesis of the role of contractive projections in Banach triple systems lies in the following theorem. This was preceded by a similar principle in a geometric context \cite{Stacho82}, which was unknown to the authors of \cite{FriRus85JFA}.

\begin{theorem}[Friedman and Russo 1985 \cite{FriRus85JFA}]
The range of a contractive projection on a $C^*$-algebra is linearly isometric to a $JC^*$-triple, that is, a linear subspace of $B(H)$ which is closed under the symmetric triple product $xy^*z+zy^*x$. \end{theorem}

Thus, as shown later and described below, the  category of $JB^*$-triples and contractions is stable under
contractive projections.
To put this result in proper prospective, let $\mathcal
B$ be the category of Banach spaces and contractions. 
 We  say that a
sub-category $\mathcal S$ of $\mathcal B$ is {\bf projectively stable} if
it has the property that whenever $A$ is an object of $\mathcal S$ and $X$
is the range of an idempotent morphism of $\mathcal S$ on $A$,
then $X$ is isometric (that is,  isomorphic in $\mathcal S$) to an object in
$\mathcal S$.\smallskip

 Examples of projectively stable categories are the following.

\begin{itemize}
\item $L_1$, contractions (Grothendieck 1955 \cite{Grothendieck55})
\item $L^p$, $1\le p<\infty$, contractions \\
(Douglas 1965 \cite{Douglas65},
Ando 1966 \cite{Ando66}, Bernau-Lacey 1974 \cite{BerLac74}, Tzafriri 1969 \cite{Tzafriri69})
\item $C^*$-algebras, completely positive unital maps (Choi-Effros 1977 \cite{ChoEff77})
\item $\ell_p$, $1\le p<\infty$, contractions\\
 (Lindenstrauss-Tzafriri 1977 \cite{LinTza77})
\item $JC^*$-algebras, positive unital maps\\
 (Effros-Stormer 1979 \cite{EffSto79})
\item $TROs$ (ternary rings of operators), \\
complete contractions
 (Youngson 1983 \cite{Youngson83})
\item $JB^*$-triples, contractions\\
 (Kaup 1984 \cite{Kaup84}, Stacho 1982 \cite{Stacho82}, Friedman-Russo 1985 \cite{FriRus85JFA})
\item $\ell^p$-direct sums of $C_p(H)$, $1\le p<\infty$, $H$ Hilbert space, contractions\\
(Arazy-Friedman, 1978 \cite{AraFri78}, 2000 \cite{AraFri92}) 
\item $\ell^p$-direct sums of $L^p(\Omega,H)$, $1\le p<\infty$, $H$ Hilbert space, contractions\\
(Raynaud 2004 \cite{Raynaud04}) 
\item$\ell^p$-direct sums of  $C_p(H)$, $1\le p\ne 2<\infty$, $H$ Hilbert space, complete contractions\\
(LeMerdy-Ricard-Roydor 2009 \cite{LeMRicRoy09}) 
\end{itemize}

It follows immediately that  if $\mathcal S$ is projectively stable, then
so is the category $\mathcal S_*$  of spaces whose dual spaces belong to $\mathcal S$. 
 It should be noted that  $TROs,\ C^*$-algebras and $JC^*$-algebras are not
stable under contractive projections and $JB^*$-triples are not stable
under bounded projections.

\subsubsection{More about JB*-triples}
JB*-triples  are generalizations of JB*-algebras and C*-algebras.  The axioms can be said to come from geometry in view of  Kaup's Riemann mapping theorem  \cite{Kaup83}.  
Kaup showed in 1983 that
JB*-triples are exactly those Banach spaces whose open unit ball is a  
bounded symmetric domain.  
Kaup's holomorphic  
characterization of JB*-triples directly
led to the proof of the projective stability of JB*-triples  mentioned above.
\smallskip

Many authors  have studied the interplay between JB*-triples and  
infinite dimensional holomorphy.
Contractive projections have proved to be a valuable tool for the study of problems on JB*-triples (Gelfand-Naimark theorem \cite{FriRus86}, structure of inner ideals \cite{EdwMcCRut96}, operator space characterization of TROs \cite{NeaRus03PJM}, to name a few)
They are justified both as a natural
generalization of operator algebras as well as because of their connections
with complex geometry.\smallskip

Preduals of JBW*-triples have been called pre-symmetric spaces \cite{Edwards06}
and have been proposed as mathematical models of physical systems \cite{Friedman05}.
 In this model the operations on the physical system are represented by contractive projections on the pre-symmetric space. 
\smallskip

JB*-triples first arose in Koecher's  proof \cite{Koecher69}, \cite{Loos77}  of the classification of bounded symmetric domains in $\CC^n$. 
 The original proof of this fact, done in the 1930's by Cartan, used Lie algebras and Lie groups, techniques which do not extend to infinite dimensions. 
   On the other hand, to a large extent, the Jordan algebra techniques do so extend, as shown by Kaup and Upmeier.\footnote{The opposite is true concerning cohomology.   Lie algebra cohomology is well developed, Jordan algebra cohomology is not}

\subsubsection{Application: Gelfand Naimark theorem for JB*-triples}

\begin{theorem}[Friedman and Russo 1986 \cite{FriRus86}]
Every JB*-triple is isometically isomorphic to a subtriple of a direct sum of Cartan factors.
\end{theorem}

\smallskip

The theorem was not unexpected.  However, the proof required new techniques because of the lack of an order structure on a JB*-triple.

\medskip
\begin{itemize}
\item {\bf Step 1: February 1983 Friedman-Russo \cite{FriRus85JFA}}

Let $P:A\rightarrow A$ be a linear projection of norm 1 on a JC*-triple $A$. Then $P(A)$ is a JB*-triple under $\{xyz\}_{P(A)}=P(\{xyz\})$ for $x,y,z\in P(A)$.
\item {\bf Step 2: April 1983 Friedman-Russo \cite{FriRus84}}

Same hypotheses.  Then $P$ is a conditional expectation in the sense that
$$
P\{PaPbPc\}=P\{Pa,b,Pc\}$$  and $$P\{PaPbPc\}=P\{aPbPc\}.
$$
\item {\bf Step 3: May 1983 Kaup \cite{Kaup84}}

Let $P:U\rightarrow U$ be a linear projection of norm 1 on a JB*-triple $U$.  Then $P(U)$ is a JB*-triple under $\{xyz\}_{P(U)}=P(\{xyz\}_U)$ for $x,y,z\in P(U)$. 
  Also, $P\{PaPbPc\}=P\{Pa,b,Pc\}$ for $a,b,c\in U$, which extends one of the formulas in the previous step.
    \item {\bf Step 4: February 1984 Friedman-Russo \cite{FriRus85JRAM}}
  
  Every JBW*-triple splits into atomic and purely non-atomic ideals.
    \item {\bf Step 5: August 1984 Dineen \cite{Dineen84}}
  
  The bidual of a JB*-triple is a JB*-triple.
    \item {\bf Step 6: October 1984 Barton-Timoney \cite{BarTim86}}
  
  The bidual of a JB*-triple is a JBW*-triple, that is, the triple product is separately weak*-continuous.
    \item {\bf Step 7: December 1984 Horn \cite{Horn87}}
  
  Every JBW*-triple factor of type I is isomorphic to a Cartan factor.  
   More generally, every JBW*-triple of type I is isomorphic to a direct sum of $L^\infty$ spaces with values in a Cartan factor.
\item {\bf Step 8: March 1985 Friedman-Russo \cite{FriRus86}}

 {\bf Putting it all together}
$$
\pi:U\rightarrow U^{**}=A\oplus N=(\oplus_\alpha C_\alpha)\oplus N=\sigma(U^{**})\oplus N
$$
implies that $\sigma\circ \pi:U\rightarrow A=\oplus_\alpha C_\alpha$ is an isometric isomorphism.
\end{itemize}
\medskip

{\bf Consequences of the Gelfand-Naimark theorem}

\begin{itemize}
\item Every JB*-triple is isomorphic to a subtriple of a JB*-algebra.
\item In every JB*-triple, $\|\{xyz\}\|\le \|x\|\|y\|\|z\|$
\item Every JB*-triple contains a unique norm-closed ideal $J$ such that $U/J$ is isomorphic to a JC*-triple and $J$ is purely exceptional, that is, every homomorphism of $J$ into a C*-algebra is zero.
\end{itemize}

\subsubsection{Preservation of type}

The following results are due to Chu-Neal-Russo \cite{ChuNeaRus04}, and simultaneously, to Bunce-Peralta \cite{BunPer02}.  The corresponding (classical) result for von Neumann algebras, due to Tomiyama \cite{Tomiyama59}, required, of necessity, that the range of the projection be a subalgebra, since the category of von Neumann algebras and positive normal contractions is not projectively stable.

\begin{theorem}
Let $P$ be a normal contractive projection on a $JBW^*$-triple $Z$ of
type I.  Then $P(Z)$ is of
type I.
\end{theorem}

\begin{theorem}
Let $P$ be a normal contractive projection on a semifinite $JBW^*$-triple
$Z$. Then $P(Z)$ is a semifinite $JW^*$-triple.
\end{theorem}

\subsection{Projective rigidity}\label{proj2}

 By considering the converse of projective stability, one is lead to the
following definition.
\smallskip

  A
sub-category $\mathcal S$ of $\mathcal B$ is {\bf projectively rigid} if
it has the property that whenever $A$ is an object of $\mathcal S$ and $X$
is a subspace of $A$ which is isometric to an object in $\mathcal S$, then
$X$ is the range of an idempotent morphism of $\mathcal S$ on $A$. \smallskip

Examples of projectively rigid categories are the following.

\begin{itemize}
\item $\ell_p$, $1<p<\infty$, contractions\\
 (Pelczynski 1960 \cite{Pelczynski60})
\item $L^p$, $1\le p<\infty$, contractions\\
  (Douglas 1965 \cite{Douglas65}, Ando 1966 \cite{Ando66}, Bernau-Lacey 1974 \cite{BerLac74})
\item $C_p, 1\le p<\infty$, contractions \\
(Arazy-Friedman 1977 \cite{AraFri77})
\item Preduals of von Neumann algebras, \\
contractions (Kirchberg 1993 \cite{Kirchberg93})
\item Preduals of $TROs$, complete contractions (Ng-Ozawa 2002 \cite{NgOza02})
\item Preduals of JBW*-triples, contractions\footnote{with one exception, see Theorem~\ref{thm:0204121}}
 (Neal-Russo 2008 \cite{NeaRus11})
 \item$\ell^p$-direct sums of  $C_p(H)$, $1\le p\ne 2<\infty$, $H$ Hilbert space, complete contractions\\
(LeMerdy-Ricard-Roydor 2009 \cite{LeMRicRoy09}) 
\end{itemize}

\begin{theorem}[Neal-Russo]\label{thm:0204121}
The category of preduals of $JBW^*$-triples with no
summands of the form $L^1(\Omega,H)$
where  $H$ is a Hilbert space
of dimension at least two,  is projectively rigid.
\end{theorem}


 \subsection{Structural Projections}

\subsubsection{Structure of inner ideals}

In a series of papers, mainly by Edwards and R\"uttimann (cf.\ \cite{EdwRut96JLMS}), the ideal structure of 
C*-algebras and JB*-algebras has been thoroughly studied, and a surprising geometric characterization of the closed inner ideals among the closed subtriples of a JB*-triple has been established: a norm-closed subtriple B of a JB*-triple A is an inner ideal if and only if every bounded linear functional on B has a unique norm-preserving linear extension to A.

Loos and Neher  (\cite{Loos89},\cite{LooNeh94}) have introduced the notions of complementation and structural projection in the purely algebraic setting of Jordan pairs and Jordan *-triples. Let A be an anisotropic Jordan *-triple; then for every element a in A the quadratic mapping Q(a) is defined by Q(a)x: = \{axa\}. The kernel Ker(B) of a subset B in A is the subspace of A consisting of all elements annihilated by the mappings Q(b) as b ranges over B. A subtriple B is said to be complemented if $A = B\oplus \hbox{Ker}(B)$. A linear projection P on A is said to be structural if it satisfies Q(Pa) = PQ(a)P for all a in A, which occurs if and only if its range B = P(A) is a complemented subtriple of A.

In \cite{EdwMcCRut96}, the study of structural projections on JBW*-triples is continued.  It is shown  that a structural projection on a JBW*-triple is necessarily contractive and weak* continuous, and that every weak* closed inner ideal in a JBW*-triple A is a complemented subtriple of A and therefore the range of a unique structural projection.

\subsubsection{Geometric characterization}

Since structural projections are contractive and weak* continuous, every structural projection is the adjoint of a contractive linear projection on the predual  of A. The contractive projections P on $A_*$  that arise in this way also have the property that they are neutral, in that, if x is an element of  for which $ \|Px\|$ and $\|x\|$ coincide then Px and x also coincide. It follows from the results of \cite{EdwRut92} that the mapping $ P\mapsto P^*(A)$ is a bijection from the family  of neutral projections on $A_*$  for which $P^*A $ is a subtriple of A onto the complete lattice  of weak*-closed inner ideals in A.

The L-orthogonal complement $G^\diamond$  of a subset G of a complex Banach space E is the set of
elements x in E such that, for all elements y in G,
$\|x\pm y\|=\|x\|+\|y\|$:
A contractive projection P on E is said to be a GL-projection if the
L-orthogonal complement $PE^\diamond$ of the range PE of P is contained in the kernel  of P.
It is shown in \cite{EdwHugRut03} that, for a $JBW^*$-triple $A$ with pre dual $A_*$,
a linear projection $R$ on $A$ 
 is structural
 if and only if it is the adjoint of a neutral GL-projection on $A$, thereby giving a purely geometric characterization of structural projections.

 \subsubsection{Physical interpretation}

The predual  of a JBW*-triple A, also called a pre-symmetric space,  has been proposed as a model of the state space of a statistical physical system, in which contractive linear mappings on  $A_*$ represent operations or filters on the system.  For a discussion of this model we refer to \cite{EdwHugRut03} and for an application to decoherent states we refer to \cite{EdwHug08}.

In the theory of Banach spaces and in the study of state spaces of physical systems, much effort has been devoted to the investigation of when particular closed subspaces of a Banach space are the ranges of contractive projections. Since the range of a contractive projection on a pre-symmetric space is itself a pre-symmetric space, this problem is particularly relevant in the case of pre-symmetric spaces. The main results of \cite{EdwHug08} are used to determine under what circumstances the closed subspace of $A_*$ generated by a family  $\{x_j : j\in\Lambda\}$ of elements of $A_*$ of norm one which are pairwise decoherent (see \cite{EdwHug08} for the definition)  is the range of a contractive projection. (Theorem~\ref{thm:0204121} might be relevant here.)

\section{Hilbertian operator spaces}\label{sec:9}

\subsection{Operator spaces revisited}

We are going to quote some definitions and theorems from \cite{Pisier03}.
Operator space theory may also be considered as noncommutative Banach space theory. Although the subject had its genesis in the papers \cite{Stinespring55} and \cite{Arveson69} (announcement in \cite{Arveson69BAMS}), the first impetus appeared in the thesis of Ruan \cite{Ruan88}.
An operator space is a Banach space together with an isometric embedding in $B(H)$.  While the objects of this category are Banach spaces,  it is the morphisms, namely, the completely bounded maps, which are more important. 

An operator space is a subspace $X$ of $B(H)$.
  Its {\it
operator space structure} is given by the sequence of norms on
the set of matrices $M_n(X)$ with entries from $X$, determined by
the identification $M_n(X) \subset M_n(B(H))=B(H\oplus H \oplus
\cdots \oplus H)$. 
A linear mapping $\varphi:X\rightarrow Y$ between
two operator spaces is {\it completely bounded} if the induced
mappings $\varphi_n:M_n(X)\rightarrow M_n(Y)$ defined by
$\varphi_n([x_{ij}])=[\varphi(x_{ij})]$ satisfy
$\|\varphi\|_{\mbox{cb}}:=\sup_n\|\varphi_n\|<\infty$.

 Some simple examples of completely bounded maps are restrictions of $^*$-homomorphisms and multiplication operators.

\begin{theorem}
Every completely bounded map is the product of the above two mentioned maps.
\end{theorem}

Some special types of completely bounded maps which we will consider are complete contractions, complete isometries, and complete isomorphisms.  Additionally, a  complete semi-isometry is defined to be a complete contraction which is isometric.

\begin{remark} The space of completely bounded maps with the completely bounded norm
$CB(X,Y),\ \|\cdot \|_{\mbox{cb}}$ is a Banach space.
\end{remark}

\begin{definition} The completely bounded Banach-Mazur distance is defined by
\[
\mbox{d}_{\mbox{cb}}(E,F)=\inf\{\|u\|_{\mbox{cb}}\cdot
\|u^{-1}\|_{\mbox{cb}}
:u:E\rightarrow F\mbox{ complete isomorphism }\}.
\]
\end{definition}

Besides the pioneering work of Stinespring and Arveson mentioned above, many tools were developed in the 1970s and
1980s by a number of operator algebraists. As noted above, an abstract framework
was developed in 1988 in the thesis of Ruan. 
Besides \cite{Pisier03} we have already mentioned four other monographs on the subject (\cite{EffRua00},\cite{Paulsen02},\cite{BleLeM04},\cite{Helemskii10}).

\subsection{Classical operator spaces}
  ``Neoclassical'' and ``modern'' operator spaces will appear later.

Two important examples of Hilbertian operator spaces (:= operator
spaces isometric to Hilbert space) are the row and column spaces
 $R,\ C$, and their finite-dimensional versions $R_n,\ C_n$.
In $B(\ell_2)$,  {\it column Hilbert space}
$C:=\overline{\mbox{sp}}\{e_{i1}:i\ge 1\}$ and {\it row Hilbert
space} $R:=\overline{\mbox{sp}}\{e_{1j}:j\ge 1\}$.
 $R$ and $C$ are Banach isometric, but
not completely isomorphic:  $\mbox{d}_{\mbox{cb}}(R,C)= \infty$.
 $R_n$ and $C_n$ are completely isomorphic,
 but not completely isometric: $\mbox{d}_{\mbox{cb}}(R_n,C_n)= n$.

$R,\ C,\ R_n,\ C_n$ are examples of {\it homogeneous} operator
spaces, that is, operator spaces $E$ for which $\forall
u:E\rightarrow E$, $\|u\|_{\mbox{cb}}=\|u\|$.
Another important
example of an Hilbertian homogeneous operator space is $\Phi(I)$.
$\Phi(I)=\overline{\mbox{sp}}\{V_i:i\in I\}$, where the $V_i$ are
bounded operators on a Hilbert space satisfying the canonical
anti-commutation relations.  In some special cases, the notations
 $\Phi_n:=\Phi(\{1,2,\ldots,n\})$, and
$\Phi=\Phi(\{1,2,\ldots \})$ are used.

 If $E_0\subset B(H_0)$ and $E_1\subset B(H_1)$ are
operator spaces whose underlying Banach spaces form a compatible
pair in the sense of interpolation theory, then the Banach space
$E_0\cap E_1$, with the norm $$\|x\|_{E_0\cap E_1}=\max
(\|x\|_{E_0},\|x\|_{E_1})$$ equipped with the operator space
structure given by the embedding $E_0\cap E_1 \ni x\mapsto (x,x)\in
E_0\oplus E_1\subset B(H_0\oplus H_1)$ is called the {\it
intersection} of $E_0$ and $E_1$ and is denoted by $E_0\cap E_1$.
Examples are  $R\cap C,\ \Phi=\cap_0^\infty H_\infty^{m,L},\  \Phi_n =\cap_1^n H_n^k.
$ $
(H_\infty^{m,L}$ and $H_n^k$ are defined below.)
The definition of intersection
extends easily to arbitrary families of compatible operator spaces
 
We already noted that 
$R,\ C,\ R_n,\ C_n$ are homogenous operator spaces.  So are
$\min(E)$ and $\max(E)$
and
$\Phi(I)$, where  $\min(E)$ is defined by
 $E\subset C(T)\subset B(H)$,
$\|(a_{ij})\|_{M_n(\min(E))}=\sup_{\xi\in
B_{E^*}}\|(\xi(a_{ij}))\|_{M_n}$
and
 $\max(E)$ is defined by
$\|(a_{ij})\|_{M_n(\max(E))}=\sup\{\|(u(a_{ij}))\|_{M_n(B(H_u))}:u:E\rightarrow
B(H_u),\ \|u\|\le 1 \}$. These two spaces have a stronger property:
$$F\stackrel{u}{\longrightarrow}\min(E)\Rightarrow \rightarrow \|u\|_{\mbox{cb}}=\|u\|
\hbox{ and }\max(E)\stackrel{v}{\longrightarrow}G\Rightarrow
 \|v\|_{\mbox{cb}}=\|v\|,$$
 and the identity map $\min(E)\rightarrow \max(E)$ is completely contractive.\smallskip

The following is an application of the Russo-Dye theorem \cite{RusDye66}.

\begin{proposition} Let $E$ be a Hilbertian operator space.  Then $E$ is homogeneous if and only if 
$\|U\|_{\mbox{cb}}=1$ $\forall \mbox{ unitary } U:E\rightarrow E$.
\end{proposition}

The following are considered to be among the ``classical''  Banach spaces
$\ell_p,\ c_0,\ L_p,\ C(K)$. A ``second generation'' of such spaces would contain Orlicz, Sobolev, and Hardy spaces, the disk algebra, and the Schatten $p$-classes.  By analogy, some classical operator spaces are $R,\ C,\ \min(\ell_2),\ \max(\ell_2),\ OH,\ \Phi$ and their finite dimensional versions
$R_n,\ C_n,\ \min(\ell_2^n),\ \max(\ell_2^n),\ OH_n,\ \Phi_n$. Note that the classical operator spaces are all Hilbertian.

 \begin{proposition}
 The classical operator spaces are mutually completely non-isomorphic.  If $E_n,F_n$ are $n$-dimensional versions, then
$\mbox{d}_{\mbox{cb}}(E_n,F_n)\rightarrow\infty$.
\end{proposition}

\subsection{Neoclassical operator spaces}\label{sec:proj}

The ``neoclassical'' operator spaces $H_n^k$ appeared in \cite{NeaRus03TAMS} (announcement in \cite{NeaRus00}), from which the following three theorems were taken.

\begin{theorem}
 There is a family of 
1-mixed injective Hilbertian operator spaces 
$H_n^k$, $1\le k\le n$, of finite dimension $n$,
with the following properties:
\begin{description}
\item[(a)] $H_n^k$ is a  subtriple  of the Cartan factor of type 1 
consisting of
all $n\choose k$ by $n\choose n-k+1$ complex matrices.
\item[(b)] Let $Y$ be a \jwst\ of rank 1 {\rm (}necessarily 
atomic{\rm )}.
\begin{description}
\item[(i)] If $Y$ is of finite dimension $n$ then it is isometrically
completely contractive to some
$H_n^k$.  
\item[(ii)] If $Y$ is infinite dimensional then it is isometrically completely
contractive to $B(H,\CC)$ or $B(\CC,K)$.
\end{description} 
\item[(c)] $H_n^n$ {\rm (}resp. $H_n^1${\rm )}
coincides with $R_n$ {\rm (}resp. $C_n${\rm )}.
\item[(d)] For $1<k<n$, $H_n^k$ is not completely semi-isometric to $R_n$
or $C_n$.
\end{description}
\end{theorem}

Example 1: 
$H_3^2$ is the subtriple of $B(\CC^3)$ consisting of all matrices of the
form
\[
\left[\begin{array}{rrr}
0&a&-b\\
-a&0&c\\
b&-c&0
\end{array}\right]
\]
and is  completely semi-isometric to the 
Cartan factor $A(\CC^3)$ of 3 by 3 anti-symmetric complex matrices.
\smallskip

Example 2: 
$H_4^3$ is the
subtriple of $B(\CC^6,\CC^4)$ consisting of all matrices of the form
\[
\left[\begin{array}{rrrrrr}
0&0&0&-d&c&-b\\
0&d&-c&0&0&a\\
-d&0&b&0&-a&0\\
c&-b&0&a&0&0
\end{array}\right].
\]

\begin{theorem}
Let $X$ be a 1-mixed injective operator space which is atomic.
Then $X$ is completely semi-isometric 
to a direct sum of Cartan factors of types 1 to
4 and the spaces $H_n^k$.
\end{theorem}

\begin{theorem}
Let $Y$ be an 
atomic w$^*$-closed $JW^*$-subtriple of a $W^*$-algebra.
\begin{description}
\item[(a)] If $Y$ is irreducible and
 of rank at least 2, then it  is completely isometric
to a Cartan factor of type 1--4 or the
space \newline $\mbox{Diag}\, (B(H,K),B(K,H))$.
\item[(b)] If $Y$ is of finite dimension $n$
 and of rank 1, then it is completely
isometric to \newline 
$\mbox{Diag}\, (H_n^{k_1},\ldots,H_n^{k_m})$, for appropriately
chosen bases, and where
$k_1>k_2>\cdots>k_m$.
\item[(c)] $Y$ is completely semi-isometric to a direct sum of the spaces in (a) and (b).
If $Y$ has no infinite dimensional rank 1 summand, then it is completely
isometric to a direct sum of the spaces in (a) and (b).
\end{description}
\end{theorem}

We now describe, in more detail, the spaces $H_n^k$, including the next two theorems,  \cite{NeaRus06PAMS}.

  A frequently mentioned  result of Friedman and Russo  states that if a
subspace $X$ of a C*-algebra $A$  is the range of a contractive
projection on $A$, then $X$ is isometric to a JC*-triple, that is,
a norm closed subspace of $B(H,K)$ stable under the triple product
$ab^*c+cb^*a$.
If $X$ is atomic (in particular, finite-dimensional), then it is
isometric to a direct sum of Cartan factors of types 1 to 4.
This latter result fails, as it stands, in the category of
operator spaces.

Nevertheless, there exists a family of $n$-dimensional Hilbertian operator
spaces $H_n^k$, $1\le k\le n$, generalizing the row and column
Hilbert spaces $R_n$ and $C_n$ such that, in the above result, if
$X$ is atomic, the word ``isometric'' can be replaced by
``completely semi-isometric,'' provided the spaces $H_n^k$ are
allowed as summands along with the Cartan factors.
The space $H_n^k$ is contractively complemented in some $B(K)$,
and for $1<k<n$, is not completely (semi-)isometric to either of
the Cartan factors $B(\CC,\CC^n)=H_n^1$ or $B(\CC^n,\CC)=H_n^n$.
These spaces appeared in a slightly different form and context in
a memoir of Arazy and Friedman \cite{AraFri78}.

The construction of  $H_n^k$ is as follows.
Let $I$ denote a subset of $\{1,2,\ldots,n\}$ of cardinality
$|I|=k-1$. The number of such $I$ is $q:={n\choose k-1}$.
Let $J$ denote a subset of $\{1,2,\ldots,n\}$ of cardinality
$|J|=n-k$. The number of such $J$ is $p:={n\choose n-k}$.
 The
space $H_n^k$ is the linear span of matrices $b_i^{n,k}$, $1\le
i\le n$, given by
\[
b_i^{n,k}=\sum_{I\cap J=\emptyset,(I\cup
J)^c=\{i\}}\epsilon(I,i,J)e_{J,I},
\]
where $e_{J,I}=e_J\otimes e_I=e_Je_I^t\in
M_{p,q}(\CC)=B(\CC^q,\CC^p)$, and $\epsilon(I,i,J)$  is the
signature of the permutation taking
$(i_1,\ldots,i_{k-1},i,j_1,\ldots,j_{n-k})$ to $(1,\ldots,n)$.
Since the $b_i^{n,k}$ are the image under a triple isomorphism
(actually ternary isomorphism) of a rectangular grid in a
JC*-triple of rank one, they form an orthonormal basis for
$H_n^k$.

\begin{theorem}
 $H_n^k$ is a homogeneous operator space.
\end{theorem}

\begin{theorem}\label{thm:0219121}
$\mbox{d}_{\mbox{cb}}(H_n^k,H_n^1)=\sqrt{\frac{kn}{n-k+1}}$, for
$1\le k\le n$.
\end{theorem}
 
The latter theorem is proved using the relation with creation operators.
Let $ C^{n,k}_{h}$ denote the wedge (or creation) operator from
$\wedge^{k-1}\CC^n$ to $\wedge^{k}\CC^n$ given by
$$C^{n,k}_{h}(h_1\wedge\cdots\wedge h_{k-1})=h\wedge
h_1\wedge\cdots\wedge h_{k-1}.
$$

Letting ${\mathcal C}^{n,k}$ denote the space
$\mbox{sp}\{C_{e_i}^{n,k}\}$, we have the following.

\begin{proposition}
 $H_n^k$ is completely isometric to ${\mathcal
C}^{n,k}$.
\end{proposition}

\subsection{Modern operator spaces}
 
This subsection is a description of the results of \cite{NeaRicRus06}, which considers the following two problems, leading to the notion of ``modern operator space.''

\begin{itemize}
\item  Classify all infinite
dimensional rank 1 JC*-triples up to complete isometry ({\bf algebraic})\\
ANSWER:  $\Phi,\quad H_\infty^{m,R},\quad H_\infty^{m,L},\quad H_\infty^{m,R}\cap H_\infty^{m,L}$
\item Give a
suitable ``classification'' of all Hilbertian operator spaces which
are contractively complemented in a C*-algebra or normally
contractively complemented in a W*-algebra  ({\bf analytic})\\
ANSWER: $\Phi,\quad {\bf C},\quad {\bf R},\quad {\bf C}\cap {\bf R}$
\end{itemize}

It is known that  ${\bf R}$ and ${\bf C}$ are the only  {\bf completely} contractively complemented Hilbertian operator spaces \cite{Robertson91}.

\subsubsection{Rank one JC*-triples}\label{TRO}

In order to attack the operator space structure theory of rank 1 JC*-triples, let us first review some of the basic concepts about them, as well as those concerning ternary rings of operators (TROs).

A \jcst\ is a norm closed complex linear subspace  of $B(H,K)$
(equivalently, of a $C^*$-algebra) which is closed under the
operation $a\mapsto aa^*a$.  \jcst s were defined and studied (using
the name $J^*$-algebra)
 as a generalization of \csa s by Harris \cite{Harris74}  in
connection with function theory on infinite dimensional bounded
symmetric domains. 
 By a polarization identity (involving $\sqrt{-1}$), any \jcst\ is closed
under the triple product
\begin{equation}\label{eq:product}
(a,b,c)\mapsto \tp{a}{b}{c}:=\frac{1}{2}(ab^*c+cb^*a),
\end{equation}
under which it becomes a Jordan triple system.  
 A linear map which
preserves the triple product (\ref{eq:product}) will be called a
{\it triple homomorphism}.
  Cartan factors are examples of \jcst s,
as are \csa s, and Jordan \csa s.
 We shall only make use of
Cartan factors of type 1, that is, spaces of the form $B(H,K)$ where
$H$ and $K$ are complex Hilbert spaces.

A special case of a \jcst\ is a {\it ternary algebra}, that is, a
subspace of $B(H,K)$ closed under the {\it ternary product}
$(a,b,c)\mapsto ab^*c$. 
 A {\it ternary homomorphism} is a linear map
$\phi$ satisfying $\phi(ab^*c)=\phi(a)\phi(b)^*\phi(c)$. 
 These
spaces are also called ternary rings of operators and abbreviated
TRO. 
 TROs have come to play a key role in operator space theory, serving
as the algebraic model in the category.
(The algebraic
models for the categories of order-unit spaces, operator systems,
and Banach spaces, are respectively Jordan $C^*$-algebras, \csa s,
and \jbst s.) 
For TROs, a ternary isomorphism is the same as
a complete isometry.

 Every \jwst\ of
rank one is isometric to a Hilbert space and every maximal collinear
family of partial isometries corresponds to an orthonormal basis.
 Conversely, every Hilbert space with the abstract triple product
$\tp{x}{y}{z}:=(\ip{x}{y}z+\ip{z}{y}x)/2$ can be realized as a
\jcst\ of rank one in which every orthonormal basis forms a maximal
family of mutually collinear minimal partial isometries.
Collinearity of $v$ and $w$ means: $$vv^*w+wv^*v=w\hbox{ and }ww^*v+vw^*w=v$$

 \subsubsection{Operator space structure of Hilbertian JC*-triples}
We shall now outline the proof  from \cite{NeaRicRus06} of the classification of infinite dimensional Hilbertian operator spaces up to complete isometry .

\begin{itemize}

\item The general setting:
$Y$ is a $JC^*$-subtriple of $B(H)$ which is Hilbertian in the
operator space structure arising from $B(H)$, and
$\{u_i:i\in\Omega\}$ is an orthonormal basis consisting of a maximal
family of mutually collinear partial isometries of $Y$. 

\smallskip

\item We let $T$ and $A$ denote the TRO and the $C^*$-algebra respectively
generated by $Y$. For any subset $G\subset\Omega$,
$(uu^*)_G:=\prod_{i\in G}u_iu_i^*$ and $(u^*u)_G:=\prod_{i\in
G}u_i^*u_i$. The elements $(uu^*)_G$ and $(u^*u)_G$ lie in the weak
closure of $A$ and more generally in  the left and right linking von
Neumann algebras of $T$.

Fix $m\ge 0$.  To construct $H_\infty^{m,R}$ we make a temporary  assumption on the ternary envelope of $Y$.

\smallskip

\item {\bf Assume} $(u^*u)_G\ne 0$ for $|G|\le m+1$ and $(u^*u)_G=0$ for
$|G|\ge m+2$.

\smallskip

\item Define elements which are indexed by an arbitrary pair of subsets $I,J$ of
$\Omega$ satisfying
\begin{equation}\label{eq:12}
|\Omega-I|=m+1,\ |J|=m,\ 
\end{equation} 
 as follows: 
 
 $u_{IJ}=$\ 
$
(uu^*)_{I- J}u_{c_1}u_{d_1}^*u_{c_2}u_{d_2}^*\cdots
u_{c_s}u_{d_s}^*u_{c_{s+1}}(u^*u)_{J- I}, \ 
$
where  

 \medskip

$I\cap J=\{d_1,\ldots,d_s\}$  and
 $ (I\cup J)^c=\{c_1,\ldots,c_{s+1}\}$.

\smallskip

\item After an appropriate choice $\epsilon(IJ)$ of signs, the map $\epsilon(IJ)u_{IJ}
\rightarrow E_{JI}$ is a ternary isomorphism (and hence complete
isometry) from the norm closure of $\mbox{sp}_{C}\, u_{IJ}$ to the
norm closure of $\mbox{sp}_{C}\, \{E_{JI}\}$, where $E_{JI}$ denotes
an elementary matrix, whose rows and columns are indexed by the sets
$J$ and $I$, with a 1 in the $(J,I)$-position. 

\smallskip

\item This map can be extended to a ternary isomorphism
from the w*-closure of $\mbox{sp}_{C}\, u_{IJ}$ onto the Cartan
factor of type I consisting of all $\aleph_0$ by $\aleph_0$ complex
matrices which act as bounded operators on $\ell_2$. By restriction, $Y$ is completely isometric to a
subtriple $\tilde Y$, of this Cartan factor of type 1.

\item
\begin{definition}
We shall denote the space $\tilde Y$ above by $H_\infty^{m,R}$.  
\end{definition}

\item An
entirely symmetric argument (with $J$ infinite and $I$ finite) under
an entirely symmetric assumption on $Y$  defines the space
$H_\infty^{m,L}$.

\item Having constructed the spaces $H_\infty^{m,R}$ and $H_\infty^{m,L}$, we can  now state the following theorem.

\end{itemize}

\begin{theorem}
Let $Y$ be a $JC^*$-subtriple of $B(H)$ which is a separable infinite dimensional Hilbertian operator
space. Then $Y$
is completely isometric to one of the following spaces:
\[
\Phi,\quad H_\infty^{m,R},\quad H_\infty^{m,L},\quad H_\infty^{m,R}\cap H_\infty^{n,L}.
\]
\end{theorem}

\subsubsection{Further properties of $H_\infty^{m,R}$ and $H_\infty^{m,L}$}

We now describe, from \cite{NeaRicRus06}  the relation of these spaces with Fock space and compute some completely bounded Banach-Mazur distances.

If $H$ is a separable  Hilbert space, $ l_m(h):H^{\wedge m}\rightarrow H^{\wedge m+1}$   is the creation operator
$l_m(h)x=h\wedge
x.
$
Creation operators form a linear space
${\mathcal
C}^{m}=\overline{\mbox{sp}}\{l_m(e_i)\}$,  where $\{e_i\}$  is an orthonormal basis, and inherit an operator space structure from $B(H^{\wedge m},H^{\wedge m+1})$. 
Annihilation operators ${\mathcal
A}^{m}$ consists of the adjoints of the creation operators on $H^{\wedge m-1}$.

\begin{lemma}
 $H_\infty^{m,R}$ is completely isometric to ${\mathcal A}^{m+1}$ and $H_\infty^{m,L}$ is completely isometric to  ${\mathcal C}^{m}$.
\end{lemma}

\begin{remark}
Every finite or infinite dimensional separable Hilbertian JC*-subtriple $Y$ is completely isometric to a finite or infinite intersection of spaces of creation and annihilation operators, as follows
\begin{description}
\item[(a)] If $Y$ is infinite dimensional, then it is completely isometric to one of
  $
 {\mathcal A}^m,\quad {\mathcal C}^m,\quad  {\mathcal A}^m\cap {\mathcal C}^k,\quad \cap_{k=1}^\infty {\mathcal C}^k 
 $
 \item[(b)] If $Y$ is of dimension $n$, then  $Y$ is completely isometric to $\cap_{j=1}^m {\mathcal C}^{k_j}$, where $n\ge k_1>\cdots>k_m\ge 1$
\end{description}
 \end{remark}

\begin{theorem}
For $m$, $k\ge 1$,
\begin{description}
\item[(a)] $d_{cb}(H^{m,R}_\infty,H^{k,R}_\infty)=
d_{cb}(H^{m,L}_\infty,H^{k,L}_\infty)=\sqrt{\frac {m+1}{k+1}}$ when $m\ge k$\item[(b)] $d_{cb}(H^{m,R}_\infty,H^{k,L}_\infty)=\infty$
\item[(c)] $d_{cb}(H^{m,R}_\infty,\Phi)=d_{cb}(H^{m,L}_\infty,\Phi)=\infty$
\end{description}
\end{theorem}


\subsubsection{Contractively complemented Hilbertian operator spaces}

\begin{theorem}
Suppose $Y$ is a separable infinite dimensional Hilbertian operator
space which is contractively complemented (resp. normally
contractively complemented) in a \csa\ $A$ (resp. W*-algebra $A$) by
a projection $P$. 
Then,  replacing $Y$  with its support (defined in \cite{NeaRicRus06}),
 $Y$  is completely isometric to either $R$, $C$, $R\cap C$, or
$\Phi$.
\end{theorem}

\begin{theorem}[Converse]
The operator spaces $R,C,R \cap C$, and  $\Phi$ are each
essentially (defined in \cite{NeaRicRus06}) normally contractively complemented in a von Neumann
algebra.
\end{theorem}

\section{Quantum operator algebras}

\subsection{Operator space characterization of TROs}

A natural object to characterize in the context of operator spaces are
the so called {\it ternary rings of operators} (TRO's).
These are subspaces of $B(H)$ which
are closed under the ternary
product $xy^{\ast}z$. TRO's, like C*-algebras, carry a natural operator space structure.
In fact, every TRO is (completely)
 isometric to a corner
$pA(1-p)$ of a C*-algebra $A$. TRO's are important because, as shown by
Ruan \cite{Ruan89}, the injectives in the category of
operator spaces are TRO's (corners of injective C*-algebras).
(For the dual version of this result see
\cite{EffOzaRua01}.)
Injective envelopes of operator systems and
of operator spaces (\cite{Hamana79},\cite{Ruan89})
 have proven to be important tools,
see for example \cite{BlePau01}. The
characterization of TRO's among operator spaces is the subject of this subsection.

Closely related to TRO's are the so called JC*-triples, norm closed
subspaces of $B(H)$
which are closed under the triple product
$(xy^{\ast}z + zy^{\ast}x)/2$.
 These generalize the class of TRO's and
have the property, as shown by Harris in \cite{Harris74},
that isometries coincide with algebraic isomorphisms. It is not hard to
see this implies that the algebraic isomorphisms in the
class of TRO's are complete isometries, since for each TRO
$A$, $M_{n}(A)$ is a JC*-triple (For the converse of this, see
\cite[Proposition 2.1]{Hamana99}).
As a consequence, if an operator space $X$ is completely isometric
to a TRO, then the induced ternary product on $X$ is unique, {\it i.e.},
independent of the TRO.

Relevant to this subsection is another property shared by all JC*-triples (and
hence all TRO's). For any Banach space $X$, we denote by $X_0$ its
open unit ball: $\{x\in X:\|x\|<1\}$.
The open unit ball of every JC*-triple is a {\it bounded
symmetric domain}. This is equivalent to saying that it has a
transitive group of biholomorphic automorphisms. It was shown by Koecher
in finite dimensions (see \cite{Loos77}) and Kaup
\cite{Kaup83} in the general case that this is a defining property for the
slightly larger class of JB*-triples. The only
illustrative basic examples of  JB*-triples which are not  JC*-triples
are the
space $H_{3}({\mathcal O})$ of 3 x 3 hermitian matrices
over the octonians and a certain subtriple of $H_{3}({\mathcal O})$. These
are called  {\it exceptional} triples, and they
cannot be represented as a JC*-triple. This holomorphic characterization
has been useful because, as noted earlier,  it gives an elegant proof, due to Kaup
\cite{Kaup84}, that the range of a contractive projection on a JB*-triple
is isometric to another JB*-triple. T

Motivated by this characterization for JB*-triples, a
holomorphic characterization of TRO's
up to complete isometry is given in  \cite{NeaRus03PJM} and stated below. 
 As a
consequence, a holomorphic operator
space characterization of C*-algebras is obtained as well. It should be mentioned that
Upmeier (for the category of Banach spaces) in \cite{Upmeier84} and El Amin-Campoy-Palacios
(for the category of Banach algebras) in \cite{ElACamPal01},
gave  different but still
holomorphic characterizations of C*-algebras up to isometry.

\begin{theorem}\label{thm:0220121}
Let $A\subset B(H)$ be an operator space and suppose that $M_n(A)_0$ is a 
bounded symmetric domain for all $n\ge 2$, then
$A$ is ternary isomorphic and completely isometric to a TRO.
\end{theorem}


\begin{theorem}
Let $A\subset B(H)$ be an operator space and suppose that  $M_n(A)_0$ is a bounded symmetric domain
for all $n\ge 2$ and $A_0$ is of tube type. Then 
$A$ is completely isometric to a C*-algebra.
\end{theorem}

\subsection{Holomorphic characterization of operator algebras}

In the category of operator spaces, that is,
subspaces of the bounded linear operators $B(H)$
on a complex Hilbert space $H$ together with the
induced matricial operator norm structure, objects are
equivalent if they are completely isometric, {\it i.e.} if there is a linear
isomorphism between the spaces which preserves this
matricial norm structure. Since operator algebras, that is,
subalgebras of $B(H)$, are  motivating
examples for much of operator space theory, it is natural
to ask if one can characterize which  operator spaces are
operator algebras. One satisfying answer was given by Blecher,
Ruan and Sinclair in \cite{BleRuaSin90}, where it was shown that among operator
spaces $A$ with a (unital but not necessarily associative)
Banach algebra product, those which are
completely isometric to operator algebras are precisely the ones whose
multiplication is completely contractive with respect to the
Haagerup norm on $A\otimes A$
(For a completely bounded version of this
result, see \cite{Blecher95}).

The main result of \cite{NeaRus12} drops the algebra assumption  on $A$  in the Blecher-Ruan-Sinclair theorem   in favor of a holomorphic assumption. Using only natural conditions on holomorphic vector fields on Banach spaces,  an algebra product is constructed on $A$ which is completely contractive and unital, so that the result of \cite{BleRuaSin90} can be applied. 

We state this result as  Theorem~\ref{thm:0221121}. In this theorem, for any element $v$ in the symmetric part  (recalled below) of a Banach space $X$, $h_v$ denotes the corresponding complete holomorphic vector field on the open unit ball of $X$.

\begin{theorem}\label{thm:0221121}
An operator space $A$ is completely isometric to a unital operator
algebra if and only there exists an element $v$ in the completely symmetric part of $A$ such that:
\begin{enumerate}
\item $h_v(x+v)-h_v(x)-h_v(v)+v=-2x$ for all $x\in A$
\smallskip
\item  For all $X\in M_n(A)$,
$
\|V-h_{V}(X)\|\le \|X\|^2.
$
\end{enumerate}
\end{theorem}
\smallskip

This result  is thus an instance where the consideration of  a ternary product, called the partial triple product, which arises from the holomorphic structure via the symmetric part of the Banach space, leads to results for binary products.  Examples of this phenomenon occurred in  \cite{ArazyvN}, \cite{AraSol90} where this technique is used to describe the algebraic properties of  isometries of certain operator algebras.  The technique was also used in 
\cite{KauUpm76} to show that Banach spaces with holomorphically equivalent unit balls are linearly isometric (see \cite{Arazysurvey} for an exposition of \cite{KauUpm76}).

Suppose that $X$ is a TRO (i.e., a closed subspace of $B(H)$ closed under the ternary product $ab^*c$) which contains an element $v$ satisfying $xv^*v=vv^*x=x$ for all $x\in X$.  Then it is trivial that $X$ becomes a unital C$^*$-algebra for the product $xv^*y$, involution $vx^*v$, and unit $v$.  By comparison, the above result starts only with an operator space $X$ containing a distinguished element $v$ in the completely symmetric part of $X$ (defined below) having a unit-like property.  This is to be expected since the result of \cite{BleRuaSin90} fails in the absence of a unit element. 

The main technique in the proof of this result is to use a variety of elementary  isometries on $n$ by $n$ matrices over $A$ (most  of the time, $n=2$) and to exploit the fact that isometries of arbitrary Banach spaces preserve the partial triple product. The first occurrence  of this technique appears in the construction, for each $n$, of a contractive projection $P_n$ on $\kbta$ ($K$= compact operators on separable infinite dimensional Hilbert space) with range $M_n(A)$  as a convex combination of isometries. The completely symmetric part of $A$ is defined  to be the intersection of $A$ (embedded in $\kbta$) and the symmetric part of $\kbta$. It is  the image under $P_1$ of the symmetric part of $\kbta$. It follows from \cite{NeaRus03PJM} that the completely symmetric part of $A$ is a TRO, which is a crucial tool in this work. 
The binary product $x\cdot y$  on $A$  is constructed  using  properties of isometries on 2 by 2 matrices over $A$ and it is shown that the symmetrized product can be expressed in terms of the partial Jordan triple product  on $A$.  Specifically,
\[
\left[\begin{array}{cc}
y\cdot x& 0\\
0&0
\end{array}
\right]=2
\tp{\left[\begin{array}{cc}
x& 0\\
0&0
\end{array}
\right]}{\left[\begin{array}{cc}
0&v\\
0&0
\end{array}
\right]}{\left[\begin{array}{cc}
0&y\\
0&0
\end{array}
\right]}
\]
and \[
\tp{x}{v}{y}=\frac{1}{2}(y\cdot x+x\cdot y).
\]

According to \cite{BleZarPNAS},  ``The one-sided multipliers of an operator space $X$ are a key to the `latent operator algebraic structure' in $X$.''  The unified approach through multiplier operator algebras developed in \cite{BleZarPNAS} leads to simplifications of known results and applications to quantum $M$-ideal theory.  They also state
``With the extra structure consisting of the additional matrix norms on an operator algebra, one might expect to not have to rely  as heavily on other structure, such as the product.''  Theorem~\ref{thm:0221121} is certainly in the spirit of this statement.

Another approach to operator algebras is \cite{Kaneda2007}, in which the set of operator algebra products on an operator space is shown to be  in bijective correspondence with the space of norm one quasi-multipliers on the operator space.

\subsubsection{Symmetric part of a Banach space}

We review the construction and properties of the partial Jordan triple product in an arbitrary Banach space.  Let $X$ be a complex Banach space with open unit ball $X_0$. Every holomorphic function $h:X_0\rightarrow X$, also called a holomorphic vector field,  is locally integrable, that is, the initial value problem
\[
\frac{\partial}{\partial t}\varphi(t,z)=h(\varphi(t,z))\ ,\ \varphi(0,z)=z, 
\]
has a unique solution for every $z\in X_0$ for $t$ in a maximal open interval $J_z$ containing 0.  A {\bf complete holomorphic vector field} is one for which $J_z=(-\infty,\infty)$ for every $z\in X_0$.

It is a fact that every complete holomorphic vector field is the sum of the restriction of a skew-Hermitian bounded linear operator $A$ on $X$
and a function $h_a$ of the form
$
h_a(z)=a-Q_a(z),
$
where $Q_a$ is a quadratic homogeneous polynomial on $X$.  

The {\bf symmetric part} of $X$ is the orbit of 0 under the set of complete holomorphic vector fields, and is denoted by  $S(X)$.  It is a closed subspace of $X$ and is equal to $X$ precisely when $X$ has the structure of a $JB^*$-triple
(by \cite{Kaup83}).

If $a\in S(X)$, we can obtain a symmetric bilinear form on $X$, also denoted by $Q_a$ via the polarization formula
\[
Q_a(x,y)=\frac{1}{2}\left(Q_a(x+y)-Q_a(x)-Q_a(y)\right)
\]
and then the partial Jordan triple product  $\{\cdot,\cdot,\cdot\}:X\times S(X)\times X\rightarrow X$ is defined by $\{x,a,z\}=Q_a(x,z)$.
The space $S(X)$ becomes a $JB^*$-triple in this triple product.

It is also true that the ``main identity'' (\ref{eq:0221121})  (see subsection~\ref{5.5}) holds whenever $a,y,b\in S(X)$ and $x,z\in X$.
The partial triple product behaves well under the actions of isometries and contractive projections
 (see \cite[5.2,5.3]{Arazysurvey}).

\section{Universal Enveloping TROs and C*-algebras}
\subsection{Operator space structure of $JC^*$-triples and TROs}

The paper \cite{BunFeeTim11}, together with its sequel \cite{BunTim11}, initiates a systematic investigation of the operator space structure of $JC^*$-triples via a study of the TROs they generate.  The approach is through the introduction and development of a variety of universal objects (TROs and $C^*$-algebras). Explicit descriptions of operator space structures of Cartan factors of arbitrary dimension are obtained as a consequence.  It should be noted that Theorem~\ref{thm:0222123} and its corollary, as well as the computation of the universal enveloping TROs of all Cartan factors of finite dimension, were also obtained roughly around the same time in \cite{BohWer11}.  Much of this subsection is taken literally from the introduction to \cite{BunFeeTim11}.

It is necessary to widen the definition of $JC^*$-triple to include any complex Banach space with a surjective isometry onto a subspace of a TRO which is closed under the (Jordan) triple product $\{a,b,c\}=(ab^*c+cb^*a)/2$.
A contractively complemented subspace of $B(H)$ is an example of such a $JC^*$-triple, not necessarily a $JC^*$-subtriple  of $B(H)$.  In the reflexive case, operator space structure of $JC^*$-triples of this kind has been investigated  in \cite{NeaRus03TAMS}, \cite{NeaRus06PAMS}, and \cite{NeaRicRus06}.  

The aim of \cite{BunFeeTim11} is to explore the operator space structures of an arbitrary $JC^*$-triple, and, quoting the authors of \cite{BunFeeTim11}, ``to place the ground-breaking results of \cite{NeaRus03TAMS}, \cite{NeaRus06PAMS},  \cite{NeaRicRus06} in the general setting of a full theory.''  This new device (universal TRO of a $JC^*$-triple) facilitates a general inquiry into the operator space structure of $JC^*$-triples.

\begin{theorem}[Theorem 3.1 of \cite{BunFeeTim11}]\label{thm:0222123}
Let $E$ be a $JC^*$-triple.  Then there is a unique pair $(C^*(E),\alpha_E)$ where $C^*(E)$ is a $C^*$-algebra and $\alpha_E:E\rightarrow C^*(E)$ is an injective triple homomorphism, with the following properties:

(a) $\alpha_E(E)$ generates $C^*(E)$ as a $C^*$-algebra;

(b) for each triple homomorphism $\pi:E\rightarrow A$, where $A$ is a $C^*$-algebra, there is a unique $*$-homomorphism $\tilde\pi:C^*(E)\rightarrow A$ with $\tilde\pi\circ\alpha_E=\pi$.
\end{theorem}

\begin{corollary} [Corollary 3.2 of \cite{BunFeeTim11}]
Let $E$ be a $JC^*$-triple.  Then there is a unique pair $(T^*(E),\alpha_E)$ where $T^*(E)$ is a TRO and $\alpha_E:E\rightarrow T^*(E)$ is an injective triple homomorphism, with the following properties:

(a) $\alpha_E(E)$ generates $T^*(E)$ as a TRO;

(b) for each triple homomorphism $\pi:E\rightarrow T$, where $T$ is aTRO, there is a unique TRO-homomorphism $\tilde\pi:T^*(E)\rightarrow T$ with $\tilde\pi\circ\alpha_E=\pi$.
\end{corollary}

In a technical tour de force, $T^*(E)$ is explicitly calculated for the Cartan factors of any dimension: Hilbert space \cite[Theorem 5.1]{BunFeeTim11}, spin factor \cite[Lemma 5.2]{BunFeeTim11}, rectangular of rank 2 or more \cite[Theorem 5.4]{BunFeeTim11}, symmetric and antisymmetric \cite[Theorem 5.5]{BunFeeTim11}. The corresponding references in \cite{BohWer11} to these results, but only for the finite dimensional cases as noted above, are:  Hilbert space \cite[Theorem 3.15]{BohWer11}, spin factor \cite[Theorem 3.5]{BohWer11}, rectangular of rank 2 or more  \cite[Theorem 3.13]{BohWer11}, symmetric and antisymmetric \cite[Theorem 3.6]{BohWer11}).

As a technical aid of possibly independent interest, the notion of a reversible $JC^*$-triple is also introduced in \cite{BunFeeTim11}. and the universal reversibility (or not) of  all Cartan factors, is also  determined.
Building on the definition of reversibility for Jordan algebras leads to the following definition.

\begin{definition} [Definition 4.1 of \cite{BunFeeTim11}]
A $JC^*$-subtriple $E$ of a TRO $T$ is said to be reversible in $T$ if
\[
a_1a_2^*a_3\cdots a_{2n}^*a_{2n+1}+a_{2n+1}a_{2n}^*\cdots a_2^*a_1\in E
\]
whenever $a_1,\ldots a_{2n+1}\in E$.   A $JC^*$-triple $E$ is universally reversible if $\alpha_E(E)$ is reversible in $T^*(E)$.
\end{definition}

\begin{theorem}
[Theorem 5.6 of \cite{BunFeeTim11}]\label{10.4}
Spin factors of dimension greater than 4 and Hilbert spaces of dimension greater than 2 are not universally reversible.  All other Cartan factors are universally reversible.
\end{theorem}

In the next theorem, operator space structures of $JC^*$-triples arising from concrete triple embeddings in $C^*$-algebras are considered, making use of universal TROs and establishing links with injective envelopes and triple envelopes, and concentrating on Cartan factors. It is necessary first to refine the definition of operator space structure.

\begin{definition} [Definition 6.1 of \cite{BunFeeTim11}]
A $JC$-operator space structure on a $JC^*$-triple $E$ is an operator space structure
determined by a linear isometry from $E$ onto a $JC^*$-subtriple of $B(H)$.
\end{definition}

Another important definition is that of an operator space ideal.  This is a norm closed ideal $\mathcal I$ of $T^*(E)$ for which ${\mathcal I}\cap \alpha_E(E)=\{0\}$.
For each operator space ideal $\mathcal I$, we have  the $JC$-operator space structure $E_{\mathcal I}$ on $E$   determined by the isometric embedding $E\mapsto T^*(E)/{\mathcal I}$ given by $x\mapsto \alpha_E(x)+{\mathcal I}$.

\begin{theorem}[Theorem 6.8 of \cite{BunFeeTim11}]\label{thm:0222121}
Let $E$ be a Cartan factor which is not a Hilbert space.   Then ${\mathcal I}\leftrightarrow E_{\mathcal I}$ is a bijective correspondence between the operator space ideals of $T^*(E)$ and the $JC$-operator space structures of $E$.
\end{theorem}

The case of a Hilbert space, excluded from Theorem~\ref{thm:0222121} is covered in \cite[Theorem 2.7]{BunTim11}.  It follows from this  that the triple envelope,
in the sense of \cite{Hamana99}, of a $JC^*$-triple $E$ which is isometric to a Cartan factor  is identified with the TRO generated by $E$.
The following theorem
shows that distinct $JC$-operator space structures on a Hilbert space cannot be completely isometric.

\begin{theorem}
[Theorem 2.6 of \cite{BunTim11}]
For a Hilbert space $E$ and operator space ideals $\mathcal I$ and $\mathcal J$ of $T^(E)$, the following are equivalent:

(a) $E_{\mathcal I}=E_{\mathcal J}$

(b) ${\mathcal I}={\mathcal J}$

(c) $E_{\mathcal I}$ and $E_{\mathcal J}$ are completely isometric.
\end{theorem}

It is also shown in \cite[Theorem 3.4]{BunTim11} how infinite-dimensional Hilbertian $JC$-operator spaces are determined explicitly by their finite-dimensional  subspaces, and that, in turn, they impose very rigid constraints upon the operator space structure of their finite-dimensional subspaces.  Finally, the operator space ideals  of the universal TRO of a Hilbert space are identified \cite[Theorem 3.7]{BunTim11}, as well as the corresponding injective envelopes \cite[Theorem 4.4]{BunTim11}.

Building on \cite{BunTim11,BunFeeTim11}  the paper \cite{BunTim13QJM}  establishes more tools for the
study of JC*-triples via the TRO and C*-algebra they generate.  
 The paper focusses on
two basic technical tools introduced in the earlier work of the authors,
namely, the universal TRO and the property of universal reversibility.
The properties developed involve on the one hand ideals and on the other
hand tensor products.

In turn, the results of \cite{BunFeeTim11,BunTim11,BunTim13QJM} are used in \cite{BunTim13JLMS} to extend Theorem~\ref{10.4} to general $JC^*$-triples (\cite[Theorem 4.7]{BunTim13JLMS}), provide the structure of universally reversible JW$^*$-triples and W$^*$-TROs  (\cite[Theorem 4.14]{BunTim13JLMS}), and show under mild conditions that a triple homomorphism on a TRO decomposes into a sum of a TRO homomorphism and a TRO antihomorphism (\cite[Proposition 5.14]{BunTim13JLMS}).

Finally we make four remarks about the paper \cite{BohWer11}.  First, the methods for computing the universal enveloping TROs of the Cartan factors (in finite dimensions only) are different from those of \cite{BunFeeTim11} in that they use grids and depend strongly on the results on grids in \cite{NeaRus03TAMS}. This approach is potentially very useful in future applications (for example, cohomology). 

Second, because the paper \cite{BohWer11} works  in the slightly more general context of $JB^*$-triples, they are able to obtain a new proof of a version of the Gelfand-Naimark theorem for $JB^*$-triples (\cite[Theorem 2]{FriRus86}.

\begin{proposition}[Corollary 2.6 of \cite{BohWer11}]
Any $JB^*$-triple $Z$ contains a unique purely exceptional ideal $J$, such that $Z/J$ is $JB^*$-triple isomorphic to a $JC^*$-triple.
\end{proposition}

Third, they define a radical for universally reversible $JC^*$-triples and use it to give the structure
of the universal enveloping TRO of a universally reversible TRO in a $C^*$-algebra with a TRO antiautomorphism of order 2, generalizing \cite[Corollary 4.5]{BunFeeTim11}. 
More generally, radicals have been defined and studied for Jordan
triple systems in, for example, \cite{Meyberg72,Loos73,Neher92}.

Finally, the results of \cite{BohWer11}  are used in \cite{BohWer11bis} to show that the classification of the symmetric spaces, derived initially from Lie theory and then from Jordan theory,  can be achieved by $K$-theoretic methods.  This in spite of the fact that $JB^*$-triples do not, in general, behave well under the formation of tensor products.

\bibliographystyle{amsplain}

\end{document}